\documentclass[12pt]{article}

\usepackage{amsmath,amsfonts,amssymb,amscd,amsthm}
\usepackage[all]{xy}
\usepackage{amsmath}

\title{Cohomological Hall algebra, exponential Hodge structures and motivic Donaldson-Thomas invariants}
\author{Maxim Kontsevich \and Yan Soibelman}

\newcommand{\op}[1]{\operatorname{#1}}
\newcommand{\C}{{\mathbb C}}
\newcommand{\R}{{\mathbb R}}

\newcommand{\Z}{{\mathbb Z}}
\newcommand{\Q}{{\mathbb Q}}

\newcommand{\M}{{\mathsf{M}}}
\newcommand{\G}{{\mathsf{G}}}
\newcommand{\Gr}{{\mathsf{Gr}}}

\newcommand{\univ}{{\rm{univ}}}

\newcommand{\B}{{\mathrm{B}}}
\renewcommand{\H}{{\mathcal{H}}}
\newcommand{\T}{{\mathsf{T}}}

\newcommand{\E}{{\mathcal E}}

\newcommand{\AC}{{\mathbb A}_\C^1}

\newcommand{\D}{{\mathcal D}}

\newcommand{\s}{{\op{sp}}}

\newcommand{\epi}{\twoheadrightarrow}
\newcommand{\mono}{\hookrightarrow}

\newcommand{\kk}{\mathbf{k}}

\newtheorem{thm}{Theorem}
\newtheorem{lmm}{Lemma}
\newtheorem{conj}{Conjecture}
\newtheorem{dfn}{Definition}
\newtheorem{cor}{Corollary}
\newtheorem*{question}{Question}
\newtheorem{prp}{Proposition}

\newtheorem{rmk}{Remark}

\begin{document}

\maketitle

\tableofcontents

\section{Introduction}

There is an old proposal in String Theory (see \cite{HarveyMoore}) which says that with a certain class of $4$-dimensional quantum theories with $N=2$ spacetime supersymmetry one should be able to associate
an algebra graded by the charge lattice, called the algebra of BPS states.
The goal of this paper is to propose a rigorous mathematical definition of an associative algebra which is presumably related to
the algebra of BPS states. We call it {\it Cohomological Hall algebra} (which we will sometimes abbreviate as COHA). Our construction can be applied to a  wide class of situations, including representations of quivers (or more generally,  smooth algebras) and  arbitrary finitely  presented algebras. 
It is also related to   multi-matrix integrals. More general framework for our theory would be the class of ind-constructible $3$-dimensional Calabi-Yau categories (see \cite{KS}) whose objects form an ind-Artin stack.

 Our definition is similar in the spirit to the conventional definition of the
 Hall algebra of an abelian category over a finite field. The main difference is that we use
 cohomology of moduli stacks of objects instead of constructible functions. In particular, Cohomological Hall algebra is different from the motivic Hall algebra introduced in \cite{KS}.
In  present paper  virtual fundamental classes
 of  proper (maybe singular) varieties classifying subrepresentations of a given representation serve as structure constants.

In the case of representations of a quiver  our construction is similar to the  one of G.~Lusztig (see \cite{Lusztig}) of the canonical basis of the   negative part of quantum enveloping algebra (see also a generalization to the affine case by I.~Grojnowski \cite{Grojnowski}). An analogous
 idea appeared long ago in the  work of H.~Nakajima (see \cite{Nak}), although in a slightly different context. That idea was used in \cite{HarveyMoore} (see also \cite{NekLosSh}) in an (unsuccessful) attempt to define the associative algebra of BPS states.
We remark that Lusztig's construction uses certain categories
of equivariant $\D$-modules on representations spaces of quivers (or rather corresponding perverse sheaves). Our approach is based morally on the category of coherent sheaves, and it is  more direct and elementary. It would be interesting to clarify the relation between quantum algebras and COHA.

In this paper we present two different versions of Cohomological Hall algebra. The most simple one is the ``off-shell'' version, related to the so-called rapid decay cohomology 
of an algebraic variety endowed with function (``potential"). We develop a generalization of the theory of mixed Hodge structures in this ``exponential"\footnote{Informally speaking, we study integrals of the form $\int \exp(f)$  where $f$ is the potential.} setting.
An advantage of rapid decay cohomology is that one can use  elementary topological methods
and obtain very strong factorization properties for the motivic Donaldson-Thomas series of the corresponding COHA\footnote{Motivic DT-series are the generating series for the Serre polynomials of graded components of COHA.}.  As a corollary, we define integer numbers which (hypothetically) count   BPS states, and satisfy wall-crossing formulas from \cite{KS}.

Another goal of the paper is to give an alternative approach to the theory of motivic Donaldson-Thomas invariants introduced in \cite{KS}. Having a quiver with potential one can define an ind-constructible $3CY$  category. Its $t$-structure has the heart consisting of finite-dimensional representations of the quiver which are critical points of the potential. This is the framework in which the theory from \cite{KS} can be applied.
In order to compare the invariants from the loc. cit. with those introduced in present paper one needs an ``on-shell'' version of Cohomological Hall algebra defined 
in terms of sheaves of vanishing cycles. The corresponding cohomology theory is related to  asymptotic expansions of exponential integrals in the formal neighborhood of the 
critical locus of the potential. We call it the {\it critical cohomology}. We prove in Section 7 that the motivic DT-series for the ``critical COHA'' essentially coincides with the one introduced in \cite{KS}. 
 In order to define critical COHA  we need several foundational results concerning sheaves of vanishing cycles. In particular we need a version of Thom-Sebastiani theorem for mixed Hodge modules
 (proven by M.~Saito in \cite{Saito-TS}) as well as the integral identity from \cite{KS}, which we prove here in the sheaf-theoretic framework.

As we pointed out above, we expect that the approach of this paper can be generalized to a wider class of $3$-dimensional Calabi-Yau categories (after a choice of generator any such category can be in a sense described ``locally" by a quiver with potential). Restriction to the case of a smooth algebra with potential allows us to avoid a discussion of some technical difficulties of the general definition from \cite{KS}. In particular there is no problem with orientation data. Furthermore,
we can prove the integrality property of generating series which is in a sense  stronger than the one conjectured in \cite{KS}. The new integrality property has a nice algebraic meaning for our
Cohomological Hall algebra  as a kind of Poincar\'e-Birkhoff-Witt theorem. The proof is based on the consideration of equivariant cohomology with respect to the action of maximal tori.

Summarizing, in all cases we define a graded vector space $\mathcal{H}$ which is:
\begin{itemize}
\item[(i)] An associative algebra (more precisely,  with a twisted associativity in the motivic setting). The algebra itself depends on a quiver with potential, not on a stability function.
\item[(ii)] For given stability function, the  motivic DT-series of the algebra admits a factorization (in the clockwise order) into a product of factors over all possible slopes. Each factor is defined in terms of the moduli space of semistable objects with a given slope. This factorization property is a consequence of an existence of a spectral sequence converging to the graded components of $\mathcal{H}$.
\item[(iii)] The motivic  DT-series admits an infinite product expression in terms of standard $q$-special functions (basically, shifted quantum dilogarithms) with integer exponents (the latter should be called  ``refined  DT-invariants" by analogy with ``refined BPS invariants'' in physics).
\end{itemize}
This product structure is based on the theory of factorization systems introduced in the paper.
The reader should keep in mind all the goals (i)-(iii) while reading about various generalizations and ramifications, which we discuss in the paper.

Finally, we remark that it would be  highly desirable to relate the classical limits of our motivic Donaldson-Thomas invariants with the invariants introduced by Joyce and Song (see \cite{JoyceSong}). More precisely, we expect that the category of coherent sheaves used in the loc. cit. carries orientation data and our motivic DT-invariants (see Section 6)  give in the classical limit the invariants introduced there.

Let us briefly discuss the content of the paper.

Section 2 is devoted to the definition of COHA in the case of a quiver without potential. We also present the  formula for the product based on the toric 
localization technique and apply it to several  examples. We observe that the canonical grading of COHA  by the ``charge lattice" can be extended to a grading by the Heisenberg group
 which is non-commutative.

Section 3 is devoted to various generalizations of COHA. In particular, the path algebra of a quiver with the set of vertices $I$ can be replaced by an $I$-bigraded algebra, which is smooth in the sense of Cuntz and Quillen. We also remark that there should be a natural generalization of COHA as an $A_{\infty}$-algebra in the dg-enhancement of the Voevodsky category of mixed motives.

In Section 4 we discuss COHA in the case of smooth $I$-bigraded algebra with potential. The definition is based on the properties of the category $EMHS$ of exponential mixed Hodge structures introduced in that section. Objects of $EMHS$  admit weight filtration which allows us to introduce Serre polynomials of graded components of COHA. There are several realizations of COHA depending on a choice of cohomology theory (i.e. a cohomology functor from the tensor category $EMHS$ to a Tannakian category). For example, one can choose rapid decay cohomology theory, which corresponds to the Betti realization of $EMHS$.   Another choice is the de Rham realization. The comparison between Betti and de Rham realizations is related to the theory of matrix integrals, since the conventional pairing between algebraic differential forms and algebraic cycles (period map) now becomes a pairing between ``exponential forms'' of the type $exp(f)\alpha$ and non-compact cycles $C$ such that the restriction $Re(f)_{|C}$ is bounded.
We end  Section 4 with a result which sometimes allows us to reduce computations of COHA from (cohomological) dimension 3 to dimension 2.

In Section 5 we introduce the notion of stability condition and discuss the corresponding spectral sequence converging to COHA. We also introduce motivic DT-series of a smooth $I$-bigraded algebra with potential and prove the Factorization Formula (a.k.a. wall-crossing formula) similar to the one from \cite{KS}. In the case of a quiver with potential we study how the DT-series changes with respect to a mutation (the result agrees with the one from Section 8 of \cite{KS}).

In Section 6 we introduce motivic DT-invariants and discuss the product structure of motivic DT-series.
Integrality of the exponents in the product formula (refined DT-invariants), see Section 6.1, is a non-trivial fact which follows from the so-called admissibility property of the motivic DT-series. Corresponding technique of factorization systems and proofs are contained in Section 6. Although COHA and motivic DT-series do not depend on the stability function (a.k.a. central charge), the motivic DT-invariants do.

Section 7 is devoted to the critical COHA defined in terms of vanishing cycles functor. The corresponding cohomology theory is called critical 
cohomology. Some important properties of cohomology theories used in the previous sections do not hold for critical cohomology, e.g. Thom isomorphism. As a result, the very definition of COHA requires much more work. The situation is in a sense similar to the categorical one from \cite{KS}. In particular, we need to prove in Section 7.8 the integral identity sketched in the loc. cit. Importance of the critical COHA is explained in Section 7.10, where we prove that the corresponding motivic DT-series is basically the same as the ``categorical one'' from \cite{KS}.

In the last Section 8 we discuss two speculative applications: categorification of COHA and motivic DT-invariants arising in the Chern-Simons theory with complex gauge group. There are many more applications and speculations which we decided not to include in order to save the  space. For example, one can give an interpretation of line operators from \cite{GMN-3} in terms of certain representations of COHA.

{\bf Acknowledgments}. We thank to Pierre Deligne, Michael Douglas, Davide Gaiotto, Ian Grojnowski, Greg Moore, Nikita Nekrasov, Andrei Okounkov,
Joerg Sch\"urmann for useful conversations and communications. Also we are grateful to referees for several remarks and corrections. We thank to Morihiko Saito for sending us a preliminary version of \cite{Saito-TS}. This work was started during the workshop at Simons Center for Geometry and Physics in January, 2008. Y.S. thanks to the University of Miami and IHES for the hospitality. His work was partially supported by an NSF grant.

\section{Cohomological Hall algebra of a quiver without potential}

 \subsection{Stacks of representations and their cohomology}

 Let us fix a finite quiver $Q$, with the set $I$ of vertices, and $a_{ij}\in \Z_{\geqslant 0}$ arrows from $i$ to $j$
 for $i,j\in I$.
 For any dimension vector
$$\gamma=(\gamma^i)_{i\in I} \in \Z_{\geqslant 0}^I$$
consider the space of representations of $Q$ in complex coordinate
 vector spaces of dimensions $(\gamma^i)_{i\in I}$
$$\M_\gamma=\M_{\gamma}^Q\simeq \prod_{i,j\in I} \C^{a_{ij}\gamma^i \gamma^j}$$
endowed with the action by conjugation of a complex algebraic group
$$\G_\gamma:=\prod_{i\in I} \op{GL}(\gamma^i,\mathbb{C})\,\,.$$
The quotient stack $\M_\gamma/\G_\gamma$ is the stack of representations of $Q$ with dimension vector $\gamma$.
We will consider the cohomology of this stack, i.e. the equivariant
cohomology of $\M_\gamma$ with $\G_\gamma$-action.  We use the standard model
	$$\op{Gr}(d,\C^\infty):=\varinjlim\op{Gr}(d,\C^N)\,,\,\,N\to +\infty$$ of
the classifying space of $\op{GL}(d,\C)$ for $d\geqslant 0$, and define
  $$\B\G_\gamma:=\prod_{i\in I} \op{BGL}(\gamma^i,\C)=\prod_{i\in I} \op{Gr}(\gamma^i,\C^\infty)\,\,.$$
  Stack $\M_\gamma/\G_\gamma$ gives the universal family over  $\rm{B}\G_\gamma$
  $$\M_\gamma^{\univ}:=\left(\mathrm{E}\G_\gamma\times \M_\gamma\right)/{\G_\gamma}\,\,, $$ where
$\rm{E}\G_\gamma\to \B\G_\gamma$ is the standard universal
  $\G_\gamma$-bundle.

  We introduce a $\Z_{\geqslant 0}^I$-graded abelian group
  $$\mathcal{H}:=\oplus_\gamma \mathcal{H}_\gamma\,\,,$$
  where each component is defined as an equivariant cohomology
  $$\mathcal{H}_\gamma:=H^\bullet_{\G_\gamma} (\M_\gamma):=
  H^\bullet(\M_\gamma^{\univ})=\oplus_{n\geqslant 0} H^n(\M_\gamma^{\univ})\,\,.$$
Here by cohomology of a complex algebraic variety (or of an inductive limit of varieties) we mean the usual (Betti)
cohomology with coefficients in $\Z$.

Notice that the natural Hodge structure on $H^n(\M_\gamma^{\univ})$ is pure of weight $n$,
 i.e. incidentally the cohomological degree coincides with the weight.

\subsection{Multiplication}
Fix any $\gamma_1,\gamma_2\in \Z_{\geqslant 0}^I$ and denote $\gamma:=\gamma_1+\gamma_2$. Denote by ${\M}_{\gamma_1,\gamma_2}$
 the space of representations of $Q$ in
  coordinate spaces of dimensions $(\gamma_1^i+\gamma_2^i)_{i\in I}$
  such that the standard coordinate subspaces of dimensions $(\gamma_1^i)_{i\in I}$  form a subrepresentation. Obviously ${\M}_{\gamma_1,\gamma_2}$ is an affine space, and also a  closed subspace of $\M_\gamma$. The group $\G_{\gamma_1,\gamma_2}\subset \G_\gamma$ consisting of transformations
 preserving subspaces $\left(\C^{\gamma_1^i}\subset \C^{\gamma^i}\right)_{i\in I}$ (i.e. the group of block upper-triangular matrices), acts on ${\M}_{\gamma_1,\gamma_2}$. In what follows we will
use the model of $\B\G_{\gamma_1,\gamma_2}$ which is the total space of the bundle over
 $\B\G_\gamma$ with the fiber $\G_\gamma/\G_{\gamma_1,\gamma_2 }$.

 Let us consider a morphism
 $$\mathsf{m}_{\gamma_1,\gamma_2}:\mathcal{H}_{\gamma_1}\otimes \mathcal{H}_{\gamma_2}\to \mathcal{H}_{\gamma}=\mathcal{H}_{\gamma_1+\gamma_2},$$
 which is the composition of  the multiplication morphism (which becomes K\"unneth isomorphism after the extension of coefficients for cohomology from $\Z$ to $\Q$)
$$\otimes:H^\bullet_{\G_{\gamma_1}}(\M_{\gamma_1})\otimes H^\bullet_{\G_{\gamma_2}}(\M_{\gamma_2})
\to H^\bullet_{\G_{\gamma_1}\times\G_{\gamma_2 }}(\M_{\gamma_1}\times \M_{\gamma_2})\,,$$
and of the following sequence of 3 morphisms:
$$H^\bullet_{\G_{\gamma_1}\times \G_{\gamma_2}}({\M}_{\gamma_1}\times {\M}_{\gamma_2} )
\stackrel{\simeq}{\to} H^\bullet_{\G_{\gamma_1,\gamma_2 } } ({\M}_{\gamma_1,\gamma_2} )
\to H^{\bullet+2c_1 }_{\G_{\gamma_1,\gamma_2 } } (\M_\gamma)\to H^{\bullet+2c_1+2c_2}_{\G_{\gamma } } (\M_\gamma)\,,$$
where
\begin{enumerate}
\item the first arrow is an isomorphism induced by natural projections
of spaces and groups, inducing homotopy equivalences
$$  \M_{\gamma_1}\times \M_{\gamma_2}\stackrel{\sim}{\twoheadleftarrow}
 \M_{\gamma_1,\gamma_2 }
		\,,\,\,\,\, \G_{\gamma_1}\times \G_{\gamma_2} \stackrel{\sim}{\twoheadleftarrow} \G_{\gamma_1,\gamma_2 }
\,\,,  $$
\item the second arrow is the pushforward map associated with the closed
$\G_{\gamma_1,\gamma_2 } $-equivariant embedding
${\M}_{\gamma_1,\gamma_2}\mono {\M}_{\gamma}$ of complex manifolds,
\item the third arrow is the pushforward map associated with the fundamental class
of compact complex manifold $\G_\gamma/\G_{\gamma_1,\gamma_2 }$, which is the product of Grassmannians $\prod_{i\in I}\op{Gr}(\gamma^i_1,\C^{\gamma^i})$.
\end{enumerate}

Shifts in the cohomological degrees are given by
$$c_1=\dim_\C \M_\gamma-\dim_\C {\M}_{\gamma_1,\gamma_2}\,,\,\,c_2=-\op{dim}_\C\G_\gamma/\G_{\gamma_1,\gamma_2 }\,\,.$$

We endow  $\mathcal{H}$ with a product $\mathsf{m}:\mathcal{H}\otimes \mathcal{H} \to \mathcal{H},\,\,\,\mathsf{m}:=\sum_{\gamma_1,\gamma_2}\mathsf{m}_{\gamma_1,\gamma_2}$.

\begin{thm} The product $\mathsf{m}$ on $\mathcal{H}$ is associative.
\end{thm}

Moreover, it is straightforward to see that the element $1_{\M_0}\in \mathcal{H}_0$ is the unit for the product. The proof of Theorem 1
will be given
in Section 2.3.

\begin{dfn} The associative unital $\Z^I_{\geqslant 0}$-graded
 algebra $\mathcal{H}$ with the product $\mathsf{m}$  is called the {\bf Cohomological Hall algebra} associated with the quiver $Q$.
\end{dfn}

The multiplication does not preserve the cohomological grading. The shift is given by
$$2(c_1+c_2)= -2\chi_Q(\gamma_1,\gamma_2)\,,$$
where
$$\chi_Q(\gamma_1,\gamma_2):=-\sum_{i,j\in I} a_{ij}\gamma_1^j \gamma_2^i+\sum_{i\in I}\gamma^i_1\gamma^i_2$$
is the Euler form
 on the $K_0$ group of the  category of finite-dimensional representations of $Q$:
 $$\chi_Q(\gamma_1,\gamma_2)=\dim{\op{Hom}}(E_1,E_2)-\dim{\op{Ext}^1}(E_1,E_2)=
\chi\left({\op{Ext}}^\bullet(E_1,E_2)\right)\,\,, $$
where $E_1,E_2$ are arbitrary representations of $Q$, of dimension vectors $\gamma_1,\gamma_2$.

The multiplication in $\mathcal{H}$ can be defined in a slightly different way, which is maybe   more intuitively clear. 
Namely, consider the following  map of manifolds endowed with $\G_\gamma$-action:
   $$\pi:\Gr_{\gamma_1,\gamma}:=\G_\gamma\times_{\G_{\gamma_1,\gamma_2}}\M_{\gamma_1,\gamma_2}\to \M_\gamma\,,\,\,\,\,(g,m)\mapsto gm\,\,.$$
Notice that $\dim \Gr_{\gamma_1,\gamma}- \dim \M_\gamma= \chi_Q(\gamma_1,\gamma_2)$. 
The map $\pi$ is proper, hence it induces the pushforward (Gysin) morphism:
   $$\pi_*:H^\bullet_{\G_\gamma}\left( \Gr_{{\gamma_1,\gamma}} \right)
   \to H^{\bullet-2\chi_Q(\gamma_1,\gamma_2)  }_{\G_\gamma}\left( \M_\gamma \right)\,\,.$$
   There are two natural $\G_\gamma$-equivariant bundles of representations of $Q$ of dimension vectors $\gamma_1$ and $\gamma_2$ on $\Gr_{\gamma_1,\gamma}$. Combining the pullback
   morphism
   $$H^\bullet_{\G_{\gamma_1}}(\M_{\gamma_1})\otimes H^\bullet_{\G_{\gamma_2}}(\M_{\gamma_2})\to H^\bullet_{\G_\gamma}(\Gr_{\gamma_1,\gamma}   )$$
   with $\pi_*$, we obtain an equivalent definition of the multiplication morphism.

One can make still another reformulation using language of stacks. The natural morphism of  stacks
$$\M_{\gamma_1,\gamma_2}/\G_{\gamma_1,\gamma_2}\to \M_\gamma/\G_\gamma$$
is a {\it proper} morphism of {\it smooth} Artin stacks, hence it induces the {\it pushforward} map
on cohomology. Combining it with the pullback by the homotopy equivalence
 $  \M_{\gamma_1}/\G_{\gamma_1}\times
 \M_{\gamma_2}/\G_{\gamma_2}   \leftarrow
 \M_{\gamma_1,\gamma_2}/\G_{\gamma_1,\gamma_2}$, we obtain  $\mathsf{m}_{\gamma_1,\gamma_2}$.

\subsection{Proof of associativity}

For given $\gamma_1,\gamma_2,\gamma_3\in \Z_{\geqslant 0}^I$,
consider the following commutative diagram where we omit  for the convenience shifts
 in the cohomological degree:
\[\xymatrix{
 H^\bullet_{\begin{picture}(11,13)
\put(3,8){\circle*{2}}
                   \put(6,5){\circle*{2}}
                                            \put(9,2){\circle*{2}}
\end{picture}}\bigl(\M_{\begin{picture}(11,13)
\put(3,8){\circle*{2}}
                   \put(6,5){\circle*{2}}
                                            \put(9,2){\circle*{2}}
\end{picture}}\bigr)  \ar[d]^-{\simeq}  \ar[r]^-{\simeq}
  & H^\bullet_{\begin{picture}(11,13)
\put(3,8){\circle*{2}}
                      \put(6,5){\circle*{2}}\put(9,5){\circle*{2}}
                                            \put(9,2){\circle*{2}}                                \end{picture}}\bigl(\M_{\begin{picture}(11,13)
\put(3,8){\circle*{2}}
                      \put(6,5){\circle*{2}}\put(9,5){\circle*{2}}
                                            \put(9,2){\circle*{2}}                                \end{picture}}\bigr)\ar[d]^-{\simeq}\ar[r]&
H^\bullet_{\begin{picture}(11,13)
\put(3,8){\circle*{2}}
                      \put(6,5){\circle*{2}}\put(9,5){\circle*{2}}
                                            \put(9,2){\circle*{2}}                                 \end{picture}}\bigl(\M_{\begin{picture}(11,13)
\put(3,8){\circle*{2}}
                      \put(6,5){\circle*{2}}\put(9,5){\circle*{2}}
                      \put(6,2){\circle*{2}}\put(9,2){\circle*{2}}                                 \end{picture}}\bigr)\ar[d]^-{\simeq}\ar[r]&
H^\bullet_{\begin{picture}(11,13)
\put(3,8){\circle*{2}}
                      \put(6,5){\circle*{2}}\put(9,5){\circle*{2}}
                      \put(6,2){\circle*{2}}\put(9,2){\circle*{2}}                                 \end{picture}}\bigl(\M_{\begin{picture}(11,13)
\put(3,8){\circle*{2}}
                      \put(6,5){\circle*{2}}\put(9,5){\circle*{2}}
                      \put(6,2){\circle*{2}}\put(9,2){\circle*{2}}                                 \end{picture}}\bigr)\ar[d]^-{\simeq}\\
  H^\bullet_{\begin{picture}(11,13)
\put(3,8){\circle*{2}}\put(6,8){\circle*{2}}
                      \put(6,5){\circle*{2}}
                                            \put(9,2){\circle*{2}}
\end{picture}}\bigl(\M_{\begin{picture}(11,13)
\put(3,8){\circle*{2}}\put(6,8){\circle*{2}}
                      \put(6,5){\circle*{2}}
                                            \put(9,2){\circle*{2}}
\end{picture}}\bigr) \ar[d]\ar[r]^-{\simeq} &
     H^\bullet_{\begin{picture}(11,13)
\put(3,8){\circle*{2}}\put(6,8){\circle*{2}}\put(9,8){\circle*{2}}
                      \put(6,5){\circle*{2}}\put(9,5){\circle*{2}}
                                            \put(9,2){\circle*{2}}
\end{picture}}\bigl(\M_{\begin{picture}(11,13)
\put(3,8){\circle*{2}}\put(6,8){\circle*{2}}\put(9,8){\circle*{2}}
                      \put(6,5){\circle*{2}}\put(9,5){\circle*{2}}
                                            \put(9,2){\circle*{2}}
\end{picture}}\bigr)
   \ar[d]\ar[r] &
   H^\bullet_{\begin{picture}(11,13)
\put(3,8){\circle*{2}}\put(6,8){\circle*{2}}\put(9,8){\circle*{2}}
                      \put(6,5){\circle*{2}}\put(9,5){\circle*{2}}
                                            \put(9,2){\circle*{2}}                                 \end{picture}}\bigl(\M_{\begin{picture}(11,13)
\put(3,8){\circle*{2}}\put(6,8){\circle*{2}}\put(9,8){\circle*{2}}
                      \put(6,5){\circle*{2}}\put(9,5){\circle*{2}}
                      \put(6,2){\circle*{2}}\put(9,2){\circle*{2}}                                 \end{picture}}\bigr)\ar[r]\ar[d]&
H^\bullet_{\begin{picture}(11,13)
\put(3,8){\circle*{2}}\put(6,8){\circle*{2}}\put(9,8){\circle*{2}}
                      \put(6,5){\circle*{2}}\put(9,5){\circle*{2}}
                      \put(6,2){\circle*{2}}\put(9,2){\circle*{2}}                                 \end{picture}}\bigl(\M_{\begin{picture}(11,13)
\put(3,8){\circle*{2}}\put(6,8){\circle*{2}}\put(9,8){\circle*{2}}
                      \put(6,5){\circle*{2}}\put(9,5){\circle*{2}}
                      \put(6,2){\circle*{2}}\put(9,2){\circle*{2}}                                 \end{picture}}\bigr)\ar[d]\\
 H^\bullet_{\begin{picture}(11,13)
\put(3,8){\circle*{2}}\put(6,8){\circle*{2}}
                      \put(6,5){\circle*{2}}
                                            \put(9,2){\circle*{2}}
\end{picture}}\bigl(\M_{\begin{picture}(11,13)
\put(3,8){\circle*{2}}\put(6,8){\circle*{2}}
\put(3,5){\circle*{2}}\put(6,5){\circle*{2}}
                                            \put(9,2){\circle*{2}}
\end{picture}}\bigr) \ar[d]\ar[r]^-{\simeq} &
\ar[d]
   H^\bullet_{\begin{picture}(11,13)
\put(3,8){\circle*{2}}\put(6,8){\circle*{2}}\put(9,8){\circle*{2}}
                      \put(6,5){\circle*{2}}\put(9,5){\circle*{2}}
                                            \put(9,2){\circle*{2}}
\end{picture}}\bigl(\M_{\begin{picture}(11,13)
\put(3,8){\circle*{2}}\put(6,8){\circle*{2}}\put(9,8){\circle*{2}}
\put(3,5){\circle*{2}}\put(6,5){\circle*{2}}\put(9,5){\circle*{2}}
                                            \put(9,2){\circle*{2}}
\end{picture}}\bigr)
   \ar[d] \ar[r]
   &
   H^\bullet_{\begin{picture}(11,13)
\put(3,8){\circle*{2}}\put(6,8){\circle*{2}}\put(9,8){\circle*{2}}
                      \put(6,5){\circle*{2}}\put(9,5){\circle*{2}}
                                            \put(9,2){\circle*{2}}
\end{picture}}\bigl(\M_{\begin{picture}(11,13)
\put(3,8){\circle*{2}}\put(6,8){\circle*{2}}\put(9,8){\circle*{2}}
\put(3,5){\circle*{2}}\put(6,5){\circle*{2}}\put(9,5){\circle*{2}}
\put(3,2){\circle*{2}}\put(6,2){\circle*{2}}                                 \put(9,2){\circle*{2}}
\end{picture}}\bigr)
    \ar[r] \ar[d] &
    H^\bullet_{\begin{picture}(11,13)
\put(3,8){\circle*{2}}\put(6,8){\circle*{2}}\put(9,8){\circle*{2}}
                      \put(6,5){\circle*{2}}\put(9,5){\circle*{2}}
                      \put(6,2){\circle*{2}}\put(9,2){\circle*{2}}                                 \end{picture}}\bigl(\M_{\begin{picture}(11,13)
\put(3,8){\circle*{2}}\put(6,8){\circle*{2}}\put(9,8){\circle*{2}}
\put(3,5){\circle*{2}}\put(6,5){\circle*{2}}\put(9,5){\circle*{2}}
\put(3,2){\circle*{2}}\put(6,2){\circle*{2}}\put(9,2){\circle*{2}}                                 \end{picture}}\bigr)\ar[d]\\
    H^\bullet_{\begin{picture}(11,13)
\put(3,8){\circle*{2}}\put(6,8){\circle*{2}}
\put(3,5){\circle*{2}}\put(6,5){\circle*{2}}
                                            \put(9,2){\circle*{2}}
\end{picture}}\bigl(\M_{\begin{picture}(11,13)
\put(3,8){\circle*{2}}\put(6,8){\circle*{2}}
\put(3,5){\circle*{2}}\put(6,5){\circle*{2}}
                                            \put(9,2){\circle*{2}}
\end{picture}}\bigr)\ar[r]^-{\simeq}
   &
   H^\bullet_{\begin{picture}(11,13)
\put(3,8){\circle*{2}}\put(6,8){\circle*{2}}\put(9,8){\circle*{2}}
\put(3,5){\circle*{2}}\put(6,5){\circle*{2}}\put(9,5){\circle*{2}}
                                            \put(9,2){\circle*{2}}
\end{picture}}\bigl(\M_{\begin{picture}(11,13)
\put(3,8){\circle*{2}}\put(6,8){\circle*{2}}\put(9,8){\circle*{2}}
\put(3,5){\circle*{2}}\put(6,5){\circle*{2}}\put(9,5){\circle*{2}}
                                            \put(9,2){\circle*{2}}
\end{picture}}\bigr)\ar[r]&
H^\bullet_{\begin{picture}(11,13)
\put(3,8){\circle*{2}}\put(6,8){\circle*{2}}\put(9,8){\circle*{2}}
\put(3,5){\circle*{2}}\put(6,5){\circle*{2}}\put(9,5){\circle*{2}}
                                            \put(9,2){\circle*{2}}
\end{picture}}\bigl(\M_{\begin{picture}(11,13)
\put(3,8){\circle*{2}}\put(6,8){\circle*{2}}\put(9,8){\circle*{2}}
\put(3,5){\circle*{2}}\put(6,5){\circle*{2}}\put(9,5){\circle*{2}}
\put(3,2){\circle*{2}}\put(6,2){\circle*{2}}                                 \put(9,2){\circle*{2}}
\end{picture}}\bigr)\ar[r]
&
    H^\bullet_{\begin{picture}(11,13)
\put(3,8){\circle*{2}}\put(6,8){\circle*{2}}\put(9,8){\circle*{2}}
\put(3,5){\circle*{2}}\put(6,5){\circle*{2}}\put(9,5){\circle*{2}}
\put(3,2){\circle*{2}}\put(6,2){\circle*{2}}                                 \put(9,2){\circle*{2}}
\end{picture}}\bigl(\M_{\begin{picture}(11,13)
\put(3,8){\circle*{2}}\put(6,8){\circle*{2}}\put(9,8){\circle*{2}}
\put(3,5){\circle*{2}}\put(6,5){\circle*{2}}\put(9,5){\circle*{2}}
\put(3,2){\circle*{2}}\put(6,2){\circle*{2}}                                 \put(9,2){\circle*{2}}
\end{picture}}\bigr)
}\]
Here we use a shorthand notation for various equivariant cohomology groups, with dots denoting non-trivial blocks for operators in $\oplus_i (\C^{\gamma_1^i+\gamma_2^i+\gamma_3^i})$, e.g.
$$
 H^\bullet_{\begin{picture}(11,13)
\put(3,8){\circle*{2}}\put(6,8){\circle*{2}}
                      \put(6,5){\circle*{2}}
                                            \put(9,2){\circle*{2}}
\end{picture}}\bigl(\M_{\begin{picture}(11,13)
\put(3,8){\circle*{2}}\put(6,8){\circle*{2}}
\put(3,5){\circle*{2}}\put(6,5){\circle*{2}}
                                            \put(9,2){\circle*{2}}
\end{picture}}\bigr)=H^\bullet_{\G_{\gamma_1,\gamma_2}\times \G_{\gamma_3}}(\M_{\gamma_1+\gamma_2}\times\M_{\gamma_3}) \,\,.      $$
Traveling first down then right, we obtain the map
$a_1\otimes a_2 \otimes a_3\mapsto (a_1 \cdot a_2)\cdot a_3$,
whereas traveling first right then down we get $a_1\otimes a_2 \otimes a_3\mapsto  a_1\cdot(a_2\cdot a_3)$.
The associativity of the product on $\mathcal{H}$ is proven. $\blacksquare$

         \subsection{Explicit formula for the product}

Here we calculate the product on  $\mathcal{H}$ using the toric localization formula.

First of all, for any $\gamma$ the  abelian group $\mathcal{H}_\gamma$ is the cohomology of the classifying space $\B\G_\gamma$, as the manifold
 $\M_\gamma$ is contractible. It is well-known that $H^\bullet({\rm{BGL}}(n,\C))$ can be canonically identified with the algebra of symmetric
 polynomials with integer coefficients in $n$ variables of degree $+2$ for any $n\geqslant 0$, via the embedding
$$H^\bullet({\mathrm{BGL}}(n,\C))\mono H^\bullet(\B(\C^\times)^n)\simeq \Z[x_1,\dots,x_n]\,\,,$$
induced by the diagonal embedding $(\C^\times)^n\mono \mathrm{GL}(n,\C)$.
Therefore, $\mathcal{H}_\gamma$ is realized as the abelian group of polynomials in variables $(x_{i,\alpha})_{i\in I, \alpha\in\{1,\dots,\gamma^i\}}$ symmetric under the group $\prod_{i\in I} {\mathrm{Sym}}_{\gamma^i}$ of permutations preserving index $i$ and
permuting index $\alpha$.

\begin{thm}
The product $f_1\cdot f_2$ of elements $f_i\in \mathcal{H}_{\gamma_i},\,i=1,2$  is given by the symmetric function
$g((x_{i,\alpha})_{i\in I, \alpha\in \{1,\dots,\gamma^i\}})$, where $\gamma:=\gamma_1+\gamma_2$,  obtained from
the following function in variables $({x}'_{i,\alpha})_{i\in I, \alpha\in \{1,\dots,\gamma^i_1\}}$ and
  $({{x}}''_{i,\alpha})_{i\in I, \alpha\in \{1,\dots,\gamma^i_2\}}$:
$$f_1(({x}'_{i,\alpha}))\,
f_2(({x}''_{i,\alpha})) \,\,
\frac{ \prod_{i,j\in I}\prod_{\alpha_1=1}^{\gamma^i_1}\prod_{\alpha_2=1}^{\gamma^j_2} ({{x}}''_{j,\alpha_2}
-{x}'_{i,\alpha_1}  )^{a_{ij}
  }}{\prod_{i\in I}
\prod_{\alpha_1=1}^{\gamma^i_1}\prod_{\alpha_2=1}^{\gamma^i_2}({{x}}''_{i,\alpha_2}-{x}'_{i,\alpha_1})
}\,\,,$$
by taking the sum over all shuffles for any given $i\in I$ of the variables
   ${x}'_{i,\alpha},{{x}}''_{i,\alpha}$ (the sum is over $\prod_{i\in I}\binom{\gamma^i}{\gamma_1^i}$ shuffles).
\end{thm}
 {\it Proof.}
The pushforward map
 $$  H^{\bullet}(B\G_{\gamma_1}\times B\G_{\gamma_2})\simeq
 H^\bullet_{\G_{\gamma_1,\gamma_2 } } ({\M}_{\gamma_1,\gamma_2} )
\to H^{\bullet+2c_1 }_{\G_{\gamma_1,\gamma_2 } } (\M_\gamma)\simeq  H^{\bullet+2 c_1}(B\G_{\gamma_1}\times B\G_{\gamma_2})  $$
 is just the multiplication by the equivariant Euler class
 $$e_{\gamma_1,\gamma_2}\in H^{2 c_1}(B\G_{\gamma_1}\times B\G_{\gamma_2})\,,\,\,\,c_1=\sum_{i,j\in I} a_{ij}\gamma_1^i\gamma_2^j$$ of the fiber at $0\in\M_{\gamma_1,\gamma_2} $ of the normal bundle to $\M_{\gamma_1,\gamma_2}\subset\M_{\gamma}$.

 The formula for the multiplication implies
that $f_1\cdot f_2$ for $f_i\in \mathcal{H}_{\gamma_i},\,i=1,2$ is obtained in the following way.
 First, we consider two canonical $\G_{\gamma_i}$-bundles over $\G_\gamma/\G_{\gamma_1,\gamma_2}$ for $i=1,2$, which are equivariant with respect to $\G_\gamma$-action.
 Then we take the product of their $\G_\gamma$-equivariant characteristic classes corresponding to symmetric polynomials
  $f_1,f_2$, multiply the result by the $\G_\gamma$-equivariant  class corresponding to $e_{\gamma_1,\gamma_2}$, and then take the
 integral over the $\G_\gamma$-equivariant fundamental class of $\G_{\gamma}/\G_{\gamma_1,\gamma_2}=\prod_{i\in I}\mathrm{Gr}(\gamma_1^i,\C^{\gamma^i})
 $. The result is an element of the cohomology ring
 of $\B\G_\gamma$.

 We can use the equivariant cohomology with respect to the maximal torus $\mathsf{T}_\gamma$ of $\G_\gamma$, instead of $\G_\gamma$. The set of fixed points
 of $\mathsf{T}_\gamma$ on $\G_\gamma/\G_{\gamma_1,\gamma_2}$ is the set of collections of coordinate subspaces, i.e. exactly
the set of shuffles of variables ${x}'_{i,\alpha}$ and ${{x}}''_{i,\alpha}$.
  The numerator in the formula in the theorem is the product of
 $\mathsf{T}_\gamma$-weights corresponding to the class
 $e_{\gamma_1,\gamma_2}$.
 The denominator is the product of $\mathsf{T}_\gamma$-weights in the tangent space of  $\G_\gamma/\G_{\gamma_1,\gamma_2}$.
The classical fixed point formula of  Bott
 gives the result.
  $\blacksquare$

\begin{rmk} The above algebra is a special case of Feigin-Odesskii shuffle algebra (see \cite{Feigin}, \cite{Enriq}).

\end{rmk}

\subsection{Example: quivers with one vertex}

Let $Q=Q_d$ be now a quiver with just one vertex and $d\geqslant 0$ loops. Then the product formula from the previous section specializes
to
\begin{multline*}
(f_1\cdot f_2)(x_1,\dots,x_{n+m}):=\\
\sum_{\substack{ i_1<\dots<i_n\\ j_1<\dots<j_m\\
\{i_1,\dots,i_n,j_1,\dots,j_m\}=\\
=\{1,\dots,n+m\}}} f_1(x_{i_1},\dots,x_{i_n})\,f_2(x_{j_1},\dots,x_{j_m})\,
\left(\prod_{k=1}^n\prod_{l=1}^m(x_{j_l}-x_{i_k})\right)^{d-1}
\end{multline*}
for symmetric polynomials, where $f_1$ has $n$  variables, and  $f_2$ has $m$ variables. The product $f_1\cdot f_2$ is a symmetric polynomial
in $n+m$ variables.

We introduce a double grading on  algebra $\H$, by declaring that a homogeneous symmetric polynomial of degree $k$ in $n$ variables
 has bigrading $(n,2k+(1-d)n^2)$. Equivalently, one can shift the cohomological grading in $H^\bullet(\mathrm{BGL}(n,\C))$ by
$[(d-1)n^2]$.

It follows directly from the product formula that the bigraded algebra is commutative for odd $d$, and supercommutative for
 even $d$. The parity in this algebra is given by the parity of the shifted cohomological degree.

It is easy to see that for $d=0$ the algebra $\mathcal{H}$ is  an exterior algebra (i.e. Grassmann algebra) generated by
 odd elements $\psi_1,\psi_3,\psi_5,\dots$ of bidegrees $(1,1),(1,3),(1,5),\dots$.
 Generators $(\psi_{2i+1})_{i\geqslant 0}$ correspond to the additive generators $(x^i)_{i\geqslant 0}$ of
$$H^\bullet(\C P^\infty)=H^\bullet(\mathrm{BGL}(1,\C))\simeq \Z[x]\simeq \Z[x_1]\,.$$
 A monomial in the exterior algebra
$$\psi_{2i_1+1}\,{\cdot}\,{\dots}\, {\cdot}\, \psi_{2i_n+1}\in \mathcal{H}_{n,\sum_{k=1}^n (2i_k+1)}\,,\,\,\,0\leqslant i_1<\dots< i_n$$
corresponds to the Schur symmetric function $s_\lambda(x_1,\dots,x_n)$,
where $$\lambda=(i_n+(1-n),i_{n-1}+(2-n),\dots,i_1)$$
 is a partition of length $\leqslant n$.

Similarly, for $d=1$ algebra $\mathcal{H}$ is isomorphic (after tensoring by $\Q$) to the algebra of symmetric functions in infinitely many variables, and it is a polynomial
 algebra generated by even elements $\phi_0,\phi_2,\phi_4,\dots$ of bidegrees $(1,0),(1,2),(1,4),\dots$.
 Again, the generators $(\phi_{2i})_{i\geqslant 0}$ correspond to the additive generators $(x^i)_{i\geqslant 0}$ of $H^\bullet(\C P^\infty)\simeq \Z[x]$.
For any $n\geqslant 0$ the collection  of all
 monomials of degree $n$ in generators $(\phi_{2i})_{i\geqslant 0}$  coincides (up to non-zero factors) with the monomial basis
 of the algebra of symmetric functions in $n$ variables.

 Notice that the underlying additive group of the algebra $\mathcal{H}$ is equal to $\oplus_{n\geqslant 0} H^\bullet({\rm{BGL}}(n,\C))$) and hence does not depend on $d$. The  isomorphism (after tensoring by $\Q$) between the underlying additive groups of the free polynomial algebra ($d=1$) and of the free exterior algebra ($d=0$), is in fact a part of the well-known boson-fermion correspondence.

For general $d$, the Hilbert-Poincar\'e series $P_d=P_d(z,q^{1/2})$ of bigraded algebra $\mathcal{H}$ twisted by the sign $(-1)^{\op{parity}}$ is the generalized $q$-exponential function:
$$\sum_{n\geqslant 0, m\in \Z}(-1)^m
\dim(\mathcal{H}_{n,m})\,z^n q^{m/2}=\sum_{n\geqslant 0} \frac{(-q^{1/2})^{(1-d)n^2}}
{(1-q)\dots(1-q^n)} z^n\in \Z((q^{1/2}))[[z]]\,\,.$$
In cases $d=0$ and $d=1$ this series decomposes in an infinite product:
$$P_0=(q^{1/2}z;q)_\infty=\prod_{i\geqslant 0} (1-q^{i+1/2}z)\,,\,\,\,\,\,\,\,\,\,\,\,\, P_1=\frac{1}{
(z;q)_\infty}=\prod_{i\geqslant 0} \frac{1}{1-q^{i}z}\,\,,$$
where we use the standard notation for the $q$-Pochhammer symbol:
$$(x;q)_\infty:=(1-x)(1-qx)(1-q^2x)\dots\,\,\,\,.$$
 Our general results in Section 6 will imply the following numerical result.
\begin{thm}\label{thm:inf-product} For any $d\geqslant 0$ there exist integers $\delta^{(d)}(n,m)$ for all $n\geqslant 1$ and $m\in (d-1)n+2\Z=(1-d)n^2+2\Z$, such that for a given number $n$
 we have $\delta(n,m)\ne 0$ only for {\bf finitely many} values of $m$, and
$$P_d= \prod_{n\geqslant 1} \prod_{m\in \Z}\left( q^{m/2}z^n;q \right)_\infty^{ \delta^{(d)}(n,m)}\,\,.$$

\end{thm}

The above  implies the following decomposition
$$ \frac{P_d(z,q^{1/2})}{P_d(qz,q^{1/2})}=\prod_{n\geqslant 1}\prod_{m\in \Z}\prod_{i=0}^{n-1} \left(1-
q^{m/2+i}z^n  \right)^{
\delta^{(d)}(n,m)}  \,\,.$$
Therefore, the limit
$$ P_d^{cl}(z):= \lim_{q^{1/2}\to 1} \frac{P_d(z,q^{1/2})}{P_d(qz,q^{1/2})}\in 1+z \Z[[z]]$$ exists and has the form
$$P_d^{cl}=\prod_{n\geqslant 1}\left(1-z^n\right)^{nc^{(d)}(n)}\,,\,\,\,\,c^{(d)}(n)\in \Z\,,$$
where $c^{(d)}(n)=\sum_m \delta^{(d)}(n,m)$.
It follows from the result of Reineke (see also Section 5.6) that $P_d^{cl}$ is an algebraic series satisfying
$$P_d^{cl}(z)=1+(-1)^{(d-1)}z \left(P_d^{cl}(z)\right)^{d}\,\,.$$

\subsection{Symmetric case}

In the case when matrix $(a_{ij})_{i,j\in I}$ is symmetric one
 can repeat essentially all the considerations made for quivers with one vertex from the
 Section 2.5.

 First of all, the algebra $\mathcal{H}$ can be endowed with $\Z_{\geqslant 0}^I \times \Z$-grading:
  $$\mathcal{H}=\oplus_{\gamma,k} \mathcal{H}_{\gamma,k}\,,\,\,\,\,\mathcal{H}_{\gamma,k}:=H^{k-\chi_Q(\gamma,\gamma)}(\B\G_\gamma)\,\,,$$
where we recall $\chi_Q(\gamma_1,\gamma_2):=\sum_{i\in I} \gamma^i_1 \gamma^i_2 -\sum_{i,j\in I} a_{ij}\gamma^i_1\gamma^j_2$.
The explicit formula for the product implies that
 $$a_{\gamma,k}\cdot  a_{\gamma',k'}=(-1)^{\chi_Q(\gamma,\gamma')} a_{\gamma',k'}\cdot a_{\gamma,k}\,\,.$$

 This identity {\it does not} mean that $\mathcal{H}$ is supercommutative. We will modify the
 product on $\mathcal{H}$ in certain way, making it supercommutative.

 Define a group homomorphism $\epsilon:\Z^I\to \Z/2\Z$ by the formula
   $$\epsilon(\gamma):=\chi_Q(\gamma,\gamma)\,(\mathrm{mod}\,\,2)\,\,.$$
   The bilinear form
   $$ \Z^I\otimes \Z^I\to \Z/2\Z\,,\,\,\,\gamma_1\otimes \gamma_2\mapsto \left(\chi_Q(\gamma_1,\gamma_2)+\epsilon(\gamma_1)\epsilon(\gamma_2)\right)\,(\mathrm{mod}\,\,2)$$
    induces a symmetric bilinear form $\beta$ on $\Z/2\Z$-vector space $(\Z/ 2\Z)^I$
    which gives zero quadratic form.

    It is easy to see that there exists a bilinear form $\psi$ on $(\Z/ 2\Z)^I$
    such that
    $$\psi(\gamma_1,\gamma_2)+\psi(\gamma_2,\gamma_1)=\beta(\gamma_1,\gamma_2)\,\,.$$
    For example, choose an order  on $I$ and
    define $\psi$ on elements $(e_i)_{i\in I}$ of the standard basis of $(\Z/ 2\Z)^I$ by
    $$\psi(e_i,e_j)=\beta(e_i,e_j)\mbox{ if }i>j\,,\,\,\,\,\,\,\psi(e_i,e_j)=0\mbox { if } i\leqslant j\,\,.$$
    We define a modified product on $\mathcal{H}$ by
   $$a_1\star a_2:=(-1)^{\psi(\gamma_1,\gamma_2)} a_1\cdot a_2\,,\,\,\,
   \,\,\,a_1\in\mathcal{H}_{\gamma_1},\,a_2\in \mathcal{H}_{\gamma_2}\,\,.$$
   The conditions on $\psi$ ensure that the new product is again associative, and
    $(\mathcal{H},\star)$ is $\Z_{\geqslant 0}^I \times \Z$-graded supercommutative algebra,
     with the parity of the component $\mathcal{H}_{\gamma,k}$ given by $\epsilon(\gamma)$. Equivalently, one can  use the parity given by $k \,(\mathrm{mod}\,\,2)$.
      It does not change the super structure on $\mathcal{H}$ because all  non-trivial cohomology groups
       of $\B\G_\gamma$ are concentrated in even degrees, hence
       $$k=\chi_Q(\gamma,\gamma)=\epsilon(\gamma)\,(\mathrm{mod}\,\,2)\,\,.$$
One can see easily that different choices of cocycles $\psi$ lead to canonically isomorphic graded supercommutative algebras.

Here we formulate a general conjecture\footnote{After the initial version of this text was posted on arxive, A.~Efimov gave a  proof of our conjecture in \cite{Efimov}).}.
\begin{conj} In the symmetric case the $\Z_{\geqslant 0}^I\times \Z$-graded  algebra $\mathcal{H}\otimes \Q$ is a free
 supercommutative algebra generated by a graded vector space $V$ over $\Q$ of the form $V=V^{\mathrm{prim}}
\otimes \Q[x]$, where $x$ is an even variable of
bidegree $(0,2)\in \Z_{\geqslant 0}^I\times \Z$, and for any given $\gamma$ the space $V^{\mathrm{prim}}_{\gamma,k}\ne 0$ is non-zero (and finite-dimensional) only for {\bf finitely many} $k\in \Z$.
\end{conj}

In the special case of quiver with one vertex (see Section 2.5) this conjecture implies that all integers $\delta^{(d)}(n,m)$ in the decomposition in Theorem 3 in Section 2.5 have definite signs:
  $$  (-1)^{m-1}  \delta^{(d)}(n,m)\geqslant 0\, .             $$
One  evidence for the above Conjecture is the admissibility property of the generating series from Sections 6.1 and  6.2 (see e.g. Theorem 10 and Corollary 3 below). Those series
generalize the product from Theorem 3.

         \subsection{Grading by the Heisenberg group and twisted graded algebras}

In general, for not necessarily symmetric matrix $(a_{ij})_{i,j\in I}$, the interaction of the product in $\mathcal{H}$ with the cohomological grading can be treated in the following way. Let $\Gamma$ be an abelian group and $B:\Gamma\otimes \Gamma\to \Z$
 be a bilinear form  (in our case $\Gamma=\Z^I$ is the group of dimension vectors and $B=\chi_Q$). We associate with $(\Gamma, B)$ the discrete Heisenberg group
 $\mathrm{Heis}_{\,\Gamma,B}$ which is the set $\Gamma\times \Z$ endowed with the multiplication
$$(\gamma_1,k_1)\cdot (\gamma_2,k_2):=(\gamma_1+\gamma_2,k_1+k_2-2B(\gamma_1,\gamma_2))\,\,.$$
Our Cohomological Hall algebra $\mathcal{H}$ can be graded by the non-commutative group $\mathrm{Heis}_{\,\Gamma,B}$:
  $$\mathcal{H}=\oplus_{\gamma,k} \mathcal{H}_{(\gamma,k)}\,,\,\,\,\mathcal{H}_{(\gamma,k)}=H^k(\B\G_\gamma)\,\,.$$

This can be generalized such as follows. Let $(\mathcal{T},\otimes)$ be a symmetric monoidal category endowed with an {\bf even} invertible object $\mathbb{T}_\mathcal{T}$, i.e. such an object that the commutativity morphism $\mathbb{T}_\mathcal{T}\otimes \mathbb{T}_\mathcal{T}\to
\mathbb{T}_\mathcal{T}\otimes \mathbb{T}_\mathcal{T}$ is the identity morphism. Let us fix $(\Gamma,B)$  as above. We define a {\bf twisted graded monoid} \footnote{We use the word ``algebra" instead of ``monoid" for {\it additive} monoidal categories.} in $\mathcal{T}$
 as a collection of objects $(\mathcal{H}_\gamma)_{\gamma\in \Gamma}$ together with
  a collection of morphisms
 $$m_{\gamma_1,\gamma_2}: \mathcal{H}_{\gamma_1}\otimes \mathcal{H}_{\gamma_2}\to \mathbb{T}_\mathcal{T}^{\otimes B(\gamma_1,\gamma_2)} \otimes \mathcal{H}_{\gamma_1+\gamma_2}\,\,,$$
and a unit morphism $1_{\mathcal{C}}\to \mathcal{H}_0$, satisfying an obvious extension of the usual associativity and the unity axioms.
 In our situation $(\mathcal{T},\otimes)$ is the tensor category of $\Z$-graded abelian groups (with the usual Koszul sign rule for the
 commutativity morphism), and $\mathbb{T}_\mathcal{T}=\Z[-2]$ is the group $\Z$ placed in degree $+2$, i.e. the cohomology of the pair $H^\bullet(\C P^1,pt)$.
 Object $\mathcal{H}_\gamma$ is just $H^\bullet(\M^\univ_\gamma)=H^\bullet(\B\G_\gamma)$.

Let us assume additionally that the form $B$ is symmetric, and we are given an invertible object $\mathbb{T}^{\otimes 1/2}_\mathcal{T}\in \mathcal{T}$ such that
$\left(\mathbb{T}_\mathcal{T}^{\otimes 1/2}\right)^{\otimes 2}\simeq \mathbb{T}_\mathcal{T}$. Then for
 any twisted graded monoid $\H=(\H_\gamma)_{\gamma \in \Gamma}$ we define the modified graded monoid by
$$\H^{mod}_\gamma:=\H_\gamma\otimes \left(\mathbb{T}_\mathcal{T}^{\otimes 1/2}\right)^{\otimes B(\gamma,\gamma)},$$
which is an ordinary (untwisted) $\Gamma$-graded monoid in $\mathcal T$. In particular, for the category of $\Z$-graded abelian groups we take
 $\mathbb{T}^{\otimes 1/2}_\mathcal{T}$ to be $\Z$ in degree $+1$. In this way we obtain another description of the modified product $\star$ from the previous section. Notice that in this case object $\mathbb{T}^{\otimes 1/2}_\mathcal{T}$ is not even.

\subsection{Non-symmetric example: quiver $A_2$}

The quiver $A_2$ has two vertices $\{1,2\}$ and one arrow $1\leftarrow 2$. The Cohomological Hall algebra $\mathcal{H}$
 contains two subalgebras $\mathcal{H}_L,\,\mathcal{H}_R$ corresponding to representations supported at the vertices $1$ and $2$ respectively.
 Clearly each subalgebra $\mathcal{H}_L,\,\mathcal{H}_R $ is isomorphic to the  Cohomological Hall algebra for  the quiver $A_1=Q_0$. Hence it is
an infinitely generated exterior algebra (see Section 2.5). Let us denote the generators by $\xi_i,\, i=0,1,\dots$ for the vertex $1$ and by $\eta_i,\, i=0,1,\dots$ for the vertex $2$. Each generator $\xi_i$ or $\eta_i$ corresponds to an additive generator of the group $H^{2i}(BGL(1,\mathbb{C}))\simeq \mathbb{Z}\cdot x^i$. Then one can check that  $\xi_i,\eta_j, \,i,j\geqslant 0$ satisfy the relations
$$\xi_i\xi_j+\xi_j\xi_i=\eta_i\eta_j+\eta_j\eta_i=0\,,\,\,\,\,\eta_i\,\xi_j=\xi_{j+1}\eta_i-\xi_j\eta_{i+1}\, .$$
Let us introduce the elements $\nu_i^1=\xi_0\eta_i\,,\, i\geqslant 0$ and $\nu_i^2=\xi_i\eta_0\,, \,i\geqslant 0$. It is easy to see that $\nu_i^1\nu_j^1+\nu_j^1\nu_i^1=0$, and similarly the generators $\nu_i^2$ anticommute. Thus we have two infinite Grassmann subalgebras in $\mathcal{H}$ corresponding to these two choices: $\mathcal{H}^{(1)}\simeq \bigwedge(\nu_i^1)_{i\geqslant 0}$ and
$\mathcal{H}^{(2)}\simeq \bigwedge(\nu_i^2)_{i\geqslant 0}$. One can directly check the following result.
\begin{prp} The multiplication (from the left to the right) induces  {\bf isomorphisms} of graded abelian groups
$$\mathcal{H}_L\otimes \mathcal{H}_R\stackrel{\sim}{\longrightarrow}\mathcal{H}, \,\,\,\,\,\,
\mathcal{H}_R\otimes \mathcal{H}^{(i)}\otimes \mathcal{H}_L\stackrel{\sim}{\longrightarrow} \mathcal{H}\,,\,\,i=1,2\,. $$
\end{prp}

In other words, we obtain a canonical isomorphism (depending on $i=1,2$) of graded abelian groups
$$H^\bullet(\mathrm{BGL}(n,\C)\times \mathrm{BGL}(m,\C))\simeq $$
$$ \simeq\bigoplus_{0\leqslant k\leqslant \op{min}(n,m)} H^\bullet(\mathrm{BGL}(m-k,\C)\times
\mathrm{BGL}(k,\C)\times \mathrm{BGL}(n-k,\C))[-2(n-k)(m-k)]\,. $$
Passing to generating series we obtain the standard identity
$$ (q^{1/2}\widehat{\mathbf{e}}_1;q)_\infty \cdot (q^{1/2}\widehat{\mathbf{e}}_2;q)_\infty=(q^{1/2}\widehat{\mathbf{e}}_2;q)_\infty\cdot (q^{1/2}\widehat{\mathbf{e}}_{12};q)_\infty\cdot (q^{1/2}\widehat{\mathbf{e}}_1;q)_\infty\, ,$$
where non-commuting variables $\widehat{\mathbf{e}}_1,\widehat{ \mathbf{e}}_2,\widehat{ \mathbf{e}}_{12}$  satisfy relations of the Heisenberg group (with $-q^{1/2}$ corresponding to the central element):
$$\widehat{\mathbf{e}}_1\cdot \widehat{\mathbf{e}}_2= q^{-1}\,\widehat{\mathbf{e}}_2\cdot \widehat{\mathbf{e}}_1=-q^{-1/2}\widehat{\mathbf{e}}_{12}\, .$$
The same example of quiver $A_2$ was studied in \cite{KS}, section 6.4. The generating series $\mathbf{E}(q^{1/2},x)$ used there is related with the $q$-Pochhammer symbol 
by the following change of variables:
$$\mathbf{E}(-q^{1/2},x)=(q^{1/2}x;q)_\infty\,\,,$$
where both sides are viwed as formal series in $x$ with coefficients in the field $\Q(q^{1/2})$ of rational functions in variables $q^{1/2}$.

\subsection{Questions about other algebraic structures on $\mathcal{H}$}

One may wonder whether $\mathcal{H}$ carries other natural algebraic structures. For example,
it would be nice to have a structure of Hopf algebra on $\mathcal{H}$, as a step towards Conjecture 1 from Section 2.6, and also by analogy with quantum enveloping algebras.

One can see immediately that $\mathcal{H}$ carries a structure of  (super) cocommutative coassociative coalgebra with
  the coproduct
$$\mathcal{H}_{\gamma_1+\gamma_2}\to \mathcal{H}_{\gamma_1}\otimes \mathcal{H}_{\gamma_2}$$
 induced by the operation of the direct sum on representations:
 $$\M^\univ_{\gamma_1}\times \M^\univ_{\gamma_2}\to \M^\univ_{\gamma_1+\gamma_2}\,\,,$$
which is realized geometrically via any embedding $\C^\infty\oplus \C^\infty\mono \C^\infty$.
 Unfortunately, this coproduct is not compatible with the product on $\mathcal{H}$.
 In Section 6.5 we will define  a supercommutative associative {\it product} on $\mathcal{H}$, which is
 also related to the operation of the direct sum, and used in the study of motivic  DT-series of
 $\mathcal{H}$. We do not know  any reasonable relation between this supercommutative product
 and the product introduced in Section 2.2.

Another potentially interesting structure  is associated with the operation of taking multiple copies of an object.
 This gives a morphism of cohomology
$$\mathcal{H}_{n\gamma}\to \mathcal{H}_\gamma\otimes H^\bullet(\mathrm{BGL}(n,\C))$$ for any $n\geqslant 0$. It comes from the map
  $$\M^\univ_{\gamma}\times \mathrm{BGL}(n,\C)\to \M^\univ_{n \gamma}\,\,,$$
and is realized geometrically via any embedding $\C^\infty\otimes \C^n\mono \C^\infty$.

\section{Generalization of COHA for smooth algebras}

Up to now, homotopy types  of the moduli stacks of representations were all very simple, namely,
 just products of classifying spaces of groups $\mathrm{GL}(n,\C),\,n\geqslant 0$.
In this section we are going to discuss several generalizations of the Cohomological Hall algebra, in particular for  representations of smooth algebras.
Also, we will use not only the usual cohomology with coefficients in $\Z$, but  other cohomology groups of algebraic varieties, like e.g. de Rham or \'etale cohomology.

\subsection{Cohomology theories}

Here we collect preliminary material about cohomology theories for algebraic varieties and fix some notations.
\begin{dfn} Let $\kk$ be a field, and $K$ be another field, $\op{char}(K)=0$. Assume that we are given a $K$-linear Tannakian category
 $\mathcal C$. A {\bf pre-cohomology} theory over $\mathbf{k}$ with values in $\mathcal C$ is a contravariant tensor
functor ${\bf H}^\bullet$ from the category of schemes of finite type over $\kk$ to the tensor category ${\mathcal C}^{\Z-gr}$ of $\Z$-graded objects in  $\mathcal C$ endowed with Koszul rule of signs.
\end{dfn}

Here is the list of standard examples of (pre-)cohomology theories, $\mathcal{C}=K-\op{mod}$:
\begin{itemize}
\item (case $\op{char}\kk=0, \,\,\kk\subset \C $): rational Betti cohomology, $K=\Q$,
\item (case $\op{char}\kk=0$): de Rham cohomology, $K=\kk$,
\item  (case $\op{char}\kk\ne l$ for prime $l$): \'etale cohomology, $K=\Q_l$.
\end{itemize}

In the case $\kk\subset \C$ and Betti cohomology one can enhance $\mathcal C$  from $\Q-\op{mod}$ to the Tannakian category of polarizable mixed Hodge structures.
Similarly, in the \'etale case we can take $\mathcal C$ to be the category of  continuous $l$-adic representations of the absolute Galois group of $\kk$.

\begin{dfn} A pre-cohomology theory
 ${\bf H}^\bullet$ is called  {\bf cohomology} theory  if there exists a fiber functor $\mathcal{C}\to K'-\op{mod}$ for some field extension  $K'\supset K$
 such that the resulting pre-cohomology theory is obtained by the extension of scalars from one of the  standard theories.
\end{dfn}

All the above examples give in fact cohomology theories.
Notice that in our definition one can mention only \'etale cohomology, since in the case $\op{char}\kk=0$ we have a chain of comparison isomorphisms
 $$\mbox{ de Rham }\simeq \mbox{ Betti }\simeq \mbox{ \'etale } .$$
In fact, in what follows we will need only a pre-cohomology theory satisfying the usual properties of cohomology, and the definition
 above is given in order to avoid  a long discussion of cohomology theories. The optimal framework should be that of triangulated category of motives, see
Section 3.5.

We denote $\mathbf{H}^2(\mathbb{P}^1_\kk)$ by  $K(-1)$ understood as an element of $\mathcal C$, and set
 $$\mathbb{T}=\mathbb{T}_{\mathbf{H}^\bullet}:=K(-1)[-2]\in \mathcal{C}^{\Z-gr} .$$

         \subsection{Smooth algebras}

Here we recall  the notion of a  smooth algebra from \cite{KR}.
\begin{dfn} An associative unital algebra $R$ over a field $\kk$ is called {\bf smooth} if it is finitely generated and formally smooth in the sense
of D.~Quillen and J.~Cuntz, i.e. if
the bimodule
 $\Omega^1_R:=\op{Ker}(R \otimes_\kk R\stackrel{\op{mult}}{\longrightarrow}R)$ is projective. Here
 $\op{mult}:R \otimes_\kk R\to R$ is the product.
\end{dfn}
 For a finitely presented algebra $R$ it is easy to check smoothness using a finite amount of calculations, see \cite{CQ} or \cite{KR}, Section 1.1.3.
The property of formal smoothness is equivalent to the lifting property for non-commutative nilpotent extensions:
for any associative unital algebra $A$ over $\kk$ and a nilpotent two-sided ideal $J\subset A,\,\,J^n=0$ for some $n>0$, and any homomorphism
$\phi:R\to A/J$ there exists a lifting of $\phi$ to a homomorphism $R\to A$.

Basic examples of smooth algebras are matrix algebras, path algebras of finite quivers, and algebras of functions on smooth affine curves.

The notion of a smooth algebra is invariant under Morita equivalence, and closed under free products and finite localizations.
 The latter means that we can add two-sided inverses to a finite collection of elements in $R$ (or, more generally, to a finite
 collection of rectangular matrices with coefficients in $R$). For example, the group ring of a free finitely generated group
$$\kk[\op{Free}_d]=\kk\langle x_1^{\pm 1},\dots, x_d^{\pm 1}\rangle$$
is a smooth algebra.

For a smooth algebra $R$ and for any finite-dimensional algebra $T$ over $\kk$ (e.g. $T=\op{Mat}(n\times n,\kk)$), the scheme of homomorphisms $R\to T$ is a smooth affine scheme. 
Indeed, this is a closed subscheme of an affine space $\mathbb{A}_\kk^{g_R \dim T}$ where $g_R$ is number of generators of $R$, hence
 it is an affine scheme of finite type over $\kk$. Also this scheme enjoys the lifting property for (commutative)
nilpotent extensions (as follows from  formal smoothness), hence it is smooth.

If $R$ is smooth then the abelian category of $R$-modules is hereditary, i.e.
   $\op{Ext}^i(E,F)=0$ for for any two objects $E,F$ and  any $i\geqslant 2$. If
  $E,F$ are finite-dimensional over $\kk$ then spaces $\op{Hom}(E,F),\,\op{Ext}^1(E,F)$ are also finite-dimensional.

\subsection{Cohomological Hall algebra for an $I$-bigraded smooth algebra}

For a finite set $I$,
 we call
 an unital associative algebra $R/\kk$  {\bf $I$-bigraded} if $R$ is decomposed (as a vector space) into the direct sum
$R=\oplus_{i,j\in I} R_{ij}$ in such a way that $R_{ij}\cdot R_{jk}\subset R_{ik}$. Equivalently, $R$ is $I$-bigraded if we are given a morphism of unital algebras $\kk^I\to R$.

Let now $R$ be an $I$-bigraded {\it smooth} algebra. It follows that  any finite-dimensional representation $E$ of $R$ decomposes into a direct sum of finite-dimensional vector spaces $E_i, i\in I$.

 For any dimension vector $\gamma=(\gamma^i)_{i\in I}\in {\Z}_{\geqslant 0}^I$ the scheme $\M_\gamma=\M_\gamma^R$ of representations of $R$ in coordinate spaces $E_i=\kk^{\gamma^i}, {i\in I}$ is a smooth affine scheme. Any choice of a finite set of $I$-bigraded generators of $R$ gives a closed embedding
 of $\M_\gamma$ into the affine space $\M_\gamma^Q$ for some quiver $Q$ with the set of vertices equal to  $I$.

Let us make the following {\bf Assumption}:\newline

{\it We are given a bilinear form $\chi_R:\Z^I\otimes \Z^I\to \Z$ such that for any two dimension vectors $\gamma_1,\gamma_2\in \Z_{\geqslant 0}^I$ and for any two representations $E_i\in \M_{\gamma_i}(\overline{\kk})$ we have the equality}
$$\op{dim}\op{Hom}(E_1,E_2)-\op{dim}\op{Ext}^1(E_1,E_2)=\chi_R(\gamma_1,\gamma_2)\, .$$
Here $\overline{\kk}$ is an algebraic closure of $\kk$, and $E_1, E_2$ are considered as representations of algebra $R\otimes_\kk \overline{\kk}$ over $\overline{\kk}$.
The assumption implies that the smooth scheme $\M_{\gamma}$ is {\it equidimensional} for any given $\gamma$  and
$$ \op{dim} \M_{\gamma}=-\chi_R(\gamma,\gamma)+\sum_i (\gamma^i)^2 .$$

For example, if the algebra $R$ is obtained from the path algebra of a finite quiver $Q$ by a finite number of localizations, then $\chi_R=\chi_Q$ and the scheme  $\M_\gamma$ is a Zariski  open subset of the affine space $\M_\gamma^Q\simeq {\mathbb{A}}^{\sum_{ij} a_{ij} \gamma^i\gamma^j}_\kk$.

Now we can use any cohomology theory $\mathbf{H}^\bullet$ and obtain the corresponding Cohomological Hall algebra, which is
 a twisted associative algebra (in the sense of Section 2.7) in tensor category $\mathcal{C}^{\Z-gr}$
with
$$\mathbb{T}_{\mathcal T}:=\mathbb{T}=K(-1)[-2]\, .$$
 The main point is that spaces $\M_{\gamma_1,\gamma_2}$ of representations of $R$
  in block upper-triangular matrices, are {\it smooth closed} subvarieties of $\M_{\gamma_1+\gamma_2}$. The definition of the product and the proof of associativity given in Sections 2.2 and 2.3
 work without  modifications in the  general case of smooth algebras and cohomology theories. In this way we obtain COHA $\mathcal{H}:=\mathcal{H}_R$ for a smooth $I$-bigraded algebra $R$ which satisfies the above Assumption.

\subsection{Square root of Tate motive and modified product}

Let $\mathcal T$ be a rigid tensor category over field $K$, and $E$ be an even invertible object of $\mathcal T$.
 Then one can formally add to $\mathcal T$ an {\it even} tensor square root $E^{\otimes 1/2}$ (see Section 2.7).  Moreover, if $\mathcal T$
is Tannakian, then the new tensor category $\mathcal{T}(E^{\otimes 1/2})$ is also Tannakian. This can be explained explicitly as follows.
 Let us assume that we are given a fiber functor $\mathcal{T}\to K-\op{mod}$, i.e. identify  $\mathcal T$ with the tensor category
 of finite-dimensional  representations of a pro-affine algebraic group $G$ over $K$. Object $E$ gives a homomorphism
 $t:G\to \mathbb{G}_{m,K}\simeq Aut(E)$.
Denote by $G^{(2)}$ the fibered product
$$G^{(2)}:=\varprojlim \left(\begin{CD}@.\mathbb{G}_m\\
         @. @V{\lambda\mapsto \lambda^2}VV \\
         G @>t>> \mathbb{G}_m\end{CD}\right)\,\,.$$
Then $\mathcal{T}(E^{\otimes 1/2})$ is canonically equivalent to the tensor category of representations of $G^{(2)}$.
 It is easy to see that the ring $K_0(G^{(2)}-\op{Rep})$ is canonically isomorphic
            to $K_0(G-\op{Rep})[x]/(x^2-[E])$ where  $x\mapsto [E^{\otimes 1/2}]$.

One can apply these considerations to the  Tannakian category $\mathcal C$ associated with
 a cohomology theory ${\bf H}^\bullet$. We take $E:=K(-1)$. The square root $E^{\otimes 1/2}$ can be denoted by $K(-1/2)$.
 Usually there is no square root of $K(-1)$ in $\mathcal C$. For classical cohomology it exists only when $\kk$ has positive characteristic
 and a square root is one of summands of $\mathbf{H}^1(C)$ where $C$ is a supersingular elliptic curve.
 In Section 4 we will develop a generalization of cohomology for varieties with potential. In this more general framework a square root of $K(-1)$
 exists much more often, e.g. it is enough to assume that $\kk$ contains a $\sqrt{-1}$.
Anyhow, in what follows we will assume that the Tannakian category $\mathcal C$ for the cohomological theory {\it contains}  a chosen
 square root $K(-1/2)$ (if it can not be found in the original category, we will add it formally as explained above).
 We will also denote
$$\mathbb{T}^{\otimes 1/2}:=K(-1/2)[-1]\in \mathcal{C}^{\Z-gr}.$$

Let $\mathcal{H}=\oplus_{\gamma\in \Z^I_{\geqslant 0}}\,\mathcal{H}_{\gamma}$ be the Cohomological Hall algebra of an $I$-bigraded smooth algebra $R$.

\begin{dfn} The {\bf modified} Cohomological Hall algebra $\mathcal{H}^{mod}$
is given by
$$\mathcal{H}^{mod}:=\oplus_{\gamma\in \Z^I_{\geqslant 0}}\,\mathcal{H}^{mod}_\gamma\,,\,\,\,\,\,
\mathcal{H}^{mod}_\gamma:=\mathcal{H}_\gamma\otimes \left(\mathbb{T}^{\otimes 1/2}\right)^{\otimes \chi_R(\gamma,\gamma)}  .$$
\end{dfn}

For any sublattice $\Gamma\subset \Z^I$ such that the restriction of $\chi_R$ to ${\Gamma}$ is symmetric, the restricted modified algebra
$$\bigoplus_{\gamma\in \Gamma\cap \Z^I_{\geqslant 0}} \mathcal{H}^{mod}_\gamma$$
is an associative $\Gamma$-graded algebra in tensor category $\mathcal{C}^{\Z-gr}$.

\subsection{Algebra in the triangulated category of motives}

We can try to use  a more refined object instead of graded cohomology spaces (for some classical
cohomology theory). A natural candidate is the motive of the ind-schemes $\M_\gamma^\univ$
 understood as an element of Voevodsky triangulated category of motives (see \cite{Voevodsky}). More precisely, we should maybe use a dg enhancement
 constructed in \cite{Bondarko} and in \cite{BeilinsonVologodsky}. We expect that the Cohomological Hall algebra has a natural enhancement  to an $A_\infty$-algebra
 in the tensor dg-category of mixed motives, with higher multiplications somehow governed by manifolds of flags
 of representations. Those are the spaces $\M_{\gamma_1,...,\gamma_n}$ from the proof of Theorem 5 in Section 5.2 below.
The invertible object $\mathbb{T}$ used in the definition of a twisted associative algebra (see Section 2.7)
is the graded space $\Z(-1)[-2]$, the Tate motive of the cohomology of pair $(\mathbb{P}^1,pt)$.

\subsection{Equivariant parameters in arrows and families of quivers}

Let us consider for simplicity the case of the path algebra of a quiver.
One can introduce a torus action on the path algebra by rescaling independently individual arrows. Then the representation
varieties will be endowed with the action of an additional torus $\mathbf{T}$, and one can use the $\mathbf{T}$-equivariant cohomology. In this way we define COHA as $H^{\bullet}(B\mathbf{T})$-module.
This is a particular case of a more general construction, when one has a family over some base $\mathcal{B}$ of smooth $I$-bigraded algebras. The  case of $\mathbf{T}$-action corresponds to $\mathcal{B}:=B\mathbf{T}$ (the classifying space of $\mathbf{T}$).

         \subsection{Complex cobordisms}

The natural framework for the theory of virtual fundamental classes in the case of almost complex manifolds
 and complex obstruction bundles is not the usual cohomology theory but the complex cobordism theory. For an element $(f:Y\to X)\in \Omega_n(X)$ of the $n$-th bordism group of $X$, where $Y$ is a compact manifold with stable complex structures, and a homogeneous monomial of characteristic classes $\alpha=c_1(T_Y)^{a_1}...c_k(T_Y)^{a_k}$, one can define the corresponding element $f_{\ast}[P.D.(\alpha)]\in H_{\bullet}(X)$. Here $P.D.(\alpha)$ is the Poincar\'e dual homology class to $\alpha$.   This means that after tensoring by $\Q$
 we replace the cohomology by the tensor product of cohomology with $H^\bullet(BU(\infty))$ (the graded version of the algebra of symmetric polynomials). In order to pass back to  complex bordisms groups we should keep track of the action of characteristic classes of the virtual tangent bundles
of the Grassmannians and of the moduli stacks on their virtual fundamental classes.
Therefore, if we define the Cohomological Hall algebra using complex bordisms instead of cohomology groups, we obtain $H^\bullet(BU(\infty),\Q)$-linear product on $H^\bullet(BU(\infty),\Q)\otimes \H$. It can be considered as a deformation of the original COHA with the base $\op{Spec}H^\bullet(BU(\infty),\Q)\simeq \mathbb{A}^\infty_\Q$.

\section{Exponential Hodge structures and COHA for a smooth algebra with potential}

\subsection{Short synopsis}

\begin{dfn} For a complex algebraic variety $X$ and a function $f\in \mathcal{O}(X)$ regarded as a regular map $f:X\to \C$, we define
 the {\bf rapid decay} cohomology $H^\bullet(X,{f})$ as the limit of the cohomology of the pair $H^\bullet(X,f^{-1}(S_t))$ for real $t\to -\infty$, where
$$S_t:=\{ z\in \C\,|\,\op{Re}z<t\} \,\,.$$
\end{dfn}
The cohomology stabilizes at some $t_0\in \R,\,\,t_0\ll 0$ (also in the definition one can replace $f^{-1}(S_t)$ by $f^{-1}(t)$).
 The origin of the term  ``rapid decay'' will become clear later in Section 4.5.

  The cohomology $H^\bullet(X,f)$ behaves similarly to the usual cohomology. In particular, for a map $\pi:Y\to X$ compatible with functions
 $f_Y\in \mathcal{O}(Y),\,f_X\in \mathcal{O}(X)$ in the sense that $f_Y=\pi^* f_X$,  we have the pullback
$\pi^*:H^\bullet(X,f_X)\to H^\bullet(Y,f_Y)$. If $\pi$ is proper and both $X$ and $Y$ are smooth, then we have the pushforward morphism
$\pi_*:H^\bullet(Y,f_Y)\to H^{\bullet+2(\dim_\C(X)-\dim_\C(Y))}(X,f_X)$.

Similarly to the usual cohomology, there is a  multiplication (K\"unneth) morphism
  $$\otimes:H^\bullet(X,{f_X})\otimes H^\bullet(Y,{f_Y})\to H^\bullet(X\times Y, {f_X\boxplus f_Y})\,\,,$$
where the Thom-Sebastiani sum $\boxplus$ is given by
$$ f_X\boxplus f_Y:=\mathrm{pr}_{X\times Y\to X}^* f_X+  \mathrm{pr}_{X\times Y\to Y}^* f_Y \,\,.$$

For any smooth $I$-bigraded algebra $R$ as in Section 3.3 and any element
 $W\in R/[R,R]$ we obtain a function $W_\gamma$ on $\M_\gamma$ which is invariant under $\G_\gamma$-action. It is given by the trace of $W$ in a representation.
Hence we can apply the formalism of Sections 2.1-2.3 to rapid decay cohomology and define
 the Cohomological Hall algebra as
$$\H=\oplus_\gamma \H_\gamma\,,\,\,\,\,\H_\gamma:=H^\bullet_{\G_\gamma}(\M_\gamma,W_\gamma)\,.$$
In the case when $R$ is the path algebra of a quiver $Q$, the element $W$ can be considered
as a linear combination of cyclic paths in $Q$. It can be thought of as a potential for a general multi-matrix model (see e.g. \cite{Eynard-1}).

In the next sections we will define an analog of mixed Hodge structure on $H^\bullet(X,{f})\otimes \Q$,
 in particular an analog of the {\it  weight filtration}. There is a topological formula
for the weight filtration  similar to the classical case.

\subsection{Exponential mixed Hodge structures}

The goal of this section is to define a cohomology theory for  pairs $(X,f)$ consisting of a
 smooth algebraic variety $X$ over $\C$ and a function $f\in \mathcal{O}(X)$, with values in certain
Tannakian $\Q$-linear category $EMHS$ (see Definition 7 below) which we call the category of {\it
exponential mixed Hodge structures}. Definition of the category $EMHS$ will be based on the properties of two equivalent Tannakian categories $\mathfrak{A}_0$ and $\mathfrak{B}_0$ described below.
The discussion below is essentially a reformulation of well-known results of N.~Katz (see \cite{Katz}, Chapter 12) on the additive convolution of $\D$-modules and perverse sheaves on the affine line
 $\AC$.

First, consider the abelian category $\mathfrak{A}=\op{Hol}_{rs}(\mathbb{A}^1_\C)$ of holonomic $\mathcal{D}$-modules with regular singularities on the standard affine line endowed
with coordinate $x$.
 By Fourier transform $FT$ this category is identified with the category $\mathfrak{B}\simeq
FT(\mathfrak{A})$
of holonomic $\D$-modules $M$ on the dual line (endowed with the canonical coordinate $y=-d/dx$) such that the direct image of $M$ to
$\C \mathbb{P}^1$
 has a regular singularity at $0$, no singularities at  $\mathbb{A}^1_\C-\{0\}$, and a possibly irregular singularity at $\infty$ of the {\bf exponential type}. The latter condition means
 that after the base change to $\C((y^{-1}))$ we have
 $$M\otimes_{\C[y]}\C((y^{-1}))\simeq
\bigoplus_i \mathfrak{exp}(\lambda_i y)\otimes_{\C((y^{-1}))} M_i\,,$$
where $\lambda_i\in \C$ run through a finite set $I$, and $M_i$ are regular holonomic $\D$-modules on the formal punctured disc $\op{Spf}(\C((y^{-1})))$.
 Here $\mathfrak{exp}(\lambda_i y)$ denotes the $\D$-module on $\op{Spf}(\C((y^{-1})))$ which is a free module over $\C((y^{-1}))$ with
generator $e$ satisfying the same equation
$$(d/dy -\lambda_i)\cdot e=0$$
as the exponential function $\exp(\lambda_i y)$.

Let us denote by $\mathfrak{B}_0$ the abelian category of holonomic $\D$-modules on $\mathbb{A}^1_\C-\{0\}$ with no singularities at $\mathbb{A}^1_\C-\{0\}$, such that its direct image to $\C\mathbb{P}^1$ has
 a regular singularity
at $0$ and an exponential type singularity at $\infty$. In particular, any object of $\mathfrak{B}_0$ is an algebraic vector bundle
 over $\mathbb{A}^1_\C-\{0\}=\op{Spec}\C[y,y^{-1}]$ endowed with a connection.
The inclusion
$$j:\mathbb{A}^1_\C-\{0\}\mono \mathbb{A}^1_\C$$
 gives two adjoint exact functors
$$j^*:\mathfrak{B}\to \mathfrak{B}_0\,,\,\,\,j_*:\mathfrak{B}_0\to \mathfrak{B} $$
such that $j^* \circ j_*=\op{id}_{\mathfrak{B}_0  }$ (hence $\mathfrak{B}_0$ can be considered as a full abelian subcategory of $ \mathfrak{B}$)
and $\Pi:=j_*\circ j^*$ is an exact idempotent endofunctor of $\mathfrak{B}$, given by
the tensor product  $M\mapsto \C[y,y^{-1}]\otimes_{\C[y]}M$. A $\D$-module $M\in \mathfrak{B}$ belongs to $\mathfrak{B}_0$ if and only if 
the operator of multiplication by $y$ is invertible in $M$.
In what follows we will use the same notation $\Pi$ for exact idempotent endofunctors in several abelian categories closely related with the category $\mathfrak{B}\simeq \mathfrak{A}$,
 we hope that this will not lead to a confusion.

 Obviously, the category $\mathfrak{B}_0$ is closed under the tensor product over $\C[y,y^{-1}]$, and is a Tannakian category over $\C$,
 with the fiber functor to complex vector spaces
given by the fiber at any given non-zero point in the  line $Spec\, \C[y]$. More precisely,
we obtain a local system of fiber functors over $\C-\{0\}$. By homotopy invariance, it gives rise to
 a local system of fiber functors over the circle $S^1$ of tangent directions at $\infty\in \C\mathbb{P}^1$ with coordinate $y$.

Translating everything back by the inverse Fourier transform, we obtain a full abelian subcategory
 $$\mathfrak{A}_0\subset \mathfrak{A}\,,\,\,\,FT(\mathfrak{A}_0)=\mathfrak{B}_0$$
equivalent to $\mathfrak{B}_0$. This category
consists of holonomic $\D$-modules $N$ on $\mathbb{A}^1_{\C}$ with regular singularities,
 such that the operator $d/dx$ is invertible on $N$. In other words,
$$\op{RHom}_{D^b_{hol}(\D_{\mathbb{A}^1_{\C}})}(\mathcal{O}_{\mathbb{A}^1_{\C}},N)=H^\bullet_{DR}(N)=0\,.$$
 The category $\mathfrak{A}_0$ is closed under the additive convolution
 in $D^b(\D_{\mathbb{A}^1_{\C}}-\op{mod})$ given by
$$N_1 *_+ N_2:= {sum}_*(N_1\boxtimes N_2)\,,$$
where ${sum}:\C\times \C\to \C$ is the addition morphism $(x_1,x_2)\mapsto x_1+x_2$.
Moreover, $\mathfrak{A}_0\subset \mathfrak{A}$ is the image of an exact idempotent functor $\Pi:\mathfrak{A}\to \mathfrak{A}$
given by the convolution with $j_! \mathcal{O}_{\mathbb{A}^1_\C-\{0\}}[1]=FT^{-1}(\C[y,y^{-1}])$,
$$\Pi(N)= N*_+j_!(\mathcal{O}_{\mathbb{A}^1_\C-\{0\}}[1])\,.$$

Duality functor on the tensor category $\mathfrak{A}_0$ is given by
$$D_{\mathfrak{A}_0} N=\Pi(r_*(D_{\D_{\mathbb{A}^1_{\C}}-\op{mod}} N))\,,$$
where $r:\AC\to \AC,\,\,r(x)=-x$ is the antipodal map, and $D_{\D_{\mathbb{A}^1_{\C}}-\op{mod}}$ is the standard duality in $D^b(\D_{\mathbb{A}^1_{\C}}-\op{mod})$.

Tensor category $\mathfrak{A}_0$ also admits a local system of fiber functors over the circle $S^1$ of tangent directions at the point $\infty$ in the complex line endowed with the coordinate $x$.
Namely, any object $N\in \mathfrak{A}_0$ gives an analytic $\D$-module  $N^{\op{an}}$, which is a complex analytic vector bundle
 with connection $\nabla$
on a complement to a finite subset of $\C$. The fiber functor associated with an angle $\phi\in \R/2\pi\Z$ assigns to $N$ the space of flat sections of $N^{\op{an}}$ on the ray
 $[R,+\infty)\cdot \exp(i\phi)$ for sufficiently large $R\gg 0$.

Here we describe a canonical isomorphism between the two local systems of fiber functors: the one described in terms of the category $\mathfrak{A}_0$ and another one described in terms of
the category $\mathfrak{B}_0$. What will follow is basically a reformulation of the classical results of B.~Malgrange on the comparison isomorphism between the topological Fourier-Laplace transform and the algebraic Fourier transform for regular holonomic $\mathcal{D}$-modules  on $\mathbb{A}^1(\C)$, see his book \cite{Malgrange} devoted to this subject.

For any $N\in \mathfrak{A}$ and a given non-zero point $y_0\in \C-\{0\}$ on the dual line,
the fiber of $FT(N)$ at $y_0$ is equal to
$$H^1_{DR}(\op{exp}_{y_0}\cdot N )=\op{Coker}(d/dx+y_0:N\to N)\,,$$
where $\op{exp}_{y_0}$ is analytic function $x\mapsto \op{exp}(xy_0)$ on $\C$.
Suppose that we are given a vector $v$ in this fiber, represented in the de Rham complex by $\op{exp}_{y_0} P dx$ for some $P\in N$. Then for any point $x_t\in \C$ lying outside of singularities
 of $N$, the same expression $ \op{exp}_{y_0} P dx$ can be considered as representing a class in
 $$ H^1_{DR}( \op{exp}_{y_0} \cdot N[x_t])\, ,\,\,\,N[x_t]:=
N\otimes_{\C[x]}\C[x,1/(x-x_t)]\,.$$
Passing to analytic functions, we obtain a class in
$$H^1_{DR}(\op{exp}_{y_0}\cdot N[x_t]^{\op{an}})\simeq H^1_{DR}(N[x_t]^{\op{an}})\simeq H^1_{DR}(N[x_t])\,,$$
where the first isomorphism is given by multiplication by $\op{exp}_{-y_0}$, and the second one follows from the assumption that $N$ has regular singularities. Now, if we assume that $N\in \mathfrak{A}_0$, i.e. $H^\bullet_{DR}(N)=0$, then
$H^1_{DR}(N[x_t])$ coincides with the fiber of $N$ at $x_t$. Hence, varying $x_t$ we obtain a
 holomorphic section of the $D$-module $N^{\op{an}}$ outside of the set of singularities. This section is not flat, but is exponentially close
to a unique flat section along the ray $\mathbb{R}_{\leqslant 0}\cdot y_0^{-1} $ at infinity. Moreover, this flat section does not depend on the choice of a representative $P$ of the cohomology class, and gives a vector  $$v'\in \Gamma([R,+\infty)\cdot \exp(-i\phi),(N^{\op{an}})^\nabla)$$ corresponding to $v$, where $R\gg 0$ and $\phi=\op{Arg}(-y_0)$. The correspondence $v\mapsto v'$ gives an identification of two fiber functors evaluated at the object $N$.

 An informal meaning of the above construction is that we integrate differential $1$-form $\op{exp}_{y_0} P dx$ over an appropriate  non-compact chain with coefficients in the constructible sheaf associated with the dual object $D_{\D_{\mathbb{A}^1_{\C}}-\op{mod}}(N)$ via the Riemann-Hilbert correspondence. The integral is convergent because of the exponential decay of $\op{exp}_{y_0}$ along the chain. The role of the point 
$x_t$ in the above construction is to approximate the integral over a non-compact chain by an integral over a compact chain with boundary
 at $x_t$ where the exponential factor $\op{exp}(x_t y_0)$ is very small.

Let us apply the Riemann-Hilbert correspondence. The category $\mathfrak{A}$ is equivalent to the abelian category
 of middle perversity constructible sheaves of $\C$-vector spaces on $\mathbb{A}^1(\C)$. The subcategory $\mathfrak{A}_0$ goes to
 the  full subcategory $\op{Perv}_0(\mathbb{A}^1(\C),\C)$ consisting\footnote{In \cite{Katz} the subcategory $\op{Perv}_0(\mathbb{A}^1(\C),\C)$  was denoted by
 $\op{Perv}_A(\C)$. Moreover, in \cite{KaKoPa} it was proven that any object of $\mathfrak{A}_0$ is in fact a usual constructible sheaf shifted by $[1]$.} of objects $E^\bullet$ such that
 $R\Gamma(E^\bullet)=0$.
The additive convolution $*_+$
$$D^b(\op{Perv}(\mathbb{A}^1(\C),\C))\times D^b(\op{Perv}(\mathbb{A}^1(\C),\C))\to D^b(\op{Perv}(\mathbb{A}^1(\C),\C))$$
is given by the same formula as for $\D$-modules.
 Projector $\Pi$ in the category of perverse sheaves of $\C$-vector spaces is given by
$$ \Pi(F)=F*_+ j_!(\C_{\mathbb{A}^1(\C)}[1])\,.$$
The conclusion is that the Tannakian category $\mathfrak{B}_0$ is equivalent to the abelian Tannakian $\C$-linear
category $\op{Perv}_0(\mathbb{A}^1(\C),\C))$ which is defined purely topologically, endowed with the tensor product  given in terms of constructible sheaves. The latter category has a local system of fiber functors over the circle $S^1$ about the point $\infty$ given by stalks (shifted by [-1]) at points approaching $\infty$ along straight rays.

It is clear that one can define an exact idempotent functor 
$$\Pi(F)=F*_+ j_!(\Q_{\mathbb{A}^1(\C)}[1])$$
 on the  $\Q$-linear category $\op{Perv}_0(\mathbb{A}^1(\C),\Q)\subset D^b(\op{Perv}(\mathbb{A}^1(\C),\Q))$ of  perverse sheaves of  $\Q$-vector spaces.
The convolution on $\op{Perv}_0(\mathbb{A}^1(\C),\Q)$ is exact in each argument.
Also the category $\op{Perv}_0(\mathbb{A}^1(\C),\Q)$ is a Tannakian category, with the local system of fiber functors over $S^1$.

Now we can use the theory by M.~Saito of mixed Hodge modules (see \cite{Saito} or \cite{Steen}). For a smooth complex algebraic variety $X/\C$ we denote
 by $MHM_X$ the  category of  mixed Hodge modules on  $X$.
The category $MHM_{\mathbb{A}^1_\C}$ is a
$\Q$-linear abelian category endowed with an exact  faithful forgetting functor to $\op{Perv}(\mathbb{A}^1(\C),\Q)$ (Betti realization).
 This category is also endowed with an exact idempotent endofunctor $\Pi$ given by the addtitve convolution with $j_!(\Q_{\mathbb{A}^1(\C)}[1])$ considered as an object of 
$MHM_{\mathbb{A}^1_\C}$.

 \begin{dfn}
 The  category of {\bf exponential} mixed Hodge structures $EMHS$ is the full subcategory of the category $MHM_{\mathbb{A}^1_\C}$
consisting of objects $M$ such that the corresponding perverse sheaf belongs to $\op{Perv}_0(\mathbb{A}^1(\C),\Q)$. Equivalently, it is the image of endofunctor $\Pi$.
 It is a Tannakian category with the tensor product given by the additive convolution.
\end{dfn}

The tensor category of ordinary mixed Hodge structures $MHM_{\op{pt}}$ can be identified with a full Serre
 subcategory of $EMHS$ closed under the tensor product and duality. Namely, with any $M_0\in 
MHM_{\op{pt}}$
we associate  an object
$\Pi(s_* M_0)\in EMHS$ where $s:\{0\}\to \mathbb{A}^1(\C)$ is the obvious embedding.

\subsection{Realization functors for exponential mixed Hodge structures}

We have already mentioned that the $\C$-linear tensor category $\mathfrak{A}_0\simeq \mathfrak{B}_0$
 has {\it two} canonically equivalent local systems of realization functors, defined
 in terms of stalks of bundles with connections either for the original holonomic $\D$-module
with regular singularity or for its Fourier transform.

\begin{dfn} For an object $N\in EMHS$ and a non-zero complex number $u$ we define {\bf de Rham}
 realization of $N$ at the point $u$ to be the fiber of $\C[y]$-module $FT(N_{DR}^{MHM})$ at the point $u^{-1}\in \op{Spec}\C[y]$, where $N_{DR}^{MHM}$ is the  {\bf algebraic} holonomic $\D$-module underlying $N$, where $N$ is considered as an object of $MHM_{\mathbb{A}^1_\C}$. We denote this fiber by $N_{DR,u}$.
\end{dfn}

It follows directly from the definition that we have a natural isomorphism
$$N_{DR,u}\simeq H^1_{DR}(\mathbb{A}^1_\C, N_{DR}^{MHM}\otimes \op{exp}(u^{-1}x))=\op{Coker}\left(d/dx+u^{-1}:
N_{DR}^{MHM}\to N_{DR}^{MHM}\right).$$
Varying point $u$ we obtain an algebraic vector bundle of realization functors $DR_u$
 with connection over the punctured line of $u\in \C-\{0\}$. The connection has exponential
type at $u\to 0$ and a regular singularity at $u\to \infty$.
It is easy to see that for ordinary mixed Hodge structures (considered as objects of $EMHS$)
the fiber functor $DR_u$ is canonically isomorphic to the usual de Rham realization.

\begin{dfn} For an object $N\in EMHS$ and a non-zero complex number $u$ we define {\bf Betti}
 realization of $N$ at point $u$ to be the space of flat sections on the ray $[R,+\infty)\cdot u$
 for sufficiently large $R\gg 0$ of the constructible sheaf $N_{Betti}^{MHM}[-1]$. Here $N_{Betti}^{MHM}$ is the object of $\op{Perv}(\mathbb{A}^1(\C),\Q)$ underlying $N$, where $N$ is considered as an object of $MHM_{\mathbb{A}^1_\C}$. We denote this fiber by $N_{Betti,u}$.
\end{dfn}

Varying point $u$ we obtain a local system of realization functors $Betti_u$ in vector spaces over $\Q$. Again, for usual mixed Hodge structures considered as objects of $EMHS$, realization
$Betti_u$ has trivial monodromy in parameter $u$ and is canonically isomorphic to the usual Betti realization.

The canonical isomorphism of fiber functors mentioned above gives a {\bf comparison} isomorphism
$$N_{DR,u}\simeq N_{Betti,u}\otimes_\Q \C$$
extending the standard comparison isomorphism for usual mixed Hodge structures.

Finally, following an idea of M.~Saito,
one can define an analog of Dolbeault (or Hodge) realization for an object of $EMHS$.
 This realization can be considered as a ``limit'' at $u\to 0$ of de Rham realizations $DR_u$.

Recall that for any object $M\in MHM_X$ for any smooth algebraic variety $X/\C$ we have a canonical
 admissible increasing filtration $F_i M_{DR}$ of the underlying {\it algebraic} holonomic $\D$-module $M_{DR}$.
In particular, for $X=\AC$ the admissibility of filtration means that
 $$(d/dx) F_i M_{DR}=F_{i+1} M_{DR}$$
for all $i\gg 0$.

\begin{dfn} For an object $N\in EMHS$ we define complex vector space $N_{DR,0}$ to be the quotient
$$(d/dx)^{-j} F_j N_{DR}^{MHM} / (d/dx)^{-j} F_{j-1} N_{DR}^{MHM}$$
for sufficiently large $j\gg 0$.
\end{dfn}

Here we use the fact that the operator $d/dx$ is invertible for $N\in EMHS$. Equivalently,
 $N_{DR,0}$ is canonically isomorphic (by applying power of $d/dx$) to the quotient
$F_{j}N_{DR}^{MHM}/F_{j-1}N_{DR}^{MHM}$ for large $j\gg 0$.

For the usual mixed Hodge structures the realization $DR_0$ is canonically isomorphic
to the usual Dolbeault/Hogde realization, i.e. to the functor
$$M\in MHM_{\op{pt}}\mapsto \oplus_i \,\op{gr}^i_{\mathcal F} M_{DR}\,,$$
where  $\mathcal{F}$ is the usual (decreasing) Hodge filtration.

 The following result is not difficult to prove. Since we will not use it, the proof is omitted.
\begin{thm} Functor $DR_0$ from $EMHS$ to the category of $\C$-vector spaces is exact, faithful and commutes with the tensor product.
\end{thm}

One can see easily that functors $DR_u$ for $u\ne 0$ and $DR_0$ glue together to an
{\it algebraic vector bundle} of fiber functors on $EMHS$ over the affine line $\op{Spec}\C[u]$, with flat connection having pole of order $2$ at $u=0$.
Hence, we obtain a so-called {\it nc Hodge structure} in the sense of \cite{KaKoPa},
 i.e. a germ of a holomorphic vector bundle near $0\in \C$ together with a meromorphic connection
 with second order pole at $0$, and a rational structure on the corresponding local system on the circle
 $S^1$ of directions about the point $0$.
 Moreover, all terms of the Stokes filtration are rational vector subspaces with respect to the $\Q$-structure.

We should warn the reader that there are natural examples of nc Hodge structures appearing in
 mirror symmetry and in singularity theory which {\it do  not} come from $EMHS$.
 In particular, for a general (not quasi-homogeneous) isolated singularity of a holomorphic function
 the corresponding nc Hodge structure does not come from an exponential one. For example, it is easy to construct a non-trivial family
 of isolated singularities with constant Milnor number over $\AC$, giving a non-trivial variation of non-commutative Hodge structures\footnote{See e.g. section 3.3 in 
the third part of 
\cite{Kulikov}, on the period map for $E_{12}$ singularities. In this book  
the nc Hodge structure for an isolated singularity is encoded by the equivalent notion of a Brieskorn lattice.}.
 On other side, any variation of $EMHS$ over $\C^1$ (with tame singularity at infinity) is trivial.
 Similarly, for the nc Hodge structure associated with a non-K\"ahler compact symplectic manifold
 with symplectic form rescaled by a sufficiently large factor, there is no reasons  to expect that there exists an underlying object of $EMHS$.

\subsection{Weight filtration in the exponential case}

Let us denote by $i$ the natural inclusion functor $EMHS\mono MHM_{\AC}$ and by $p:MHM_{\AC}\to EMHS$ its left adjoint.
Both functors are exact, and
$$p\circ i=\op{id}_{EMHS}\,,\,\,\,\,i\circ p=\Pi\,,$$
where we use the same notation $\Pi$ (as for categories $\mathfrak A\simeq \op{Perv}(\AC)$) for
 an exact idempotent functor of $MHM_{\AC}$ given  by
$$\Pi:M \mapsto M*_+ j_!(\Q(0)_{\AC}[1])\,.$$

\begin{dfn} The {\bf weight } filtration of an exponential mixed Hodge structure is defined by the formula
$$W_{\leqslant n}^{\op{EMHS}} M=W_{n}^{\op{EMHS}} M:= p(W_n(i(M)))\,,$$
where in the RHS we consider the usual weight filtration for mixed Hodge modules.
\end{dfn}

\begin{prp} The endofunctor $\op{gr}_\bullet^{W,EMHS}$ of $EMHS$ obtained by taking the associated graded object with respect to the filtration $W_\bullet^{EMHS}$
 is exact and faithful.
\end{prp}

{\it Proof.} The exactness follows from exactness of functors $i,p$, and of the endofunctor $\op{gr}_\bullet^{W}$ of $MHM_{\AC}$.
 To prove the faithfulness we use the existence of the fiber functor on the tensor category $EMHS$. Namely,  any object  $M\in EMHS$ has
 $\op{rk}(M)\geqslant 0$, and $M=0$ if and only if $\op{rk}(M)=0$. Endofunctor $\op{gr}_\bullet^{W,EMHS}$ preserves the rank because
$$ W_n^{\op{EMHS}} M =0\mbox{ for }n\ll 0\,,\,\,\,\, W_n^{\op{EMHS}} M =M\mbox{ for }n\gg 0\,.    $$
Let $\alpha:M\to N$ be a morphism such that $\op{gr}_\bullet^{W,EMHS}(\alpha)=0$. The exactness of $\op{gr}_\bullet^{W,EMHS}$ implies that
  $$   \op{gr}_\bullet^{W,EMHS}(\op{Ker}\alpha)=\op{Ker}(\op{gr}_\bullet^{W,EMHS}\alpha)=\op{gr}_\bullet^{W,EMHS}(M)\,.$$
Hence the rank of $\op{Ker}\alpha$ is equal to the rank of $M$, and we conclude that $\op{rk}( M/\op{Ker}\alpha)=0$. Therefore $M=\op{Ker}\alpha$ and $\alpha=0$.
$\blacksquare$.

\begin{dfn} An object $M\in EMHS$ is called {\bf pure} of weight $n\in \Z$ if $W^{EMHS}_n M=M$ and $W^{EMHS}_{n-1}M=0$.
\end{dfn}

Graded factors of the weight filtration of an object in $EMHS$ are pure.
The next proposition describes the structure of pure objects in $EMHS$.

\begin{prp} The full subcategory of direct sums of pure objects in $EMHS$ is a semisimple abelian category. Moreover, it is equivalent
to a direct summand of the semisimple category of pure objects in $MHM_{\AC}$, with the complement consisting
 of constant variations of pure Hodge structures on $\mathbb{A}^1(\C)$.
\end{prp}
{\it Proof.}
Obviously, the image of any pure $F\in MHM_{\AC}$ of weight $n\in \Z$ by projection $p$ is a pure object of $EMHS$ of the same weight.
 Conversely, any pure object $E\in EMHS$ of weight $n$ is the image of a pure object of weight $n$ in $MHM_{\AC}$. Indeed,
we have
$$p(W_n(i(E)))=E,\,\,\,p(W_{n-1} (i(E)))=0$$
and by exactness of $p$ we conclude that
$$E=p(W_n(i(E))/W_{n-1}(i(E)))=p(\op{gr}^W_n (i(E)))\,.$$
Also, it is easy to see that the kernel of the projector  $\Pi$ on $\op{Perv}(\mathbb{A}^1(\C),\Q)$ consists of constant sheaves.
Any pure object of $MHM_{\AC}$ whose underlying perverse sheaf is constant is a constant variation of pure Hodge structure.
Hence, we conclude that the set of isomorphism classes of pure objects of weight $n$ in $EMHS$ is in one-to-one correspondence with pure objects
 of weight $n$  in $MHM_{\AC}$ which have no direct summands which are constant variations.

For any pure object $F$ of weight $n$ in $EMHS$ without constant direct summand the corresponding
 pure object $E=p(F)$ has the following structure. The corresponding object $i(E)=\Pi(F)$
 has two-step weight filtration:
$$\Pi(F)=W_n \Pi(F)\supset W_{n-1} \Pi(F)\supset W_{n-2} \Pi(F)=0$$ with
$\op{gr}_n^W \Pi(F)=F$  and $\op{gr}_{n-1}^W \Pi(F)$ being a constant variation
of Hodge structures of weight $n-1$.

We have to check for any two pure objects $E_1,E_2\in EMHS$ of weights $n_1,n_2\in \Z$ the following isomorphism
$$ \op{Hom}(E_1,E_2)\simeq\op{Hom}(\op{gr}^W_{n_1}(i(E_1)),\op{gr}^W_{n_2}(i(E_2)))\,.$$
In the case $n_1\ne n_2$ the left hand side vanishes by the faithfulness of the functor $\op{gr}_\bullet^{W,EMHS}$, and the right hand side
 vanishes by the faithfulness of $\op{gr}_\bullet^{W}$.

Consider now the case of equal weights $n_1,n_2$ and denote $F_j:=\op{gr}_n^{W}(i(E_j)),\,\,j=1,2$ the corresponding pure objects of $MHM_{\AC}$.
 Adjunction morphisms
 $$F_i\to \Pi(F_i),\,\,\,\Pi=i\circ p$$
induce (after composition with $\op{gr}_\bullet^{W}$) isomorphisms $F_j\simeq \op{gr}_n^{W}(\Pi(F_j))$.
We have
$$\op{Hom}(E_1,E_2)\simeq\op{Hom}(p(F_1),p(F_2))\simeq \op{Hom}(F_1, \Pi(F_2))\mono$$
$$\mono\op{Hom}(\op{gr}_\bullet^{W}(F_1),\op{gr}_\bullet^{W}(\Pi F_2))\simeq\op{Hom}
(F_1,\op{gr}_n^W i(E_2))\simeq\op{Hom}(F_1,F_2)$$
by faithfulness of the functor $\op{gr}_\bullet^{W}$. Hence $\op{rk}\op{Hom}(E_1,E_2)\leqslant \op{rk}\op{Hom}(F_1,F_2)$.
On the other hand, we have natural transformations of functors
$$\op{id}_{MHM_{\AC}}\to \Pi\to \Pi/(W_{n-1}\circ \Pi)$$
which induce isomorphisms
$$F_j\to \Pi(F_j)\to \Pi(F_j)/W_{n-1}\Pi(F_j)\simeq\op{gr}_n^{W}(\Pi(F_j))\simeq F_j\,.$$
Hence $\op{Hom}(F_1,F_2)$ is a retract of $\op{Hom}(\Pi(F_1),\Pi(F_2))\simeq \op{Hom}(E_1,E_2)$.
We conclude that $\op{Hom}(F_1,F_2)\simeq\op{Hom}(E_1,E_2)$.
$\blacksquare$

\begin{prp}\label{lmm:tensor-weight}
The weight filtration on $EMHS$ is strictly compatible with the tensor product.
 \end{prp}
{\it Proof.} We start with the following lemma.

\begin{lmm} Let $E,F\in D^b(MHM_{\mathbb{A}^1_{\C}})$ and
$$G=\op{Cone}(sum_{!}(E\boxtimes F)\to sum_*(E\boxtimes F))\,.$$
 Then $\Pi(G)=0$.
\end{lmm}
{\it Proof of the Lemma.} We first observe that we can assume that $E$ and $F$ are just constructible sheaves (since the forgetful functor to constructible sheaves is faithful). We will show that $G$ is a  constant sheaf. The fiber of $G$ at $x\in \mathbb{A}^1(\C)$ is canonically isomorphic to $H^{\bullet}(S^1_{x,R},E\boxtimes F)$, where
$S^1_{x,R}\subset sum^{-1}(x)$ is a circle of a sufficiently large radius $R$ with any given center. Obviously the cohomology groups do not depend on $x$. Hence they form a constant sheaf on $\mathbb{A}^1(\C)$. In order to finish the proof of the lemma we observe that the functor $\Pi$ kills all constant sheaves on $\mathbb{A}^1(\C)$. $\blacksquare$

Let us return to the proof of the Proposition. It is enough to prove that  pure objects
 in $EMHS$ are closed under the tensor product,  and weights behave additively.

Assume that $E_j\in EMHS, j=1,2$ are objects of pure weights $n_k,k=1,2$. Then $E_j=p(F_j)$ for some pure objects $F_j$ in $MHM_{\AC},\,\,j=1,2$. In the category $EMHS$ we have:
$$E_1\otimes_{EMHS}E_2=p(sum_*(F_1\boxtimes F_2))=p(sum_!(F_1\boxtimes F_2))\,,$$
where the second equality follows from the Lemma.
Since both objects $F_1$ and $F_2$ are pure, we see that $F_1\boxtimes F_2$ is pure. It is well-known (see e.g. \cite{Steen}) that $sum_*(F_1\boxtimes F_2)$ has weights (in the derived sense) greater or equal than $n_1+n_2$ while $sum_!(F_1\boxtimes F_2)$ has weights less or equal than $n_1+n_2$.

Let us denote by $H_{perv}^0(F)$ the zeroth middle perverse cohomology of an object $F\in D^b(MHM_{\mathbb{A}^1_{\C}})$. Then
$$p\circ sum_*(F_1\boxtimes F_2)=p\circ (H_{perv}^0(sum_*(F_1\boxtimes F_2))\,,$$
and similarly for  $sum_!(F_1\boxtimes F_2)$ (recall that the functor $p$ is exact).
Therefore
$$E_1\otimes_{EMHS}E_2=p(H_{perv}^0(sum_*(F_1\boxtimes F_2)))=
p (H_{perv}^0(sum_!(F_1\boxtimes F_2)))\,.$$
We see that $E_1\otimes_{EMHS}E_2$ is the image by $p$ of an object of $MHM_{\AC}$ of the weight
 $\geqslant (n_1+n_2)$. But it is also the image by $p$ of an object of the weight $\leqslant (n_1+n_2)$.
Notice that the exact functor $p$ maps weight filtration in $MHM_{\AC}$ to the weight filtration in
$EMHS$, as follows from our description of pure objects in $EMHS$.
Therefore $E_1\otimes_{EMHS}E_2$ is pure of weight $n_1+n_2$. This concludes the proof of the Proposition. $\blacksquare$

Let us define {\bf Serre polynomial} of an exponential mixed Hodge structure $E$ as
$$S(E)=\sum_i rk \op{gr}_i^{{exp}}(E)q^{i/2}\in \Z[q^{\pm{1/2}}]\,,$$
where  $\op{gr}_{\bullet}^{{exp}}$ denotes the weight filtration on $EMHS$.
By additivity we extend $S$ to a functional on the $K_0$-group of the bounded derived category of $EMHS$. By Proposition \ref{lmm:tensor-weight} the map
$S$ is a ring homomorphism.

\subsection{Cohomology of a variety with function}

For any separated scheme of finite type $X/\C$ and a function $f\in \mathcal{O}(X)$ (sometimes called {\it potential of $X$}) we
 define {\bf exponential cohomology} (which will be an exponential mixed Hodge structure) $H^i_{EMHS}(X,f)$ for integers $i \geqslant 0$ in the following way:
$$H^i_{EMHS}(X,f):=\Pi(H^i_{MHM_{\AC}} f_* \Q_X) =H^i_{MHM_{\AC}} (\Pi f_* \Q_X)\in EMHS\,,$$
where we interchange $\Pi$ and cohomology in $D^b(MHM_{\AC})$ using the exactness of the functor $\Pi$.

If $f$ is the restriction of a function $f'$ defined on a bigger variety $X'\supset X$ then we will abuse the notation denoting both functions by the same symbol. The rapid decay cohomology of the pair $(X,f)$ does not depend on such an extension. We will see in  Section 7 that for the so-called critical cohomology (defined in terms of the vanishing cycles functor) the situation is different, and the result depends on the extension of the function.

As follows directly from definitions, the Betti realization of $H^i_{EMHS}(X,f)$ is given by the
cohomology of pair
$$ H^i_{{Betti},u}(X,f):=H^i_{{Betti}}(X(\C),f^{-1} (-c\cdot u);\Q)\,,\,\,\,H^i_{{Betti}}(X,f):=H^i_{{Betti},1}(X,f)\,,$$
where $c\gg 0$ is a sufficiently large positive real number. Cohomology groups are identified for different large values of $c$ via
the holonomy along the
 ray $\R_{<0}\cdot u$.

We  define Betti cohomology with {\it integer} coefficients as
$$H^i_u(X,f):=H^i(X(\C),f^{-1} (-c\cdot u);\Z),\,\,\,H^i_{{Betti},u}(X,f)\simeq H^i_u(X,f)\otimes_\Z \Q\,.$$
The cohomology $H^\bullet(X,f):=H^\bullet_{1}(X,f)$ behaves similarly to the usual cohomology. In particular, for a map $\pi:Y\to X$ compatible with functions
 $f_Y\in \mathcal{O}(Y),\,f_X\in \mathcal{O}(X)$ in the sense that $f_Y=\pi^* f_X$, we have the pullback
$\pi^*:H^\bullet(X,f_X)\to H^\bullet(Y,f_Y)$. When $\pi$ is proper and both $X$ and $Y$ are smooth we have the pushforward morphism
$\pi_*:H^\bullet(Y,f_Y)\to H^{\bullet+2(\dim_\C(X)-\dim_\C(Y))}(X,f_X)$. All these morphisms induce morphisms of exponential mixed Hodge structures
 (for the pushforward we should twist by a power of the Tate module $\mathbb T$). The proof is a straightforward application of the standard formalism of six functors. 

As for usual cohomology, there is a  multiplication  morphism
  $$\otimes:H^\bullet(X,{f_X})\otimes H^\bullet(Y,{f_Y})\to H^\bullet(X\times Y, {f_X\boxplus f_Y})\,\,,$$
where the Thom-Sebastiani sum $\boxplus$ is given by
$$ f_X\boxplus f_Y:=\mathrm{pr}_{X\times Y\to X}^* f_X+  \mathrm{pr}_{X\times Y\to Y}^* f_Y \,\,.$$
This morphism is well-defined because when the real part of the sum of two functions tends to $-\infty$, then the real part of at least one of the functions also
tends to $-\infty$.
The multiplication morphism becomes an isomorphism after the extension of coefficients for cohomology from $\Z$ to $\Q$. The multiplication morphism gives an isomorphism of exponential mixed Hodge structures
$$H^\bullet_{EMHS}(X,{f_X})\otimes H^\bullet_{EMHS}(Y,{f_Y})\stackrel{\sim}{\to} H^\bullet_{EMHS}(X\times Y, {f_X\boxplus f_Y})\,\,.$$

For smooth $X$ the de Rham realization of $H^i_{EMHS}(X,f)$ (at a point $u\in \C^*$) is given by the
 de Rham cohomology of $X$ with coefficients in the holonomic $\mathcal{D}_X$-module $\op{exp}\left(u^{-1}\cdot f\right)\cdot \mathcal{O}_X$. In other words
it is a finite-dimensional $\Z$-graded vector space over $\kk$ which is the hypercohomology in Zariski topology
 $$H^\bullet_{DR,u}(X,f):=\mathbb{H}^\bullet\left(X_{Zar},(\Omega^\bullet_X,d+u^{-1}\,df\wedge\cdot)\right),\,\,H^\bullet_{DR}(X,f):=
H^\bullet_{DR,1}(X,f).$$

The abstract comparison isomorphism between complexified Betti and de Rham realization (see Section 4.3, and also Theorem 1.1 in \cite{Sabbah} and references therein) is
  $$H^\bullet_{DR}(X,f)\simeq \C\otimes H^\bullet(X,f)\,\,.$$
  In the affine case it is given by the
  integration of complex-analytic closed forms on $X(\C)$ of the type $\op{exp}(f)\,\alpha$, where $\alpha$ is an algebraic form on $X$ such that
$$d\alpha+df\wedge \alpha=0\,\Longleftrightarrow\, d\left(\op{exp}(f)\,\alpha\right)=0\,\,,$$
 over closed real semi-algebraic chains with the ``boundary at infinity'' in the direction $\op{Re}(f)\to -\infty$. The integral
   is absolutely convergent because the form $\op{exp}(f)\,\alpha$ decays rapidly at infinity. This explains the term ``rapid decay cohomology".

 An important example of a (graded) $EMHS$ is
$$\mathbb{T}^{\otimes 1/2}:=H^\bullet(\AC, -z^2)\,,$$
which is an odd one-dimensional space in degree $+1$ corresponding to the period integral
$\int_{-\infty}^{+\infty} \op{exp}(-z^2) dz= \sqrt{\pi}$.
This structure is a ``square root'' of the shifted Tate structure $\mathbb{T}=\Q(-1)[-2]$, and is a pure $EMHS$ of weight $+1$ and degree $+1$.

Notice that de Rham cohomology can be defined for a pair $(X,f)$ where $X$ is a smooth scheme over {\it arbitrary} field of characteristic zero,
 and $f\in \mathcal{O}(X)$. If $(X_0,f_0)$ is defined over $\Q$ then the comparison isomorphism (say, at the point $u=1$)  gives {\it two} rational structures on the
same complex vector space. These two structures are related by a period matrix whose coefficients are exponential periods in the sense of Section 4.3 in \cite{KZ}.

One can also define the exponential cohomology {\bf with compact support} by
$$H^\bullet_{c,EMHS}(X,f):=\Pi(H^{\bullet}_{MHM{\AC}} f_! \Q_X)\,.$$
They satisfy the following Poincar\'e duality for  smooth $X$ of dimension $d$:
$$H^\bullet_{c,EMHS}(X,f)^\vee\simeq H^\bullet_{EMHS}(X,-f)\otimes \Q(d)[2d]\,.$$

Also, one can define exponential cohomology of {\bf pairs} by
$$H^\bullet_{EMHS}(X,Y,f):=\Pi(H^{\bullet}_{MHM_{\AC}} f_*\op{Cone}( \Q_X\to (Y\to X)_*\Q_Y))[-1]\,.$$
For a smooth closed $Z\subset X$ of codimension $d'$ in smooth $X$ we have the Thom isomorphism
$$H^\bullet_{EMHS}(X,X-Z,f)\simeq H^\bullet_{EMHS}(Z,f)\otimes \Q(d')[2d']\,.$$

If affine algebraic group $G$ acts on a scheme of finite type $X$ preserving $f\in \mathcal{O}(X)$, we define {\bf equivariant} exponential cohomology as
$$H^\bullet_{G,EMHS}(X,f):=H^\bullet_{EMHS}(X^\univ,f^\univ)\,,$$
where $X^\univ$ is the universal bundle over $\op{B}G$ endowed with the induced function $f^\univ$. This cohomology depends only on the quotient stack $X/G$
and on $f\in \mathcal{O}(X/G)$, so we will sometimes denote it by $H^\bullet_{EMHS}(X/G,f)$.

We will also need exponential cohomology with compact support which can be defined for a quotient Artin stack endowed with a function in the following way. We can always assume that our stack is represented as $X/G$
 where $G=\op{GL}(n,\C)$ for some $n$. The standard classifying space $\op{B}G=\varinjlim(\op{B}{G})_N$ (where $N\to \infty$) is the union of finite-dimensional  Grassmannians $(\op{B}G)_N:=\op{Gr}(n,\C^N)$ which are compact smooth varieties. Then we have an inductive system $X_{N}\to (\op{B}G)_{N}$
 of locally trivial bundles endowed with functions $f_N$.  Then we define the exponential {\bf cohomology with compact support of the stack} $X/G$ with function $f$ as
$${H}_{c,EMHS}^{\bullet}(X/G,f):=\varinjlim_{N} {H}_{c,EMHS}^{\bullet}(X_N,f_N)\otimes {\mathbb T}^{\otimes -\dim((\op{B}G)_{N}\times G)}\,,$$
where connecting maps are Gysin morphisms. Cohomology  groups ${H}_{c,EMHS}^{\bullet}(X/G,f)$ stabilize in each degree and are concentrated in cohomological degrees bounded from above.  For smooth equidimensional $X$ cohomology groups ${H}_{EMHS}^{\bullet}(X/G,f)$ and ${H}_{c,EMHS}^{\bullet}(X/G,-f)$
 satisfy Poincar\'e duality in dimension $\dim(X/G):=\dim(X)-\dim(G)$ \footnote{Cohomology of Artin stacks are discussed in different settings and approaches in e.g.\cite{Behr},  \cite{BerLu}, \cite{LasOl}. We use here a more direct approach which is sufficient for our purposes.}. In what follows we will often omit the word ``exponential'' if it is clear that we are talking about exponential cohomology.

\begin{dfn} A $\Z$-graded object $E^\bullet=\oplus_{n\in \Z} E^n$ of $EMHS$ is called {\bf pure}  if for any $n\in \Z$ component $E^n$ is pure of weight $n$.
\end{dfn}

Similarly to the classical situation without potential we have
\begin{prp} For a pair $(X,f)$ over $\C$ such that $X$ is smooth and
 the natural morphism $H^\bullet_{c,EMHS}(X,f)\to H^\bullet_{EMHS}(X,f)$ is an isomorphism, the cohomology $H^\bullet_{c,EMHS}(X,f)$ is pure.
\end{prp}
{\it Proof.} This follows from the fact that $f_*\Q_X\in D^b(MHM(\AC))$ has  weights $\geqslant 0$ (in the derived sense), while $f_!\Q_X$
 has weights $\leqslant 0$. $\blacksquare$

In particular, the conditions of the Proposition hold
for a pair $(X,f)$ such that $X$ is smooth and $f:X\to \AC$ is {\it proper}. More generally, $f$ could be cohomologically tame (see e.g. \cite{Dimca}, \cite{Sabbah2}),
 e.g. when $X$ is an affine space and $f$ is a quasi-homogeneous polynomial with isolated singularity. For example,
  the graded $EMHS$ given by $\mathbb{T}^{\otimes 1/2}$ is pure.

If $H^\bullet_{EMHS}(X,f)$ is pure then the Serre polynomial
$$S(H^\bullet_{EMHS}(X,f))=\sum_{n\in \Z} (-1)^n S(H^n_{EMHS}(X,f))\in \Z[q^{1/2}]$$
coincides with the Poincar\'e polynomial of $H^i_{{Betti}}(X(\C),f^{-1} (c\cdot u);\Q)$ in variable $(-q^{1/2})$, and the evaluation at $q^{1/2}=1$ coincides with the Euler characteristic
$$\chi(X,f):=\chi(H^\bullet_{EMHS}(X,f))\,.$$

Using the purity of  $H^\bullet_{EMHS}(X,f)$ for smooth $X$ and proper $f$, we can obtain an elementary 
topological description of the weight filtration in the more general case when $X$ is smooth but $f$ is not necessarily proper.
Namely, for any smooth $X/\C$ and $f\in \mathcal{O}(X)$ there exist smooth $\overline{X}\supset X$ and $f'\in \mathcal{O}(\overline{X})$
such that $f'_{|X}=f$, the map $f':\overline{X}\to \AC$ is proper, and the complement $D=\overline{X}-X$ is a divisor with normal crossings.
Then the term $W_k$ of the weight filtration on $H^i_{EMHS}(X,f)$ is given by the image of $H^i_{EMHS}(X^{(j)},f'_{|X^{(j)}})$ where $X^{(j)}$ is
the complement in $\overline{X}$ to the union of $j$-dimensional strata in $D$, for $j=k-i-1$. The same is true for $H^\bullet_{DR}$ and $H^\bullet_{Betti}$.

Also, similarly to the case of usual mixed Hodge structures, one can describe the associated graded space with respect to the weight filtration entirely in terms of pure structures.  For that purpose it is more convenient to use cohomology with compact support. Let us make a simplifying assumption that $D=\overline{X}-X$ is the union of smooth divisors $D_1,\dots ,D_N$
 intersecting transversally. Then $\op{gr}^W_k {H}^l_{c,Betti}(X,f)$ 
 is isomorphic to the cohomology group in degree $n=l-k$ of the complex
$$\bigoplus_{0\leqslant n\leqslant N}\,\, \bigoplus_{1\le i_1<\dots<i_n\leqslant N}{H}^k_{Betti}(D_{i_1}\cap \dots\cap D_{i_n},f)$$
with the differential given by the restriction morphisms.

\subsection{\'Etale case}

Let $\kk$ be any field and $l$ be a prime, $l\ne \op{char}(\kk)$.
With any scheme of finite type $X/\kk$ and $f\in \mathcal{O}(X)$ we associate
perverse $l$-adic sheaves on $\mathbb{A}_\kk^1$  by the formula
$$H^i_{l}(X,f):=\Pi(H^i_{\op{Perv}_l\mathbb{A}_\kk^1}f_*\Q_{l,X})\,,$$
where the functor $\Pi$ on the category $ \op{Perv}_l\mathbb{A}_\kk^1  $ of $l$-adic perverse sheaves on $\mathbb{A}_\kk^1$ is defined similarly to the constructible case
and $\kk=\C$. These sheaves are killed by the functor $\op{R}\Gamma$.

If $\op{char}\kk=0$ then by the comparison isomorphism we obtain just another realization
 of the ``exponential motive'' associated with $(X,f)$. All previous considerations
concerning the weight filtration hold in $l$-adic case as well. In the case $\kk=\Q$ one can define an exponential motivic Galois group (see Chapter 4 and in particular Section 4.3  in \cite{KZ}) whose representations have both $l$-adic and $ EMHS$ realizations.

The situation is more complicated if $p:=\op{char}\kk$ is non-zero.
 In this case one can apply Fourier transform using the standard exponential sheaf
depending on the choice of a primitive $p$-th root of unity in $\overline{\Q}_l$.

The restriction of $FT(f_*\Q_{l,X})$ to $\mathbb{A}^1_\kk-\{0\}$ coincides with
$FT(\Pi(f_*\Q_{l,X}))$. Perverse sheaves $FT(H^i_{l}(X,f))$ are in general not {\it lisse}.

Picking any non-zero point $u$ in $\mathbb{A}^1(\kk)$ (e.g. point $u=1$)
one obtains a cohomology theory for pairs $(X,f)$ with values in $\overline{\Q}_l$-linear
representations of $\op{Gal}(\kk^{sep}/\kk)$, by taking $H^i$ of the fiber of
$FT(f_*\Q_{l,X})$ at $u$. For this cohomology theory (which one can denote by $H^\bullet_{l,u}(X,f)$)
one can define the weight filtration, Serre polynomial, and also a multiplicative $\overline{\Q}_l$-valued
numerical invariant (in fact it is a $\Q(\exp(2\pi i/p)$-valued invariant) obtained
 by taking the supertrace of Frobenius.

For $X=\mathbb{A}^1_\kk$ and any $\lambda\in \kk^*$ the exponential $l$-adic structure
$\mathbb{\otimes T}^{1/2}_{\lambda,l}$ corresponding to $(\mathbb{A}^1_\kk,\lambda z^2)$, where $z$ is a coordinate on $\mathbb{A}^1_\kk$,
is a tensor square root of $\mathbb{T}=\mathbb{T}_l=\Q_l(-1)[-2]$, if and only if $\kk$ contains $\sqrt{-1}$.
 Different choices of $\lambda$ give in general non-isomorphic objects which differ by quadratic characters of $\op{Gal}(\kk^{sep}/\kk)$.

\subsection{Cohomological Hall algebra in the case of a smooth algebra  with potential}

Let $R$ be a smooth $I$-bigraded algebra over field $\kk$ endowed with a bilinear form $\chi_R$ on $\Z^I$ compatible with Euler characteristic
 as in Section 3.3.

Let us assume that we are given an element
$$W\in R/[R,R]$$
represented by some element $\widetilde{W}\in R,\,\,\,W=\widetilde{W}\pmod{[R,R]}$. The element $W$
(or its lifting $\widetilde{W}$) is called a {\it potential} \footnote{We also use sometimes the same word ``potential'' for a function on an algebraic variety when we consider exponential Hodge structures. We hope it will not lead to a confusion.}. 
 Then for any $\gamma\in \Z_{\geqslant 0}^I$ we obtain a function $W_\gamma$  on the  affine variety $\M_\gamma$, invariant under the action of $\G_\gamma$.
The value of $W_\gamma$ at any representation is given by the trace of the image of $\widetilde{W}$.
For any short exact sequence
$$0\to E_1\mono E\epi E_2\to 0$$ of representations of $R$
we have $W_{\gamma_1+\gamma_2}(E)=W_{\gamma_1}(E_1)+W_{\gamma_2}(E_2)$, where $\gamma_i,\,i=1,2$ are dimension vectors of $E_i,\,i=1,2$.

Let us denote by $\bf H$ one of the cohomology theories considered above (i.e. $EMHS$, Betti, de Rham or \'etale cohomology), with values in an appropriate Tannakian category $\mathcal C$.

\begin{dfn} The {\bf Cohomological Hall algebra} of $(R,W)$ (in realization $\bf H$) is an associative twisted graded  algebra in $\mathcal C$
defined by the formula
$$\mathcal{H}:=\oplus_\gamma \mathcal{H}_\gamma\,,\,\,\,\,\mathcal{H}_\gamma:={\bf H}^\bullet(\M_\gamma/\G_\gamma,W_\gamma):=
{\bf H}^\bullet(\M_\gamma^\univ,W_\gamma^\univ)\,\,\,\,\forall \gamma \in \Z_{\geqslant 0}^I$$
in the obvious notation. (The definition of the product and the proof of associativity are completely similar to the case of an algebra without potential).
\end{dfn}

More precisely, all the morphisms (the K\"unneth isomorphism, the pullback and the Gysin map) which appear  in the definition of the product in  Section 2.2 induce morphisms in realization $\bf H$.  In the case ${\bf k} \subset\C$ the commutativity of  nine small  diagrams in Section 2.3 follows from the corresponding commutativity
 in Betti realization and from the comparison isomorphism. Alternatively, one can use the formalism of six functors for general $\bf k$.

As an example, we describe the algebra $\mathcal{H}$ in the case of the quiver $Q=Q_1$ with one vertex and one loop
 $l$,
and potential $W=\sum_{i=0}^{N} c_i l^i,\,\,\,c_N\ne 0$, given by a polynomial of degree $N\in \Z_{\geqslant 0}$ in one variable.
Dimension vector $\gamma$ for such $Q$ is given by an integer $n\geqslant 0$. For simplicity, we consider Betti realization.

In the case $N=0$, the question reduces to $Q$ without potential, and was considered above. Algebra $\mathcal{H}$ is the polynomial algebra of infinitely many variables (see Section 2.5).

In the case $N=1$ the cohomology of pair vanishes for matrices of size greater than $0$, hence $\mathcal{H}=\mathcal{H}_0=\Z$.

In the case $N=2$ we may assume without loss of generality that $W=-l^2$. It is easy to see that the cohomology of
  $$(\op{Mat}(n\times n, \C)\,,\,\,\, x \mapsto -\op{Tr}(x^2))$$
can be identified
under the restriction map with the cohomology 
  $$(\op{Herm}(n)\,, \,\,\,x\mapsto -\op{Tr}(x^2))\,\,,$$
where $\op{Herm}(n)$ is the space of Hermitean matrices. The latter cohomology group is same as $H^*(\op{D}^{n^2},\partial\op{D}^{n^2})$,
where $\op{D}^{n^2}$ is the standard closed unit ball in $\R^{n^2}\simeq \op{Herm}(n)$.
 Hence $H^\bullet\left(\op{Mat}(n\times n, \C), \,\,\,x \mapsto -\op{Tr}(x^2)\right)$
 is isomorphic to $\Z$ concentrated in the cohomological degree $n^2$.
 Moreover, one can use the unitary group $U(n)$ instead of the homotopy equivalent group $\mathrm{GL}(n,\C)$ in the definition of equivariant cohomology. Group   $U(n)$ acts on the pair $(\op{D}^{n^2},\partial\op{D}^{n^2})$, and Thom isomorphism gives  a canonical isomorphism of cohomology groups
    $$\mathcal{H}_n\simeq H^\bullet(\op{BGL}(n,\C))[-n^2]\simeq H^\bullet_{U(n)}(\op{D}^{n^2},\partial\op{D}^{n^2})\,\,.$$
   Let us endow $\mathcal{H}$ with the natural bigrading which was used in the case of a symmetric incidence matrix in Section 2.6.
 Then an easy calculation shows  that $\mathcal{H}$ coincides as a bigraded abelian group with the algebra associated with the quiver $Q_0$ with one vertex and zero arrows.
Furthermore, comparing Grassmannians which appear in the definition of the multiplication for the quiver $Q_0$ (see Section 2.2), with those which arise for $Q_1$ with potential $W=-l^2$, one can check that the multiplications coincide as well. Hence the algebra $\mathcal{H}=\mathcal{H}_{(Q,W)}$ is the exterior algebra with infinitely many generators.

In the case of degree $N\geqslant 3$ one can show that the bigraded algebra $\mathcal{H}$
is isomorphic to the $(N-1)$-st tensor power of the  exterior algebra corresponding to the case $N=2$.
The proof will be given in Section 7, using so-called critical cohomology theory. In this particular case the critical cohomology coincides with the rapid decay one.

Let us assume that for the Tannakian category $\mathcal C$ (target of the cohomology functor) we have the notion of weight filtration. We will say that COHA $\mathcal{H}$ is {\it pure} if for any $\gamma$  the graded space ${\bf H}^\bullet(\M_\gamma^\univ,W_\gamma^\univ)$ is pure, i.e. its $n$-th component is of weight $n$ for any $n\in \Z$.
 In order to prove that  $\mathcal{H}$ is pure it is sufficient to check that ${\bf H}^\bullet(\M_\gamma,W_\gamma)$ is pure. Indeed,  in this case the spectral sequence
$${\bf H}^\bullet(\M_\gamma,W_\gamma)\otimes {\bf H}^\bullet(\op{B}\G_\gamma)\Longrightarrow {\bf H}^\bullet(\M_\gamma^\univ,W_\gamma^\univ)$$
collapses because ${\bf H}^\bullet(\op{B}\G_\gamma)$ is pure, and hence for every weight $n$ we have a complex supported only in degree $n$.
  Examples of pure COHA include all quivers with zero potentials as well as the above example of $Q_1$ with one vertex, one loop and an arbitrary potential.
  For $R=\C[t,t^{-1}]$ (quiver $Q_1$ with an  invertible loop), and $W=0$, the corresponding COHA is {\it not} pure.

  \begin{rmk}
a) In the case when $R$ is the path algebra of a quiver $Q$,
 the space of representations $\M_\gamma$ is the space of collections of matrices. 
Integrals of $\op{exp}(W_\gamma/u)$ over
appropriate non-compact cycles in $\M_\gamma$ (usually over the locus of Hermitean or unitary matrices) are exactly objects of study
 in the theory of matrix models in mathematical physics.
Those integrals are encoded in the comparison isomorphism between Betti and de Rham realizations of  $H^\bullet_{EHMS}(\M_\gamma/\G_\gamma)$ and can be interpreted as periods of the corresponding ``exponential motives''.

b) The non-zero constant $u$ parametrizing the comparison isomorphism corresponds to the string coupling constant $g_s$. Notice that the parameter $u$ is the same for all dimension vectors $\gamma$. Moreover, we do not have neither a distinguished integration cycle nor a volume element. As a result we do not consider
the ``large N'' (in our notation ``large $|\gamma|$'') behavior of matrix integrals.

c) S. Barannikov developed (see \cite{Baran}) a generalization of the notion of Calabi-Yau algebra in terms of BV-formalism for multi-trace functionals (i.e. products of trace functionals). It seems that in his formalism there is a distinguished volume element. The connection of his approach to our work is not clear, because multi-trace functionals do not enjoy the additivity property with respect to the direct sum of representations.

 \end{rmk}

\subsection{Example: linear potentials and cohomological dimension $2$}

We calculate an exponential motive in the following general situation.
  Let $\pi:X\to Y $ be a fibration by affine spaces between smooth equidimensional algebraic varieties over $\kk$ endowed with an action of an affine algebraic group $G$, and $f\in \mathcal{O}(X)^G$ be a function which
is affine along fibers of $\pi$ (i.e. $f$ is a polynomial of degree $\leqslant 1$ on $\pi^{-1}(y)$ for $\forall y$). Denote by $Z\subset Y$ the  closed subvariety $Z\subset Y$ of points $y\in Y$ such that $f$ is constant on $\pi^{-1}(y)$, and let $f_Z$ be the induced function on $Z$.

\begin{prp} Under the above assumptions, there is a natural isomorphism
$$\mathbf{H}^\bullet(X/G,f)\simeq \mathbf{H}^\bullet_c(Z/G,-f_Z)^\vee\otimes \mathbb{T}^{\otimes \op{dim}Y/G}\,,$$
where the cohomology of a stack is defined via the universal bundle construction as in Section 1.
\end{prp}
{\it Proof.}
Notice first that
$$\mathbf{H}^\bullet(X-\pi^{-1}(Z),f)\simeq 0$$
because $f_*(\Q_{X-\pi^{-1}(Z)})$ is a constant sheaf on $\mathbb{A}^1_\kk$, and the projector $\Pi$ used in the definition of exponential
cohomology in any realization kills constant sheaves. The same is true for equivariant cohomology.  Applying the duality we see that $\mathbf{H}^\bullet_c((X-\pi^{-1}(Z))/G,-f)$ vanishes.
 Hence (from the long exact sequence of pair $(X,X-\pi^{-1}(Z))$) we see that the direct image induces an isomorphism
$$\mathbf{H}^\bullet_c(\pi^{-1}Z/G,-f_{|\pi^{-1}(Z)})\simeq  \mathbf{H}_c^\bullet(X/G,-f)\,.$$
Applying Poincar\'e duality for $X/G$ and Thom isomorphism for the affine fibration $\pi^{-1}(Z)\to Z$, we obtain the statement of the Proposition. $ \blacksquare$

Let $R_0$ be an $I$-bigraded smooth algebra and $N$ be a finitely generated projective $R_0$-bimodule. Then the free algebra $R$ generated
 by $N$ over $R_0$,
 $$R=R_0\oplus N\oplus (N\otimes_{R_0}N) \oplus \dots$$
is again an $I$-bigraded smooth algebra. We assume that {\it both} Euler  forms for $R$ and $R_0$ are represented as pullbacks of forms $\chi_R$ and $\chi_{R_0}$ on $\Z^I$.

Let $W=W^{(0)}+W^{(1)}\in R/[R,R]$ be a potential which comes from an element $\widetilde{W}=\widetilde{W}^{(0)}+\widetilde{W}^{(1)}\in R_0\oplus N\subset R$.
 Then for any dimension vector $\gamma\in \Z^I_{\geqslant 0}$ the representation space $\M_\gamma$ for $R$ is a vector bundle over the space $\M_\gamma^0$
parametrizing $\gamma$-dimensional representations of $R_0$,
 and the potential $W_\gamma$ is affine along fibers. Hence we can apply the above Proposition.

We obtain
$$\mathcal{H}_\gamma\simeq \mathbf{H}^\bullet_c(\M^1_\gamma/\G_\gamma,-W_{\gamma}^{(0)})^\vee\otimes \mathbb{T}^{-\otimes \chi_{R_0}(\gamma,\gamma)}\,,$$
where $\M^1_\gamma\subset \M_\gamma^0$ is a closed $\G_\gamma$-invariant subvariety parametrizing fibers of fibration $\M_\gamma\to \M_\gamma^0$ along 
which $\widetilde{W}_\gamma^{(1)}$ vanishes. It is easy to see that $\M^1_\gamma$ is the variety of representations of a certain quotient algebra $R_1=R_0/J$ of $R_0$. The two-sided ideal $J$ is the image of the evaluation morphism
$$\op{Hom}_{R_0-\op{mod}-R_0}(N, R_0)\to R_0,\,\,\,\,\phi\mapsto \phi(\widetilde{W}^{(1)})\,.$$
For example, if $R_0$ is path algebra of a quiver $Q_0$ with vertex set $I$, quiver $Q$ is obtained by adding finitely many edges, and the potential $W^{(1)}$
is a cyclic polynomial in arrows of $Q$, which is linear in new edges, then $R_1$ is the quotient of $R_0$ by the relations $\partial W^{(1)}/\partial \alpha=0$ 
for all $\alpha \in \op{Edges}(Q)-\op{Edges}(Q_0)$. In particular, {\it any} finitely presented $I$-bigraded algebra can be obtained in such a way as $R_1=R_0/J$.

Algebra $R_1$ is isomorphic to $H^0(R_1')$, where $R_1'$ is a differential graded algebra concentrated in non-positive degrees. Forgetting the differential we see that the algebra $R_1'$
 is smooth. It is generated (over the algebra $R_0$ placed in degree $0$) by the bimodule $M:=\op{Hom}_{R_0-\op{mod}-R_0}(N, R_0\otimes R_0)$ placed in degree $-1$.
  Element $\widetilde{W}$ gives a differential $d:\alpha\mapsto mult\circ\alpha(\widetilde{W}), \alpha\in M$, where $mult:R_0\otimes R_0 \to R_0$ is the product. Clearly $d(\alpha)=0$ for $\alpha\in R_0$. If $R_1$ is quasi-isomorphic to $R_1'$ then $R_1$ has global cohomological dimension $\leqslant 2$.
 Examples include  preprojective algebras and algebras of functions on smooth affine surfaces.
In general, the abelian category of $R_1$-modules admits a natural dg-enhancement such that for any two objects one has $\op{Ext}$-groups only in degrees $0,1,2$.

\begin{rmk}
This construction indicates that it should be possible to define Cohomological Hall algebra for a larger class of abelian categories than $R-mod$, e.g.
 for the category of coherent sheaves on any smooth projective surface. For such categories moduli stacks of objects are not smooth (and not of finite type), 
but have virtual fundamental classes (as functionals on the cohomology with compact support) as well as  virtual dimensions.
\end{rmk}

         \section{Stability conditions, wall-crossing and motivic DT-series}

In this section we introduce the generating series which is an analog of the motivic
Donaldson-Thomas series from \cite{KS}. It is defined in terms of COHA of a smooth $I$-bigraded algebra
with potential. A choice of stability function gives rise to a ``PBW decomposition" of the generating series with factors parametrized by different rays in the upper-half plane. We also study how the series changes under a mutation of quiver with potential.

\subsection{Stability and Harder-Narasimhan filtration}
Let us fix a smooth $I$-bigraded algebra $R$ over $\kk$ (possibly endowed with a potential, or even with a closed 1-form), for a given finite set $I$.
Also we fix an algebraic closure $\overline{\kk}$ of the ground field $\kk$.
We have two abelian categories: $\mathfrak{A}$ which is $\kk$-linear and whose objects are $\kk$-points of schemes $\M_\gamma$,
and its $\overline{\kk}$-linear cousin $\mathfrak{A}_{\overline{\kk}}$ which is $\overline{\kk}$-linear and whose  objects are
$\overline{\kk}$-points of $\M_\gamma$. We  will work mostly with $\mathfrak{A}_{\overline{\kk}}$.

 \begin{dfn} For given $I$, a {\bf central charge} $Z$ (a.k.a.  stability function) is an additive map
$Z:\Z^I\to \C$ such that the image of any of the standard base vectors lies in the upper-half plane  $\mathbb{H}_+:=\{z\in \C\,|\,\op{Im}(z)>0\}$.
Central  charge $Z$ is called {\bf generic} if there are no two $\Q$-independent elements of ${\Z}_{\geqslant 0}^I$ which are mapped by $Z$ to the same straight line.
\end{dfn}
Let us fix a central charge $Z$. Then  for any non-zero object $0\ne E \in \op{Ob} (\mathfrak{A}_{\overline{\kk}})$ one defines
$$\op{Arg}(E):=\op{Arg}(Z(cl(E)))\in (0,\pi)\,\,,$$
where $\gamma=cl(E)\in \Z_{\geqslant 0}^I$ is the dimension vector of the object $E$, i.e. $E\in \M_\gamma(\overline{\kk})$.
 We will also use the shorthand notation $Z(E):=Z(cl(E))$.
\begin{dfn} A non-zero object $E$  is called {\bf semistable} (for the central charge $Z$) if there is no non-zero subobject $F\subset E$ such that
 $\op{Arg}(F)>\op{Arg}(E)$.
\end{dfn}
It is easy to see that  the set of semistable objects is the set  $\overline{\kk}$-points of
a Zariski open $\G_\gamma$-invariant subset $\M_\gamma^{ss}\subset \M_\gamma(\overline{\kk})$
 defined over $\kk$. In particular it is smooth.

Any non-zero object $E$ of $\mathfrak{A}_{\overline{\kk}}$ admits a canonical {\bf Harder-Narasimhan filtration} (HN-filtration in short), i.e. an increasing filtration $0=E_0\subset E_1\subset\dots \subset E_n=E$ with $n\geqslant 1$ such that all the quotients $F_i:=E_i/E_{i-1},\,i=1,\dots,n$ are semistable and
$$\op{Arg}(F_1)>\dots>\op{Arg}(F_n)\,\,.$$

It will be convenient to introduce a total order $\prec $ (a lexicographic order in polar coordinates) on the upper-half plane by
$$z_1\prec z_2 \mbox { iff } \op{Arg}(z_1)>\op{Arg}(z_2)\mbox{ or }\left\{\op{Arg}(z_1)=\op{Arg}(z_2)\mbox{ and } |z_1|>|z_2|
\right\}\,. $$
One can show easily that if $0=E_0\subset E_1\subset \dots$ is a HN-filtration of a non-zero object $E$, then $E$ does not contain a  subobject $F\subset E$ with
$$Z(F)\prec Z(E_1)\,\,,$$
and the unique non-zero subobject $F'$ of $E$ with $Z(F')=Z(E_1)$ is $F'=E_1$.

For a non-zero $\gamma\in \Z_{\geqslant 0}^I$ let us denote by $\mathcal{P}(\gamma)$ the set of collections
$$\gamma_\bullet=(\gamma_1,\dots,\gamma_n),\,\,n\geqslant 1\,,\mbox{ such that }\sum_{i=1}^n\gamma_i=\gamma\,\,,$$
where $\gamma_i,\,i=1,\dots,n$ are non-zero elements of $\Z_{\geqslant 0}^I$ satisfying
$$\op{Arg}(Z(\gamma_1))>\op{Arg}(Z(\gamma_2))>\dots >\op{Arg}(Z(\gamma_n))\,\,.$$
  We introduce a partial ``left lexicographic'' order on by $\mathcal{P}(\gamma)$ by the formula
$$(\gamma'_1,\dots,\gamma'_{n'})<(\gamma_1,\dots,\gamma_{n})\mbox{ if }\exists i\,,\,\,1\leqslant i \leqslant \op{min}(n,n')\mbox{ such that }$$
$$ \gamma_1'=\gamma_1,\dots, \gamma_{i-1}'=\gamma_{i-1}\mbox{ and }Z(\gamma_i')\prec Z(\gamma_i)\,\,.$$
This order is in fact a total order for generic $Z$.
We denote by $\M_{\gamma;\gamma_1,\dots,\gamma_n}$ the constructible subset of $\M_\gamma$ whose $\overline{\kk}$-points are objects
$E\in \M_\gamma(\overline{\kk})$ which admit
an increasing filtration $0=E_0\subset E_1\subset\dots \subset E_n$ such that
$$cl(E_i/E_{i-1})=\gamma_i\,\,,\,\,\,i=1,\dots,n\,\,.$$

One can  see that the subset $\M_{\gamma;\gamma_1,\dots,\gamma_n}$ is closed, because 
$\M_{\gamma;\gamma_1,\dots,\gamma_n}$ is the image under projection to $\M_\gamma$
of a closed subset in the product of $\M_\gamma$ and of the appropriate flag varieties associated with the vertices of the quiver.

\begin{lmm} For any non-zero  $\gamma\in \Z_{\geqslant 0}^I$ and any $\gamma_\bullet\in \mathcal{P}(\gamma)$ the complement
$$\M_{\gamma, \gamma_\bullet}-  \cup_{\gamma'_\bullet< \gamma_\bullet} \M_{\gamma,\gamma'_\bullet}   $$
coincides with the constructible subset $\M_{\gamma,\gamma_\bullet}^{HN}\subset\M_\gamma$ consisting of representations whose HN-filtration has type $\gamma_\bullet$.
\end{lmm}
{\it{Proof}.} First, take $E\in\M_{\gamma,\gamma_\bullet}^{HN}$ and assume that $E\in \M_{\gamma,\gamma_\bullet'}$
 for some $\gamma'_\bullet<\gamma_\bullet$. Hence we have two filtrations on $E$, the HN-filtration $E_\bullet^{HN}$ and
 a filtration $E_\bullet'$ with quotients with dimension vectors $\gamma_\bullet'$.
Denote by $i$ the minimal index for which $\gamma_i\ne \gamma_i'$.
 It is easy to show by induction that both filtrations coincide up to index $i-1$, and that the inequality $Z(\gamma_i')\prec Z(\gamma_i)$
is not compatible with the stability of the quotient $E_i^{HN}/E_{i-1}^{HN}$. Hence our assumption can not be satisfied.

Conversely, let us assume that $E\in \M_{\gamma,\gamma_\bullet}$ does not belong to $\M_{\gamma,\gamma_\bullet}^{HN}$.
  Denote by $E_\bullet^{HN}$ the HN-filtration for $E$, and by $E_\bullet$ any filtration with dimensions of quotients given by the vector
 $\gamma_\bullet$.
  Let us denote by $i$ the minimal index for which $\gamma_i^{HN}\ne \gamma_i$. We have as above
 $E_j=E_j^{HN}$ for all $j<i$, and $E_i\ne E_i^{HN}$. Then  $E_i^{HN}/E_{i-1}^{HN}$ is the lowest non-zero term
in the HN-filtration for $E/E_{i-1}^{HN}$, and $E_i/E_{i-1}^{HN}\ne E_i^{HN}/E_{i-1}^{HN}$.
Hence we have $Z(\gamma_i^{HN})\prec Z(\gamma_i)$.
Therefore we have $\gamma_\bullet ^{HN}<\gamma_{\bullet}$, and $E$ belongs to the union $\cup_{\gamma'_\bullet< \gamma_\bullet} \M_{\gamma,\gamma'_\bullet}   $.
$\blacksquare$

 For any non-zero $\gamma\in \Z_{\geqslant 0}^I$ and any $\gamma_\bullet \in \mathcal{P}(\gamma)$ we define a closed set
$$\M_{\gamma,\leqslant \gamma_\bullet}:=\cup_{\gamma'_\bullet\leqslant \gamma_\bullet} \M_{\gamma,\gamma'_\bullet}\,\,.$$
Let us assume for simplicity that $Z$ is generic, hence the partial order $<$
 is a total order.
Above considerations  imply that we get a chain of closed sets
 whose consecutive differences are exactly
 sets $\M_{\gamma,\gamma_\bullet}^{HN}$.
Taking the complements, we obtain a chain of open subspaces
$$\emptyset=U_0\subset U_1\subset\dots\subset U_{N_\gamma}=\M_\gamma\,,$$
where $N_\gamma:=\#\mathcal{P}(\gamma)$ is the number of elements of the finite set $\mathcal{P}(\gamma)$. The consecutive differences $U_i- U_{i-1}$ are exactly locally closed spaces $\M_{\gamma,\gamma_\bullet}^{HN}$.
This chain is obviously invariant under the action of $\G_\gamma$, and we obtain
 a chain of open subspaces
  $$U_l^\univ\subset \M_\gamma^\univ\,\,,\,\,l=0,\dots,N_\gamma\,\,.$$

\subsection{Spectral sequence converging to $\mathcal{H}_\gamma$}

The chain of open subspaces $(U_l)_{l=0,\dots,N_\gamma}$ introduced above, leads to
 a spectral sequence converging to $\mathcal{H}_\gamma={{\bf H}}^\bullet_{\G_\gamma}(\M_\gamma,W_\gamma)$ with the first term
$$\bigoplus_{l=1}^{N_\gamma} {\bf H}^\bullet_{\G_\gamma}(U_l,U_{l-1},W_\gamma)\,\,.$$

\begin{thm} In the above notation, the first term  is
 isomorphic to
 $$\bigoplus_{n\geqslant 0}\bigoplus_{\substack{\gamma_1,\dots,\gamma_n\in \Z_{\geqslant 0}^I- 0\\
\op{Arg}\gamma_1>\dots>\op{Arg}\gamma_n}  }  {\bf H}^\bullet_{\prod_{i\in I}\G_{\gamma_i}}(\M^{ss}_{\gamma_1}\times\dots\times \M^{ss}_{\gamma_n},W_\gamma)\otimes
\mathbb{T}^{\otimes \sum_{i<j} (-\chi_R(\gamma_i,\gamma_j))}\,\,.   $$
\end{thm}
{\it Proof.}
First, for a given decomposition $\gamma_\bullet=(\gamma_1,\dots,\gamma_n)\in \mathcal{P}(\gamma)$,
 we define  $\M_{\gamma_1,\dots,\gamma_n}\subset \M_{\gamma}$ as the space of  representations
 of $R$ in coordinate spaces, such that collections of standard coordinate subspaces $(\overline{\kk}^{\sum_{j\leqslant k}\gamma_i^j})_{i\in I}$
  form a subrepresentation of $R$ for any $k=0,\dots,n$. This space is smooth, because it is the scheme of representations of $R$
in a finite-dimensional algebra (the product over $i\in I$ of algebras of block upper-triangular matrices).
 Denote by $\M_{\gamma_1,\dots,\gamma_n}^{ss}\subset \M_{\gamma_1,\dots,\gamma_n}$ the open subspace consisting of representations
 such that the associated subquotients of dimension vectors $\gamma_i,\,\,i=1,\dots, n$ are semistable.
In particular, $\M_{\gamma_1,\dots,\gamma_n}^{ss}$ is smooth.

It follows from definitions that
$$\M^{HN}_{\gamma,\gamma_\bullet}\simeq\G_\gamma\times_{\G_{\gamma_1,\dots,\gamma_n}} \M_{\gamma_1,\dots,\gamma_n}^{ss}\,\,,$$
where $\G_{\gamma_1,\dots,\gamma_n}\subset \G_\gamma$ is the group of invertible block upper-triangular matrices, acting naturally on
$\M_{\gamma_1,\dots,\gamma_n}$ and on $\M_{\gamma_1,\dots,\gamma_n}^{ss}$.
 This formula shows that $\M^{HN}_{\gamma,\gamma_\bullet}$ is a smooth manifold. Its codimension in $\M_\gamma$ can be easily calculated:
$$\op{codim} \M^{HN}_{\gamma,\gamma_\bullet}= -\sum_{i<j} \chi_R(\gamma_i,\gamma_j)\,\,.$$
 For given $\gamma_\bullet\in \mathcal{P}(\gamma)$, the component of the first term of the spectral sequence, corresponding to $\M^{HN}_{\gamma,\gamma_\bullet}$, has the form
$${\bf H}^\bullet_{\G_\gamma}(U_l,U_{l-1},W_\gamma)$$
for some $l,\,1\leqslant l\leqslant N_\gamma$, where $\M^{HN}_{\gamma,\gamma_\bullet}$ is a closed $\G_\gamma$-invariant submanifold of the  manifold $U_l$, with the complement
  $U_l- \M^{HN}_{\gamma,\gamma_\bullet}$ equal to $U_{l-1}$.
Hence, by Thom isomorphism 
we see that this component
 can be written as
 $${\bf H}^\bullet_{\G_\gamma}(\M^{HN}_{\gamma,\gamma_\bullet},W_\gamma)\otimes \mathbb{T}^{\otimes(- \sum_{i<j} \chi_R(\gamma_i,\gamma_j))}\,\,.$$
Finally, using the isomorphism
$${\bf H}^\bullet_{\G_\gamma}(\M^{HN}_{\gamma,\gamma_\bullet},W_\gamma)\simeq {\bf H}^\bullet_{\G_{\gamma_1,\dots,\gamma_n}}(\M_{\gamma_1,\dots,\gamma_n}^{ss},W_\gamma)$$
and
$\mathbb{A}^1$-homotopies
$$\M^{ss}_{\gamma_1}\times\dots\times \M^{ss}_{\gamma_n}\sim \M_{\gamma_1,\dots,\gamma_n}^{ss},\,\,\,
\G_{\gamma_1}\times\dots\times\G_{\gamma_n}\sim \G_{\gamma_1,\dots,\gamma_n}\,\,,$$
we obtain the statement of the theorem. $\blacksquare$

Alternatively, one can use cohomology with compact support of smooth Artin stacks $\M^{HN}_{\gamma,\gamma_\bullet}/\G_\gamma$ with potential $W_\gamma$, see Section 5.4.

Let $V\subset \mathbb{H}_+$ be a  sector, i.e. $V$=$V+V$=$\R_{> 0}\cdot V$ and $0\notin V$.
 Denote by $\M_{V,\gamma}\subset \M_\gamma$ the set of all representations
 whose $HN$-factors have classes in $Z^{-1}(V)$. It is easy to see that for any $\gamma$
   the set $\M_{V,\gamma}$ is Zariski open and $\G_\gamma$-invariant. We define the {\bf Cohomological Hall vector space} by
 $$\mathcal{H}_V:=\oplus_\gamma \mathcal{H}_{V,\gamma}=\oplus_\gamma {\bf H}^\bullet_{\G_\gamma}(\M_{V,\gamma},W_\gamma)\,.$$

Similarly to Theorem 5, one has a spectral sequence converging to $\mathcal{H}_V$ where we use only rays lying in $V$.
In general the space $\mathcal{H}_V$ does not carry a product, because the corresponding Grassmannians are not necessarily {\it closed}.
Nevertheless, for $V$ being a ray $l=\op{exp}(i\phi)\cdot \R_{> 0}$ (as well as for $V$ equal to the whole upper-half plane),
 the product is well-defined (and is twisted associative) because the corresponding Grassmannians are closed.

If we assume that the restriction of the bilinear form $\chi_R$ to the sublattice $Z^{-1}\left(\op{exp}(i\phi)\cdot \R\right)$ is {\it symmetric}
(e.g. it always holds for generic $Z$),  then we have the modified untwisted associative product on
 $$\H_l^{mod}:=\oplus_{\gamma\in Z^{-1}(l)}\,\, \mathcal{H}_{l,\gamma}\otimes
\left(\mathbb{T}^{\otimes 1/2} \right)^{\otimes \chi_R(\gamma,\gamma)},\,\,\,l:=\op{exp}(i\phi)\cdot \R_{>0}$$

The decomposition in
Theorem 5 means that  the twisted algebra $\mathcal{H}=\mathcal{H}_{\mathbb{H}_+}$ ``looks like'' (i.e. related by a spectral sequence with) a
clockwise  ordered product of twisted algebras $\mathcal{H}_l$ over the set of rays $l$.
This is similar to a decomposition for universal enveloping algebras which follows from the Poincar\'e-Birkhoff-Witt theorem.
 Namely, if $\mathfrak{g}=\oplus_{\gamma \in \Z_{\geqslant 0}^I -\{0\}} \mathfrak{g}_\gamma$ is a graded Lie algebra (over any field), then
for any ray $l\subset \mathbb{H}_+$ we have Lie subalgebra $\mathfrak{g}_l:=\oplus_{\gamma \in Z^{-1}(l)} \mathfrak{g}_{\gamma}$, and
$$ \bigotimes_l^{\curvearrowright} U\mathfrak{g}_l\stackrel{\sim}{\longrightarrow}U\mathfrak{g}\,.$$
In the case of COHA the twisted algebras $\mathcal{H}_l$ are {\it not} naturally realized as twisted subalgebras of $\mathcal H$, there is only a
 homomorphism of graded spaces $\mathcal{H}\to \mathcal{H}_l $
 given by the restriction morphism for the inclusions $\M^{ss}_\gamma\mono\M_\gamma$.

In the special case of quiver $A_2$ and Betti cohomology, there are embeddings of  algebras $\mathcal{H}_l\mono \mathcal{H}$ (with cohomological $\Z/2\Z$-grading), as follows from direct calculations
 (see Section 2.8), and the clockwise ordered product indeed induces an isomorphism
$$ \bigotimes^{\curvearrowright}_l \mathcal{H}_l \stackrel{\sim}{\longrightarrow} \mathcal{H}\,.$$

Finally, we should warn the reader of a possible caveat. Namely, the calculation from section 4.8 is {\it not} applicable to $\mathcal{H}_{V,\gamma}$
 for general $V$ because corresponding spaces $\M_{V,\gamma}$ are only Zariski open subsets, and do not form themselves  fibrations by affine spaces. Therefore the stability in the ambient category of cohomological dimension $1$ is not directly related with the stability in the smaller category of cohomological dimension $2$.

\subsection{Generating series and quantum tori}

The cohomology theory ${\bf H}^\bullet$
used in the definition of Cohomological Hall algebra takes values in a Tannakian category $\mathcal C$.
 In other words, cohomology carries an action of a pro-affine algebraic group, which we call the {\bf  motivic Galois group} $\op{Gal}^{mot}_{{\bf H}}$
 for the theory $\bf H$.
We assume that there is a notion of weight filtration in $\mathcal C$. For example, in the case of rapid decay Betti (tensored by $\Q$) or de Rham cohomology we consider
 ${\bf H}^\bullet(X)$ for any variety $X$ as a vector space graded by cohomological degree and endowed
with the weight filtration. In terms of $\op{Gal}^{mot}_{{\bf H}}$ this means that we have an embedding $w:\mathbb{G}_m\mono \op{Gal}^{mot}_{{\bf H}}$
(defined up to conjugation) such that the Lie algebra of $\op{Gal}^{mot}_{{\bf H}}$
 has non-positive weights with respect to the adjoint $\mathbb{G}_m$-action.  For any representation $E$
of  $\op{Gal}^{mot}_{{\bf H}}$ the weight filtration is defined by
$$W_i E:=\oplus_{j\leqslant i} E_j\,,\,\,\,i,j\in \Z\,,$$
where $E_j\subset E$ is the eigenspace of $w(\mathbb{G}_m)$ with weight $j$.

Let us consider the $K_0$-ring $\mathcal{M}$ of tensor category $\mathcal C$ (i.e. of the category of finite-dimensional representations of $\op{Gal}^{mot}_{{\bf H}}$). It contains an invertible  element $\mathbb{L}$ corresponding to the Tate motive  ${\bf H}^2(\mathbb{P}^1)$
of weight $+2$ .  We complete $\mathcal{M}$ by adding infinite sums of pure motives with weights approaching to $ +\infty$. Notice that this completion (which we denote by $\widehat{\mathcal{M} }$)
differs from the completion used in the theory of motivic
 integration, where weights are allowed to
go to $-\infty$.

With any $I$-bigraded smooth algebra $R$ and a bilinear form $\chi_R$ compatible with the Euler form,
we associate the following series with coefficients in $\widehat{\mathcal{M} }$:
$$A=A^{(R,W)}:=\sum_{\gamma\in \Z^I_{\geqslant 0}} [\mathcal{H}_\gamma]\, \mathbf{e}_\gamma\,\,, $$
where variables $\mathbf{e}_\gamma$ are additive generators of the associative unital algebra ${\cal R}_+$ over $\mathcal{M}$ isomorphic to the subalgebra of the motivic quantum torus ${\cal R}$ introduced in \cite{KS}. The relations are given by the formulas
$$\mathbf{e}_{\gamma_1} \cdot \mathbf{e}_{\gamma_2}=\mathbb{L}^{-\chi_R(\gamma_1,\gamma_2)}\mathbf{e}_{\gamma_1+\gamma_2}\,\,\,\,\,\forall \gamma_1,\gamma_2\in \Z^I_{\geqslant 0}\,,\,\,\,
\mathbf{e}_{0}={1}\,,$$
and the coefficient of the series $A$ is given by
$$[\mathcal{H}_\gamma]:=\sum_{k\geqslant 0} (-1)^k [{\bf H}^k_{\G_\gamma}(\M_\gamma,W_\gamma)]\in \widehat{\mathcal{M} }\,\,.$$

Series $A$ belongs to the completion $\widehat{{\cal R}}_+$ consisting of infinite series in $\mathbf{e}_\gamma$. It has the form
 $$A={1}+\mbox{ higher order terms}\,\,\,$$
and is therefore  invertible.
\begin{dfn} We call $A$ the {\bf motivic Donaldson-Thomas series} of the pair $(R,W)$.

\end{dfn}

The series $A$ can be written in a more explicit form using the following result.
\begin{lmm} For any scheme $X$ of finite type endowed with an action of the affine algebraic group $G=\prod_i\op{GL}(n_i)$ and a function
$f\in \mathcal{O}(X)^G$, one has
the following identity in $\widehat{\mathcal{M} }$:
$$[{\bf H}^\bullet_G(X,f)]=[{\bf H}^\bullet(X,f)]\cdot [{\bf H}^\bullet(\op{B}G)]=\frac{[{\bf H}^\bullet(X,f)]}{\prod_i (1-\mathbb{L})\dots (1-\mathbb{L}^{n_i})
}\,\,.$$
\end{lmm}
{\it Proof.}  There is a filtration of (a finite approximations of) the classifying space $\op{BG}$ by open subschemes such that the universal
 $G$-bundle is trivial in consecutive differences. The spectral sequence associated with such a filtration gives the result. $\blacksquare$

In the basic example of the quiver $Q_0=A_1$ with one vertex and no edges, the motivic DT-series  is
$$A=1+\frac{1}{1-\mathbb{L}}\mathbf{e}_1+\frac{1}{(1-\mathbb{L})(1-\mathbb{L}^2)}\mathbf{e}_2+\dots=$$
$$=\sum_{n\geqslant 0}\frac{
{\mathbb{L}}^{n(n-1)/2} }
{(1-\mathbb{L})\dots (1-\mathbb{L}^n)}\mathbf{e}_1^n=
(-\mathbf{e}_1;\mathbb{L})_\infty=\prod_{n\geqslant 0} (1+\mathbb{L}^n \mathbf{e}_1)\,.$$

Let $Z:\Z^I\to \C$ be a central charge. Then for any  sector $V\subset \mathbb{H}_+$ we define
the {\bf motivic DT-series associated with $V$} such as follows:

$$A_V:=\sum_{\gamma\in \Z^I_{\geqslant 0}} [\mathcal{H}_{V,\gamma}]\, \mathbf{e}_\gamma\,\,. $$
The spectral sequence from Section 5.2 implies the following {\bf Factorization Formula}:
$$A=A_{\mathbb{H}_+}=\prod^{\curvearrowright}_l A_{l}\,,$$
where the product is taken in the clockwise order over all rays $l\subset \mathbb{H}_+$ containing non-zero points in $Z(\Z^I_{\geqslant 0})$ (i.e. in the decreasing order of the arguments of complex numbers). Each factor $A_l$ corresponds to semistable objects with the central charge in $l$:
$$A_l=1+\sum_{\gamma\in Z^{-1}(l)} \sum_{k\geqslant 0} (-1)^k[{\bf H}^k_{\G_{\gamma}}(\M^{ss}_{\gamma},W_{\gamma})]\,\mathbf{e}_{\gamma}\,\in\widehat{{\cal R}}_+\,.$$
Also, for any pair of disjoint  sectors $V_1,V_2 \subset \mathbb{H}_+$ whose union is also a sector, and such that $V_1$ lies on the left of
$V_2$, we have
$$A_{V_1\cup V_2}=A_{V_1} A_{V_2}\,.$$

It will be convenient also to rewrite series $A$ in rescaled variables.
 Let us add an element $\mathbb{L}^{1/2}$ to $\widehat{\mathcal{M} }$, satisfying the relation
 $\left(\mathbb{L}^{1/2}\right)^2=\mathbb{L}$
and define rescaled quantum variables by
$$\widehat{\mathbf{e}}_\gamma:=
\left(-\mathbb{L}^{1/2}\right)^{-\chi_R(\gamma,\gamma)}
\,\mathbf{e}_\gamma\,\,.$$
The rescaled variables satisfy the relation
$$\widehat{ \mathbf{e}}_{\gamma_1}\cdot \widehat{ \mathbf{e}}_{\gamma_2}=\left(-\mathbb{L}^{1/2}
\right)^{-\langle\gamma_1,\gamma_2\rangle}
\widehat{ \mathbf{e}}_{\gamma_1+\gamma_2}\,,$$
where
$$\langle\gamma_1,\gamma_2\rangle:=\chi_R(\gamma_1,\gamma_2)-\chi_R(\gamma_2,\gamma_1)$$
is a {\it skew-symmetric} bilinear form on $\Z^I$.
Then we have:
$$A= \sum_{\gamma\in \Z^I_{\geqslant 0}} \left(-\mathbb{L}^{1/2}\right)^{\chi_R(\gamma,\gamma)}[\mathcal{H}_\gamma]\,
\widehat{\mathbf{e}}_\gamma= \sum_{\gamma\in \Z^I_{\geqslant 0}} \frac{\left(-\mathbb{L}^{1/2}\right)^{\chi_R(\gamma,\gamma)}[{\bf H}^\bullet(\M_\gamma,
W_\gamma)]}
{\prod_i (1- \mathbb{L})\dots(1- \mathbb{L}^{\gamma^i}) }\,
\widehat{\mathbf{e}}_\gamma\,\,= $$
$$=\sum_\gamma \left(-\mathbb{L}^{1/2}\right)^{-\dim \M_\gamma/\G_\gamma}[{\bf H}^\bullet_{\G_\gamma}(\M_\gamma,W_\gamma)]
\widehat{\mathbf{e}}_\gamma\,\,. $$
In the basic example of quiver $Q_0$ the motivic DT-series is $A=(\mathbb{L}^{{1/2}}\widehat{\mathbf{e}}_1;\mathbb{L})_\infty$.

In the case when $\chi_R$ is symmetric (i.e. $\langle\,,\,\rangle=0$), the rescaled variables are monomials in commuting variables
 $$\widehat{ \mathbf{e}}_{\gamma_1}\cdot \widehat{ \mathbf{e}}_{\gamma_2}=
\widehat{ \mathbf{e}}_{\gamma_1+\gamma_2}\,,$$
and the motivic DT-series $A$ is the generating  series for the Serre polynomials of graded components of the associative (non-twisted) algebra $\H^{mod}$ with the modified product (see Section 2.7):
$$A= \sum_{\gamma\in \Z^I_{\geqslant 0}} [\mathcal{H}_\gamma^{mod}]\,
\widehat{\mathbf{e}}_\gamma\,.$$
In general, coefficients of the series $A$ are rational functions in the variable $\mathbb{L}^{1/2}$    with coefficients in $\mathcal{M} $.

There is a natural homomorphism given by the Serre polynomial
$$S: \mathcal{M}[\mathbb{L}^{\pm 1/2} ]\to \Z[q^{\pm 1/2}]\,,\,S([E]):=\sum_{n\in \Z}  \op{dim}\op{Gr}_n^W(E)\,q^{n/2},\,\
S(\mathbb{L}^{1/2}):=q^{1/2}\,,$$
where $\op{Gr}_w^W$ denotes the graded component of weight $w$ with respect to the  weight filtration.
Applying $S$ to $A$ we obtain a series $S(A)=\sum_{\gamma}S([\mathcal{H}_{\gamma}]){\mathbf{e}}_\gamma$ in quantum variables with coefficients in
$$\Z[q^{\pm 1/2},\left((1-q^n)^{-1}\right)_{n=1,2,\dots}]\subset \Q((q^{1/2}))\,.$$
For {\it pure} COHA the series $S([\mathcal{H}_{\gamma}])$ for any $\gamma\in \Z^I_{\geqslant 0}$ coincides with the Hilbert-Poincar\'e series of the graded space $\H_\gamma$ in variable $(-q^{1/2})$. Therefore for the quiver $Q_d$ with zero potential we obtain the same generating series as in Section 2.5.

\subsection{Reformulation using cohomology with compact support}

It is instructive to revisit the factorization formula for the generating functions using cohomology with compact support.
If we apply the involution $D:\mathcal{M} \to \mathcal{M} $ given by passing to dual motive, and extend it to the
square root of Tate motive by
 $$D(\mathbb{L}^{1/2})=\mathbb{L}^{-1/2}$$
we obtain the series
$$D(A)=\sum_{\gamma\in \Z^I_{\geqslant 0}}[\mathbf{H}_{c,\G_\gamma}(\M_\gamma,-W_\gamma) ]D(\mathbf{e}_\gamma)=
\sum_{\gamma\in \Z^I_{\geqslant 0}}\frac{[\mathbf{H}_{c}(\M_\gamma,-W_\gamma)]}{[\mathbf{H}_c(\G_\gamma)]}
D(\mathbf{e}_\gamma)=$$
$$=\sum_{\gamma\in \Z^I_{\geqslant 0}}\frac{\left(-\mathbb{L}^{1/2}\right)^{\chi_R(\gamma,\gamma)}[{\bf H}^\bullet_c(\M_\gamma,-W_\gamma)]}
{\prod_i (\mathbb{L}^{\gamma^i}- 1)\dots(\mathbb{L}^{\gamma^i}- \mathbb{L}^{\gamma^i-1}) }\,
D(\widehat{\mathbf{e}}_\gamma)\,\,, $$
where the dual quantum variables satisfy the relations:
$$  D(\mathbf{e}_{\gamma_1})\cdot D( \mathbf{e}_{\gamma_2})=\mathbb{L}^{-\chi_R(\gamma_2,\gamma_1)}
D( \mathbf{e}_{\gamma_1+\gamma_2})\,,$$
$$  D(\widehat{ \mathbf{e}}_{\gamma}) =  \left(-\mathbb{L}^{1/2}\right)^{-\chi_R(\gamma,\gamma)}
D(\mathbf{e}_\gamma)\,,\,\,\, D(\widehat{ \mathbf{e}}_{\gamma_1})\cdot D(\widehat{ \mathbf{e}}_{\gamma_2})=\left(-\mathbb{L}^{1/2}
\right)^{+\langle\gamma_1,\gamma_2\rangle}
D(\widehat{ \mathbf{e}}_{\gamma_1+\gamma_2})\,.$$

The series $D(A)$ is analogous to the generating series considered in \cite{KS} (up to the sign change $\mathbb{L}^{1/2}\mapsto - \mathbb{L}^{1/2}  $).
Namely, we can associate with $R$ an ind-constructible 3-dimensional Calabi-Yau category with the $t$-structure whose heart is equivalent to the full subcategory of the category of finite-dimensional $R$-modules, consisting
of critical points of $W_\gamma,\,\,\,\gamma\in \Z^I_{\geqslant 0}$.
This category is ind-constructible, it carries a canonical orientation (in the sense of \cite{KS}), and its motivic DT-series introduced in the loc.cit. has essentially the same form as $D(A)$, but with the rapid decay cohomology replaced by the cohomology with coefficients in the sheaf of vanishing cycles. All this will be discussed in Section 7.10.

The argument with spectral sequences in Section 5.2 becomes more direct for the compactly supported cohomology.
Namely, for any stratified Artin stack $\mathcal S$ endowed with a function $f$ we have the identity
$$[\mathbf{H}_c(\mathcal{S},f)]=\sum_\alpha[\mathbf{H}_c(\mathcal{S}_\alpha,f)]\,,$$
where $(\mathcal{S}_\alpha)$ are strata.

The calculation in Section 4.8 is also more convenient in the case of cohomology with compact support. We will illustrate it below in Section 5.6. In particular, in the case when the potential is linear in additional arrows we obtain series
$$\sum_\gamma [\mathbf{H}_c(S_\gamma)]\cdot (-\mathbb{L}^{1/2})^{-\op{dim}_{virt} S_\gamma} D(\widehat{\mathbf{e}}_\gamma)\,,\,\,\,$$
where $S_\gamma=\M^1_\gamma/\G_\gamma$ is the stack of representations of the dg-algebra $R_1'$ (in the notation from Section 4.8) of dimension vector $\gamma$. The virtual dimension of $S_\gamma$ is equal to $-\chi_{R_1'}(\gamma,\gamma)$, where $\chi_{R_1'}$ is a lift to $\Z^I$ of the Euler form of the dg-category of finite-dimensional $R_1'$-modules.

In the case when the ground field $\kk=\mathbb{F}_q$ has finite  characteristic $p>0$ we can apply the multiplicative numerical invariant from Section 4.6.
The resulting series is
$$\sum_{\gamma\in \Z^I_{\geqslant 0}}\sum_{\mathcal{E}\in \M_\gamma(\kk)/\op{iso}}\frac{\op{exp}\left(-\frac{2\pi i}{p}\op{Trace}_{\mathbb{F}_q/\mathbb{F}_p}(W(\mathcal{E}))\right)}{|\op{Aut}(\mathcal{E})|} D(\mathbf{e}_\gamma)\,.$$
Also, the multiplication of the variables $D(\mathbf{e}_\gamma)$ has a very transparent form:
$$\frac{1}{|\op{Aut}(\mathcal{E}_1)|}D(\mathbf{e}_{\gamma_1})\cdot \frac{1}{|\op{Aut}(\mathcal{E}_2)|}D(\mathbf{e}_{\gamma_2})=
\sum_{\mathcal{E}_1\to \mathcal{E}\to\mathcal{E}_2}\frac{1}{|\op{Aut}(\mathcal{E}_1\to \mathcal{E}\to\mathcal{E}_2)|}D(\mathbf{e}_{\gamma_1+\gamma_2})\,,$$
where $\mathcal{E}_1,\mathcal{E}_2$ have dimension vectors $\gamma_1,\gamma_2$ respectively.

\subsection{Classical limit}

Let us add to the quantum torus additional invertible generators $(\mathbf{f}_i)_{i\in I}$ satisfying the relations
$$ \mathbf{f}_i\cdot  \mathbf{f}_j =\mathbf{f}_j \cdot \mathbf{f}_i\,,\,\,\,\,\mathbf{f}_i\cdot\mathbf{e}_\gamma =\mathbb{L}^{\gamma^i}\mathbf{e}_\gamma\cdot \mathbf{f}_i
\,\,\,\,\forall i,j\in I,\,\gamma\in \Z^I\,.$$
For any $i_0\in I$ we define an invertible series $A^{(i_0)}$  in variables $\mathbf{e}_\gamma$ by the formula
$$A\cdot \mathbf{f}_{i_0}\cdot A^{-1}=A^{(i_0)}\cdot \mathbf{f}_{i_0} \Leftrightarrow A^{(i_0)}:=
A\cdot\left(\mathbf{f}_{i_0} A \,\mathbf{f}_{i_0}^{-1}\right)^{-1}.$$
More generally, we define an element $A^{(\gamma)}$ for any $\gamma\in \Z^I$ such as follows:
$$A^{(\gamma)}:=A\cdot\left(\mathbf{f}^\gamma A \,\mathbf{f}^{-\gamma}\right)^{-1}\,,\,\,\,\,\,\mathbf{f}^\gamma:=\prod_{i\in I}
\mathbf{f}_i^{\gamma^i}.$$
It is equal to a finite product of elements $A^{(i)}$ and of  conjugates of these elements by monomials in $\mathbf{f}_j$
(or inverses of such expressions). Also one has
$$A\cdot \mathbf{e}_\gamma\cdot A^{-1}=A^{(\gamma')}\cdot \mathbf{e}_\gamma\,,\,\,\,(\gamma')^i:=\sum_j c_{ij} \gamma^j.$$
where $(c_{ij})_{i,j\in I}$ is the matrix in the standard basis of the skew-symmetric form $\langle\,,\,\rangle=\chi_R-(\chi_R)^t$ on $\Z^I$. E.g. for the case of a quiver in notation of Section 2.1 one has 
$$c_{ij}=a_{ji}-a_{ij}\,.$$
\begin{thm} For every $\gamma\in \Z^I$  each coefficient of the series $A^{(\gamma)}$  belongs to the subring $\mathcal{M}$ of $\widehat{\mathcal{M}}$.
\end{thm}
Therefore, we can apply Euler characteristic and obtain a symplectomorphism $T$ of the classical symplectic double torus\footnote{Symplectic double torus is naturally associated with the ``double'' 
lattice $\Z^I\oplus ({\Z^I})^{\vee}$, see \cite{KS}.} with coordinates
  $$(\mathbf{e}_i^{cl},\mathbf{f}_i^{cl}\,)_{i\in I}\,,\,\,\,$$
given by the following automorphism of the algebra of functions:
$$\mathbf{e}_i^{cl}\mapsto \prod_{j\in I} \left(A^{(j),cl}
\right)^{c_{ji}} \mathbf{e}_i^{cl}\,,\,\,\,\,\,\,\,\,\mathbf{f}_i^{cl} \mapsto A^{(i),cl}\, \mathbf{f}_i^{cl}\,,$$
where $A^{(i),cl}$ is a series with integer coefficients in variables $(\mathbf{e}_j^{cl})_{j\in I}$ given by
$$A^{(i),cl}:=\chi (A^{(i)})\,.      $$
Here $\chi$ denotes the Euler characteristic. Symplectomorphism $T$ is the classical limit $\op{lim}_{q^{1/2}\to 1} \op{Ad}(S(A))$ where $S$ denotes Serre polynomial (cf. \cite{KS}).

{\it Proof}. Let us fix $i_0\in I$. It is sufficient to prove that $A^{(i_0)}$ has coefficients in $\mathcal{M}$. Let us consider a larger index set
 $\widetilde{I}=I\sqcup \{i_0^*\}$, where $i_0^*$ is a new vertex, and a larger algebra $\widetilde{R}$ bigraded by $\widetilde{I}$ in such a way that
 a representation $\widetilde{E}$ of $\widetilde{R}$ is the same as a triple $(E,V,\phi)$, where $E$ is a representation of $R$, vector space $V=E_{i_0^{\ast}}:=p_{i_0^*}(\widetilde{E})$ is
  associated with $i_0^*$, and $\phi$ is a map from $V$ to the space $p_{i_0} E$. Here $p_i$ denotes the projector corresponding to the vertex $i$. If $R$ is the path algebra of a quiver $Q$ with the set of vertices $I$, then $\widetilde{R}$ is the path algebra of a new quiver $\widetilde{Q}$
obtained from $Q$ by adding a new vertex $i_0^*$ and an arrow
from $i_0^*$ to $i_0$. We extend the potential  from $R$ to $\widetilde{R}\supset R$ in a straightforward way.

For the representations of the modified algebra $\widetilde{R}$ we consider two stability conditions $Z=(z_i)_{i\in \widetilde{I}}$: one with $z_{i_0^*}$ very far to the left of all $(z_i)_{i\in I}$, and the one with $z_{i_0^*}$ on the right. For these two stability conditions we compare the  coefficients of the monomials in the Factorization Formula from Section 5.2 with dimension at $i_0^*$ equal to one. This leads to the following equality
$$\left(\sum_{\gamma\in \Z^I_{\geqslant 0}} [\mathbf{H}^\bullet(\M_\gamma^{(i_0)}/\G_\gamma,W_\gamma)]\,\mathbf{e}_\gamma \mathbf{e}_{i_0^*}\right)\cdot A=A\cdot  \mathbf{e}_{i_0^*}$$
where $\M_\gamma^{(i_0)}$ is the set of pairs $(E,v)$, where $E$ is a representation of $R$ with dimension vector $\gamma$ and a {\it cyclic} vector
$v\in E_{i_0}$. Therefore, we have:
$$A^{(i_0)}=\sum_{\gamma\in \Z^I_{\geqslant 0}} [\mathbf{H}^\bullet(\M_\gamma^{(i_0)}/\G_\gamma,W_\gamma)]\mathbf{e}_\gamma\,.$$
 The gauge group $\G_\gamma$ acts freely on $\M_\gamma^{(i_0)}$, hence cohomology of the quotient is  finite-dimensional.
$\blacksquare$

All the above considerations can be generalized to the case of a stability condition and of a  sector  $V\subset \mathbb{H}_+$. Finiteness of
coefficients of the series $A^{(i),cl}_V$ can not be proven in the same way as above, but it can be deduced from a very powerful integrality
 result proven in the next section.

\subsection{Examples}

First of all, any quiver with zero potential gives a kind of $q$-hypergeometric series:
$$A=\sum_{\gamma\in \Z^I_{\geqslant 0}}[\mathbf{H}^\bullet(B\G_\gamma)]\,\mathbf{e}_\gamma=\sum_{\gamma\in \Z^I_{\geqslant 0}}\frac{\left(-\mathbb{L}^{1/2}\right)^{\chi_Q(\gamma,\gamma)}}
{\prod_i (1- \mathbb{L})\dots(1- \mathbb{L}^{\gamma^i}) }\,
\widehat{\mathbf{e}}_\gamma\,.$$
Applying mutations from the next section, one can get further  examples with {\it non-zero} potentials.

One can easily write an explicit formula for the symplectomorphism  $T$ in the case of a quiver with zero potential. This is a  result of M.~Reineke (see \cite{Reineke3}). We give a sketch of the proof for completeness.
\begin{thm} For a quiver $Q$ with incidence matrix $(a_{ij})_{i,j\in I}$ and zero potential, series $A^{(i),cl}$ are algebraic functions, and they form a unique solution (in invertible formal series) of the system of polynomial equations
$$A^{(i),cl}=1+\mathbf{e}_i^{cl}\cdot \prod_{j\in I} \left(A^{(j),cl}\right)^{a_{ij}}.$$
\end{thm}
{\it Proof}. We have to calculate the Euler characteristic of the non-commutative Hilbert scheme, parameterizing representations with cyclic vector (see proof of Theorem 6). First, we parametrize the set of arrows between vertices $i,j\in I$ by $\alpha, 1\leqslant \alpha\leqslant a_{ij}$. 
Then one can apply the action of the torus $\mathbb{G}_m^{\sum_{ij}a_{ij}}$ acting by rescaling  arrows of the quiver.
 A fixed point of this action is either the trivial $0$-dimensional representation, or it corresponds to  a collection of $\sum_j a_{ij}$ torus-invariant representations with cyclic vectors $u_{j,\alpha}$ for $j\in I, 1\leqslant \alpha\leqslant a_{ij}$, such that the cyclic vector $v$ in the original representation is mapped to $u_{j,\alpha}$  by the arrow $i\to j$ corresponding to $\alpha$.
 The proof follows. $\blacksquare$

If the quiver is acyclic then all the series $A^{(i),cl}$ are rational functions, and the symplectomorphism $T$ gives a birational map.
 The same is true for a quiver $Q$ given by a clockwise oriented $N$-gon for any  $N\geqslant 3$.
 It was suggested in \cite{KS}, Section 8.4, that in such a situation it is natural to compose map $T$ with the antipodal involution of the classical symplectic
double torus.
  Using the relation of symplectomorphisms and mutations as in \cite{KS} one obtains the following result,
  which we formulate using the standard terminology of cluster algebras (see e.g. \cite{FomZel},\cite{FockGonch}).
\begin{thm} Let $Q$ be an acyclic quiver which is obtained by mutations from a Dynkin quiver. Then the birational map given by
$$x_i\mapsto \widetilde{x}_i:= \left(x_i \prod_j g_j^{a_{ij}-a_{ji}}\right)^{-1},\,\,\,y_i\mapsto\widetilde{y_i}:=\left(g_i y_i\right)^{-1}$$
where $(g_i)_{i\in I}$ is the unique solution of the system of equations
$$g_i=1+x_i\prod_{j} g_j^{a_{ij}}\in \Z[[(x_j)_{j\in I}]]\,,$$
induces the birational automorphisms of $\mathcal{X}$-torus with coordinates $(x_i)_{i\in I}$ and $\mathcal{A}$-torus with coordinates $(y_i)_{i\in I}$ of the order either $h+2$ or $(h+2)/2$ depending on the Dynkin diagram, where $h$ is the Coxeter number (cf.  with Coxeter automorphism from \cite{ZhuBin}).

Here  $\mathcal{X}$-torus is a quotient of the double torus under the projection $((x_i)_{i\in I},(y_i)_{i\in I})\mapsto (x_i)_{i\in I}$, and  $\mathcal{A}$-torus is a subtorus of the double torus given by the system of equations $x_i=\prod_{j\in I}y_j^{a_{ij}-a_{ji}}, i\in I$ (see \cite{FockGonch} for the details).

\end{thm}
Here we the use alleviated notation $x_i=\mathbf{e}^{cl}_i,\,y_i=\mathbf{f}^{cl}_i,\,g_i=A^{(i),cl}$.

The formula for the transformation from variables $(x_i)_{i\in I}$ to $(\widetilde{x}_i)_{i\in I}$ can be rewritten in a more symmetric form:
$$x_i=\frac{g_i-1}{\prod_j g_j^{a_{ij}}}\,,\,\,\,\,\widetilde{x}_i=\frac{\prod_j g_j^{a_{ji}}}{g_i-1}\,.$$

In the next example $Q=Q_3$ is the quiver with one vertex and three loops $x,y,z$, and the potential is $W=xyz-zyx$.

\begin{prp}
For $A:=A^{(Q_3,W)}$ we have the following formula:
$$A=\prod_{n,m\geqslant 1}(1-\mathbb{L}^{m-2}\,\widehat{\mathbf{e}}_1^n)^{-1}.$$
Therefore $A^{(1),cl}=\prod_{n\geqslant 1} (1-(\mathbf{e}_1^{cl})^n)^{-n}$ is the MacMahon series (cf. \cite{SzBrBehr}).
\end{prp}
{\it Proof.}
It is convenient to
use the duality and the approach from Section 4.8. Potential $W$ is linear in variable $z$, hence we have:
$$D(A)=\sum_{n\geqslant 0} [\mathbf{H}^\bullet_c(X_n)]\,  t^n, \,\,\,t:= D(\widehat{\mathbf{e}}_1),$$
where $X_n$ is the stack of $n$-dimensional representations of the polynomial algebra $\kk[x,y]$. The virtual dimension of $X_n$ is zero. 

In general, for any stack $S$ endowed with a map $w:S\to \Z_{\ge 0}$ (the target is a countable set understood as
 a stack), such that for any $k\ge 0$ the fiber $w^{-1}(k)$ is an Artin stack of finite type, we denote by
$\mathbf{H}^\bullet_c(S)$ an infinite series $\sum_{k\ge 0}[\mathbf{H}^\bullet_c(w^{-1}(k))] t^k$, 
 with coefficients of the $\lambda$-ring $\widehat{\mathcal M}$. The action of $\lambda$-operations on the variable $t$ is defined as $\lambda^{\geqslant 2} (t)=0$. Hence we can write $D(A)=[\mathbf{H}^\bullet_c(X)]$ where $X:=\sqcup_{n\ge 0} X_n,\,\,w(X_n):=\{n\}$.

It is a general property of $\lambda$-rings that any series $F$ of the form $1+\dots\in \widehat{\mathcal M}[[t]]$  
can be written uniquely as $\op{Sym}(G)$, where the series $G$ has zero constant term (see e.g. \cite{Knutson} or Section 6.1 below). We will denote $G$ by $\op{Log}_\lambda(F)$. Explicitly it is given by
 Harrison complex of a commutative algebra, see Section 6.6 for a related material.

It is convenient to introduce the following notation:
 for a stack $A$ endowed with a map $w:A\to \Z_{\ge 1}$ with fibers $w^{-1}(k),\,k\ge 1$ of finite type,
 and for another stack $B$ of fintie type, we denote
$$(pt\sqcup A)^B:=\sqcup_{n\ge 0} \left( A^n\times(B^n - \op{Diag}_n^B)\right)_{\op{Sym}_n}$$
by analogy with the Newton formula
$$(1+a)^b=\sum_{n\ge 0} \frac{a^n \cdot  b(b-1)\dots (b-n+1)}{n!}$$
Here $ \op{Diag}_n^B\subset B^n$ is the union of all diagonals, and 
the ``weight'' function $w$ is defined on $(pt\sqcup A)^B$ by the formula
$$w([a_1,\dots,a_n,b_1,\dots,b_n)])=\sum_i w(a_i)$$
 Then one can check easily that
$$\op{Log}_\lambda([\mathbf{H}^\bullet_c((pt\sqcup A)^B)])=[\mathbf{H}^\bullet_c(B)]\cdot \op{Log}_\lambda([\mathbf{H}^\bullet_c(pt\sqcup A)])$$

Now we apply this formula to our situation.
Let us denote by $X^{(0)}_+\subset \sqcup_{n\geqslant 0} X_n$ the  substack consisting of {\it nilpotent} representations of $\kk[x,y]$ (i.e. for any $n\ge 1$ both $x,y$ act in the $n$-dimensional representation by operators whose $n$-th power is zero). 
Then we have a {\it constructible} isomorphism of $\Z_{\ge 0}$-graded stacks
$$(pt \sqcup X^{(0)}_+)^{\mathbb{A}^2_\kk}\simeq X$$
 The reason for this identity is that any finite-dimensional $\kk[x,y]$-module
 is canonically isomorphic to a direct sum of shifts of non-trivial nilpotent modules by $(x\mapsto x+c_1,y\mapsto y +c_2)$ with where $(c_1,c_2)\in \mathbb{G}_{a,\kk}^2$ are different points of a finite subset of $\mathbb{A}^2_\kk $.  This implies the identity
$$\op{Log}_\lambda\left(\sum_{n\geqslant 0} [\mathbf{H}^\bullet_c(X_n)] \, t^n\right)=\mathbb{L}^2\op{Log}_\lambda\left(\sum_{n\geqslant 0} [\mathbf{H}^\bullet_c(X_n^{(0)})]\,  t^n \right),$$
where the factor $\mathbb{L}^2$ is equal to $[\mathbf{H}^\bullet_c(\mathbb{A}^2_\kk)]$ and $X_n^{(0)}=X^{(0)}\cap X_n$.

Similarly, we have
$$\op{Log}_\lambda\left(\sum_{n\geqslant 0} [\mathbf{H}^\bullet_c(X_n^{(1)})] \, t^n\right)=(\mathbb{L}-1)\op{Log}_\lambda\left(\sum_{n\geqslant 0} [\mathbf{H}^\bullet_c(X_n^{(0)})] \, t^n \right),$$
where $X^{(1)}$ denotes the stack of representations of $\kk[x,y]$ in which $x$ is nilpotent and $y$ is invertible. We claim that
 $$\sum_{n\geqslant 0} [\mathbf{H}^\bullet_c(X_n^{(1)})]\,  t^n=\prod_{n\geqslant 1} (1-t^n)^{-1}\,.$$
Indeed, the contribution of any stratum of $X^{(1)}_n$ corresponding to a given conjugacy class of the nilpotent operator $x$ (i.e. corresponding to a partition of $n$) is equal
 to $[\mathbf{H}^\bullet_c(G_x/G_x)]$, where $G_x$ is the normalizer of $x$ in $\op{GL}_\kk(n)$, and we consider the stack $G_x/G_x$
 given by the {\it adjoint} action of $G_x$ on itself. Group $G_x$ is the group of invertible elements in a finite-dimensional associative algebra, hence it is an extension of the product of general linear groups by a unipotent group. We can use Lemma 3 from Section 5.3 for the calculation of the class, and obtain the equality
$$[\mathbf{H}^\bullet_c(G_x/G_x)]=1\,.$$
Therefore the generating series for $ [\mathbf{H}^\bullet_c(X_n^{(1)})],\,n\ge 0$ is just the usual generating series for partitions. 
This implies the formula
$$\op{Log}_\lambda\left(\sum_{n\geqslant 0} [\mathbf{H}^\bullet_c(X_n^{(1)})]\,  t^n\right)=\sum_{n\geqslant 1} t^n\,,$$
and hence
$$D(A)=\op{Sym}\left(\mathbb{L}^2/(\mathbb{L}-1)\cdot \sum_{n\geqslant 1} t^n\right) =\prod_{n,m\geqslant 1}(1-\mathbb{L}^{2-m}\,t^n)^{-1}\,.$$
Applying the duality $D$ we obtain the desired formula for $A$. $\blacksquare$

In the final example quiver $Q$ has three vertices $I=\{1,2,3\}$, six arrows $\alpha_{12},\alpha_{23},\alpha_{31},\beta_{12},\beta_{23},\beta_{31}$ and potential
$$W=\alpha_{31}\cdot \alpha_{23}\cdot \alpha_{12}+\beta_{31}\cdot\beta_{23}\cdot\beta_{12}\,.$$
 This example is interesting because it is invariant under mutations (cf. \cite{KS}, Section 8.4).
Again, using the fact that $W$ is linear in $\alpha_{31},\beta_{31}$, we reduce the question to the calculation of $[\mathbf{H}^\bullet_c(S_\gamma)]$,
 where $S_\gamma,\,\gamma\in\Z^3_{\geqslant 0}$ is the stack of representations of $Q$ of dimension vector $\gamma$ with removed arrows $\alpha_{31},\beta_{31}$, and relations $\alpha_{23}\cdot \alpha_{12}=
\beta_{23}\cdot\beta_{12}=0$. It has the same cohomology as a similar stack for the quiver $Q'$ with five vertices $I={1,1',2,3,3'}$, six arrows
 $\alpha_{12},\alpha_{23},\beta_{1' 2},\beta_{23},\beta_{23'}, \delta_{11'},\delta_{33'}$ with relations $\alpha_{23}\cdot \alpha_{12}=0,
\beta_{23'}\cdot\beta_{1'2}=0$ and conditions that $\delta_{11'},\delta_{33'}$ are invertible. If one removes from $Q'$ the arrows $\delta_{11'},\delta_{33'}$ (which corresponds to the division by a simple factor $[\mathbf{H}^\bullet_c(GL(\gamma^1)\times GL(\gamma^2))]$) then one obtains a tame problem of linear algebra with 13 indecomposable objects and with dimension vectors
$$\left\{(\gamma^1,\gamma^{1'},\gamma^2,\gamma^3,\gamma^{3'})\,| \,\,\gamma^i\in \{0,1\},\,\,\gamma^1\gamma^2\gamma^3=\gamma^{1'}\gamma^2\gamma^{3'}=0\right\} .     $$ The result is a complicated sum (over 13 indices) of $q$-hypergeometric type.
 We do not know yet an explicit formula for the symplectomorphism $T$.

 \subsection{Mutations of quivers with potentials}

Here we study an application of the product formula to  the tilting (mutation) of quivers with polynomial potentials.
Results of this Section will  not be used in the rest of the paper, so the reader can jump to Section 6.1.

 Let $(Q,W)$ be a quiver with potential. We fix a vertex $i_0\in I$ of $Q$ which is loop-free (i.e. $a_{i_0 i_0}=0$ for the incidence matrix
$(a_{ij})_{i,j\in I}$). We will write the potential $W$  as a finite $\kk$-linear combination of cycles (in other words, cyclic  paths) $\sigma$ in $Q$
$$W=\sum_\sigma c_\sigma \sigma\,, $$
where $\kk$ is the ground field.
For any vertex $i\in I$ we have the corresponding cycle $(i)$ of length $0$ (the image of the projector corresponding to $i$).

We define below   the (right) {\bf mutation} $(Q',W')$ of $(Q,W)$ at the vertex $i_0$. Our goal in this section is to establish a relation between generating series $A^{(Q,W)}$ and $A^{(Q',W')}$ similar to the one from \cite{KS}, Section 8.4.

The definition of $(Q',W')$ is similar to the standard one from \cite{Zel} (in the context of quivers whose potentials are formal cyclic series).

1) The set of vertices of $Q'$ is the same set $I$.

2) The new set of arrows and new matrix $(a'_{ij})_{i,j\in I}$ are defined such as follows:

$a'_{i_0 i_0}=0$;

$a'_{i_0 j}=a_{j i_0}, a'_{j i_0}=a_{i_0 j}$ for any $j\ne i_0$ (in terms of arrows: we reverse each arrow $\alpha$ which has head or tail at $i_0$, i.e. replace $\alpha$ by a new arrow $\alpha^*$);

$a'_{j_1 j_2}=a_{j_1j_2}+a_{j_1i_0}a_{i_0j_2}$ for $j_1,j_2\ne i_0$ (in terms of arrows: for every pair of arrows $i_0 \stackrel{\beta}{\to}j_2$, $j_1\stackrel{\alpha}{\to}
i_0$,   we create a new arrow $j_1\stackrel{[\beta\alpha]}{\to}j_2$) \footnote{Since we use {\it left} $\kk Q$-modules our notation for the composition of arrows is compatible with the notation for the composition of linear operators.}.

3) The mutated potential $W'$ is defined as a sum of $3$ terms:
$$W'=W_1+W_2+W_3\,,$$
where
$$W_1=\sum_{j_1\stackrel{\alpha}{\to}
i_0 \stackrel{\beta}{\to}j_2}\beta^{\ast}\cdot [\beta\alpha]\cdot \alpha^{\ast}\,\,,$$

$$W_2=\sum_{ \sigma\ne (i_0)}c_{\sigma}\sigma^{mod}\,\,,$$

$$W_3=c_{(i_0)}\left(-(i_0)+\sum_{j\stackrel{\alpha}{\to}
i_0 } (j)\right)=c_{(i_0)}\left( -(i_0)+\sum_{j\in I }a_{ ji_0} (j)\right)\,\,.$$

Let us explain the notation in the formulas for $W_i, i=1,2,3$.

The summand $W_1$ consists of cubic terms generated by cycles from $i_0$ to $i_0$ of the form $\beta^{\ast}\cdot [\alpha\beta]\cdot \alpha^{\ast}$, and it can be thought  of as a ``Lagrange multiplier".

The summand $W_2$ is obtained from  $W$
 by modifying each cycle $\sigma$ (except $(i_0)$) to a cycle $\sigma^{mod}$ such as follows:
for each occurrence of $i_0$ in the cycle $\sigma$ (there might be several of them) we replace the consecutive two-arrow product $
\left(i_0 \stackrel{\beta}{\to}j' \right) \cdot \left(j\stackrel{\alpha}{\to}i_0\right)$ of $\sigma$ by an arrow
 $\stackrel{[\beta\alpha]}{j\to j'}$. New cycle is denoted by $\sigma^{mod}$, and it is taken with the same coefficient $c_{\sigma}$. In particular, if $\sigma$ does not contain $i_0$ then $\sigma^{mod}=\sigma$.
This modification procedure is not applicable only to the cycle $(i_0)$ of zero length.

The last summand $W_3$ can be thought of as a modification of the term in $W$ corresponding to the exceptional cycle $(i_0)$.
 There is also a version of $W_3$ for the {\it left} mutation in which the same sum is taken over all arrows with the tail (and not the head) at  $i_0$.

Our goal in this section is to establish a relation between generating series $A^{(Q,W)}$ and $A^{(Q',W')}$ similar to the one from \cite{KS}, Section 8.4. Each of the two series uniquely determines the other one.

Let us choose a central charge $Z:\Z_{\geqslant 0}^I\to \overline{\mathbb{H}}_+$, where $\mathbb{H}_+$ is the (open) upper-half plane. It is completely determined by its values on the standard basis vectors, so we  identify $Z$ with a collection of complex numbers $(z_i)_{i\in I}$ belonging to $\mathbb{H}_+$. We impose the condition that
$$\op{Arg}(z_{i_0})>\op{Arg}(z_{i})\,,\,\,\,\forall i\ne i_0\,\,.$$
 Denote by $\H^{(Q,W)}$ the Cohomological Hall algebra corresponding to $(Q,W)$ and by $\H^{(Q,W)}_V$ the Cohomological Hall  vector space  corresponding to a subsector $V\subset\mathbb{H}_+$, with the apex at the origin. For $V=l$ (a ray in $ \mathbb{H}_+$) this vector space is in fact a Cohomological Hall algebra, as we explained before. Let us introduce a ray $l_{0}=\R_{> 0}z_{i_0}$,   and a sector $V_{0,+}=\{z\in \mathbb{H}_+|\op{Arg}(z)<\op{Arg}(z_{i_0})\}$. The Factorization Formula implies
$$A^{(Q,W)}=A_{l_0}^{(Q,W)} A_{V_{0,+}}^{(Q,W)}\,.$$

We observe that $\H_{l_{0}}^{(Q,W)}$ corresponds to the category of representations of $Q$ supported at the vertex $i_0$. Since it is the same as the category of representations of the quiver $A_1$ with trivial potential, we can use the results of Section 2 and see that $\H_{l_0}^{(Q,W)}$ is a free exterior algebra, and
 $$ A_{l_0}^{(Q,W)}=(-\mathbf{e}_{i_0};\mathbb{L})_\infty  = (\mathbb{L}^{1/2}\widehat{ \mathbf{e}  }_{i_0};\mathbb{L})_\infty\, .  $$

Similarly, let us choose for the quiver $Q'$ a central charge $Z'$ such that for the corresponding collection of complex numbers $(z_i')_{i\in I}\in \mathbb{H}_+$  the inequalities
$$\op{Arg}(z_{i}')>\op{Arg}(z'_{i_0})\,,\,\,\,\forall i\ne i_0$$ hold. Consider
 the ray $l_0':=\R_{>0}z'_{i_0}$ and sector $V_{0,-}':=\{z\in \mathbb{H}_+|\op{Arg}(z)>\op{Arg}(z_{i_0}')\}$.
Then we use again the Factorization Formula
$$A^{(Q',W')}=A_{V_{0,-}'}^{(Q',W')}  A_{l_{0}'}^{(Q',W')}$$
 and the equality
$$ A_{l_{0}'}^{(Q',W')}= (\mathbb{L}^{1/2}\widehat{ \mathbf{e}  }_{i_0}';\mathbb{L})_\infty\,\,.$$

Comparison of the generating functions for $(Q,W)$ and $(Q',W')$ is based on the following result.

\begin{prp} Let $\gamma,\gamma'\in \Z^I$ are related as
$$(\gamma')^i=\gamma^i\mbox{ for }i\ne i_0\,,\,\,
(\gamma')^{i_0}=\sum_j a_{i_0 j} \gamma^j -\gamma^{i_0}\,.$$
Then there is an isomorphism of graded cohomology spaces
$$\H_{V_{0,+},\gamma}^{(Q,W),\bullet}\otimes \mathbb{T}^{\sum_j a_{j i_0} \gamma^j (\gamma')^{i_0} }
\simeq \H_{V_{0,-}',\gamma'}^{(Q',W'),\bullet}\,.$$
\end{prp}
\begin{cor}
After the identification of quantum tori
$$\widehat{\mathbf{e}}_\gamma\leftrightarrow \widehat{\mathbf{e}}'_{\gamma'}\,,\,\,\,
\mathbf{e}_\gamma\leftrightarrow \mathbb{L}^{\sum_j a_{j i_0} \gamma^j (\gamma')^{i_0}  } \mathbf{e}'_{\gamma'}\,, $$
where $\gamma$ and $\gamma'$ are related as in the Proposition, the
series $ A_{V_{0,+}}^{(Q,W)}$ and  $A_{V_{0,-}'}^{(Q',W')}$ coincide with each other.
\end{cor}
Corollary implies the comparison formula for $A^{(Q,W)}$ and $A^{(Q^{\prime},W^{\prime})}$ which coincides with the one from Section 8.4 in \cite{KS} (called there the property 3) of
the generating series $\bf{E}_Q$). 

{\it Proof of Proposition.} By definition,
$$\H_{V_{0,+},\gamma}^{(Q,W),\bullet}=\mathbf{H}_{\G_\gamma}(\M_{V_{0,+},\gamma}^Q,W_\gamma)\,,$$
where $\M_{V_{0,+},\gamma}^Q \in \M_\gamma^Q$ is a $\G_\gamma$-invariant Zariski open subset consisting of representations of $Q$ which do not have
 a non-trivial subrepresentation supported at the vertex $i_0$. In other words, a representation $E=\oplus_i E_i $ (over $\overline{\kk}$)
 belongs to $\M_{V_{0,+},\gamma}^Q$ if and only if the homomorphism
$$E_{i_0}\to \oplus_{j\ne i_0}\,\oplus_{i_0\stackrel{\beta}{\to}j} E_{j}$$
given by the  direct sum of all arrows with head at $i_0$, is a monomorphism.

Similarly, we have
$$\H_{V_{0,-}',\gamma'}^{(Q',W'),\bullet}=\mathbf{H}_{\G_{\gamma'}}(\M_{V_{0,-}',\gamma'}^{Q'},W_{\gamma'}')\,,$$
where $\M_{V_{0,-}',\gamma'}^{Q'} \subset \M_{\gamma'}^{Q'}$ is given by the condition that
$$\oplus_{j\ne i_0}\,\oplus_{i_0 \stackrel{\beta^*}{\leftarrow}j} E_{j}'\to E_{i_0}'$$
is an epimorphism.

Let us denote by $Q''$ the quiver obtained from $Q'$ by removing all arrows $j\stackrel{\alpha^*}{\leftarrow}i_0$ with tail at $i_0$.
The condition on a representation $E'$ of $Q'$ to belong to $\M_{V_{0,-}',\gamma'}^{Q'}$ depends only on the restriction of the action of $Q'$ to
 $Q''$, and we get a fibration by affine spaces
$$\pi: \M_{V_{0,-}',\gamma'}^{Q'} \to \M_{V_{0,-}',\gamma'}^{Q''}\,.$$

Potential $W'_{\gamma'}$ is at most linear along fibers of $\pi$. Hence, the closed subscheme of $ \M_{V_{0,-}',\gamma'}^{Q'}$ consisting
 of points at which $W'_{\gamma'}$ has zero derivative along the fiber of $\pi$, is the pullback by $\pi$ of a closed subscheme $Z_{\gamma'}$.
It is straightforward to see (considering terms of type $W_1$ in the formula for $W'$)
that $Z_{\gamma'}\subset \M_{V_{0,-}',\gamma'}^{Q''}$ is
the closed $\G_{\gamma'}$-invariant subscheme given by the condition
$$\forall k\ne i_0\,, \,\forall k\stackrel{\alpha^*}{\leftarrow }i_0\,,\,\,\,\,\sum_{j\ne i_0} \sum_{i_0\stackrel{\beta^*}{\leftarrow} j}\beta^*\cdot [\beta \alpha]=0 \in \op{Hom}(E_k,E_{i_0}')\,.$$

In other words, $Z_{\gamma'}$ is given by the equation $\partial W_1/\partial \alpha^*=0, \,\,\,\forall k\stackrel{\alpha^*}{\leftarrow} i_0$.

\begin{lmm} Scheme $Z_{\gamma'}$ is smooth, and we have a natural equivalence of stacks
$$\M_{V_{0,+},\gamma}^Q/\G_{\gamma}\simeq Z_{\gamma'}/\G_{\gamma'}$$
identifying restrictions of potentials $W_\gamma$ and $W'_{\gamma'}$ respectively.
\end{lmm}

{\it Proof of Lemma.} This is a version of the classical Gelfand-Ponomarev correspondence. Namely, if $E$ is a representation
 of $Q$ which belongs to  $\M_{V_{0,+},\gamma}^Q/\G_{\gamma}$\footnote{We will abuse the notation identifying a representation with the corresponding point of the stack.}, then we associate with it a representation $E'$ of $Q''$ by
$$E'_j:=E_j\mbox{ for } j\ne i_0\,,\,\,\,E'_{i_0}:=\op{Coker}\left( E_{i_0}\mono \oplus_{j\ne i_0}\oplus_{i_0\stackrel{\beta}{\to} j} E_{j}\right), $$
and the action of arrows in $Q''$ given in the following way.
First, for any ``old'' arrow $j_1\to j_2$ with $j_1,j_2\ne i_0$ its action in $E$ and $E'$ is the
same.
The action of arrow $\beta^*:E_{j_0}'=E_{j_0}\to E'_{i_0}$ for any $i_0\stackrel{\beta}{\to} {j_0}$ in the original quiver $Q$ is the composition
$$E_{j_0}\mono  \oplus_{j\ne i_0}\oplus_{i_0\stackrel{\beta}{\to} j} E_{j}\epi   E'_{i_0} $$
of the natural inclusion of the direct summand, and of the projection to the cokernel.
The action of $[\beta\alpha]$ in $E'$ is defined as the composition of arrows corresponding to $\alpha$ and $\beta $ in $Q$.
The dimension   vector of $E'$ is $\gamma'$ as defined in the Proposition.

 It is easy to see that $E'$ belongs to $M_{V_{0,-}'',\gamma'}^{Q'}/\G_{\gamma'}$, and
for any arrow $\alpha:j\to i_0$ in $Q$ the equation
 $$\sum_{\beta}\beta^*\cdot  [\beta\alpha]=0$$
is satisfied. Hence $E'$ belongs to $Z_{\gamma'}/\G_{\gamma'}$.

Conversely, for any representation $E'$ of $Q''$ lying in $Z_{\gamma'}/\G_{\gamma'}$ we define the corresponding representation $E$ of $Q$ by
$$E_j:=E_j'\mbox{ for } j\ne i_0,\,\,\,E_{i_0}:=\op{Ker}\left( \oplus_{j\ne i_0}\oplus_{i_0\stackrel{\beta^*}{\leftarrow} j} E_{j}'\epi E_{i_0}'\right), $$
and the action of arrows of $Q$ is determined uniquely in such a way that $E'$ is obtained from $E$ by the formulas given at the beginning of the proof.

Hence, we conclude that {\it stacks} $\M_{V_{0,+},\gamma}^Q/\G_{\gamma}$ and $Z_{\gamma'}/\G_{\gamma'}$ are naturally equivalent, and therefore the
 {\it scheme} $Z_{\gamma'}$ is smooth (because $\M_{V_{0,+},\gamma}^Q/\G_{\gamma}$ is a smooth stack).

A straightforward check shows that on $\pi^{-1}(Z_{\gamma'})$ the component of $W'_{\gamma'}$ corresponding to the term $W_1$ vanishes, the component corresponding to $W_2$
matches $W_\gamma-c_{i_0}\op{dim}{E_{i_0}}$, and the value of constant given by $W_3$ coincides with $c_{i_0}\op{dim}{E_{i_0}}$.
 $\blacksquare$

In order to finish the proof of Proposition we  apply  Proposition 6 from Section 4.8 with the following choices of spaces and functions:
$$
X=\M_{V_{0,-}',\gamma'}^{Q'}\,,\, Y=\M_{V_{0,-}',\gamma'}^{Q''}\,,\,G=\G_{\gamma'}\,,\,Z=Z_{\gamma'}\,,\,f=W'_{\gamma'}\,,$$
together with Poincar\'e duality for the smooth stack $Z_{\gamma'}/\G_{\gamma'}$. The shift  
$$d=\sum_j a_{j i_0} \gamma^j (\gamma')^{i_0} $$
is the dimension of fibers of  fibration $\pi$.
Proposition is proven. $\blacksquare$

\section{Factorization systems and integrality of exponents}

\subsection{Admissible series and integrality properties}

We will work in the framework of $\lambda$-rings (see e.g.\cite{Groth}, \cite{Knutson}, \cite{Macdon}).
Recall that for such a ring $B$ we have a collection of operations $\lambda^n: B\to B, n\geqslant 0$ satisfying the natural properties of the operations of wedge power on the $K_0$-group of the Tannakian category $G-mod$ of finite-dimensional representations of an affine algebraic group $G$ (see loc. cit. for precise definition).
\begin{dfn} An element $x$ of  a $\lambda$-ring is called {\bf line} element if  $x$ is invertible and
$$\lambda^i(x)=0\,,\,\,\,\,\forall i\geqslant 2\,.$$
\end{dfn}
For example, for the $K_0$-ring of the tensor category of finite-dimensional representations of a pro-affine algebraic group over a field, the class
of any one-dimensional representation is a line element. The set of line elements in $B$ is a subgroup of $B^\times$.

 If $B$ is $\lambda$-ring and $t$ is a variable (which we treat as a line element) then
$$B[t],\,\,B[t,t^{-1}],\,\,B[[t]],\,\,B((t))$$
are again $\lambda$-rings.

\begin{dfn} For a $\lambda$-ring $B$ and a finite set $I$, a series
$$F\in B((q^{1/2}))[[(x_i)_{i\in I}]]\,,$$
where $q^{1/2},(x_i)_{i\in I}$ are line elements (variables) is called {\bf admissible} if it has a form
$$F=\op{Sym}\left(\sum_{\gamma\in \Z_{\geqslant 0}^I-\{0\}} \left(f_\gamma\cdot \sum_{n\geqslant 0} q^n\cdot\prod_{i\in I} x_i^{\gamma^i}\right) \right),$$
where $f_\gamma\in B[q^{\pm 1/2}]$ for any $\gamma=(\gamma^i)_{i\in I} \in \Z_{\geqslant 0}^I-\{0\}$, and
$$\op{Sym}(b):=\sum_{n\geqslant 0}\op{Sym}^n(b)=\sum_{n\geqslant 0} (-1)^n \lambda^n(-b)=\left(\sum_{n\geqslant 0} (-1)^n \lambda^n(b)\right)^{-1}.$$
\end{dfn}
Coefficients of an admissible series belong to the ring
$$B\,[q^{\pm 1/2},\left((1-q^n)^{-1}\right)_{n=1,2,\dots}]\subset B((q^{1/2}))\,.$$
In the case $B=\Z$ and $|I|=1$
a series $F$ in one variable $x$ with coefficients in $\Z((q^{1/2}))$ is  admissible if and only if it can be represented as a product
$$F=\prod_{n\geqslant 1} \prod_{i\in \Z}\left( q^{i/2}x^n;q\right)_\infty^{c(n,i)}\in 1+x \cdot \Z((q^{1/2}))[[x]]\,,$$
where $c(n,i)\in \Z$ for all $i,n$, and for any given $n$ we have $c(n,i)=0$ for $|i|\gg 0$. Equivalently,
$$F=\op{exp}\left(-\sum_{n,m\geqslant 1}\frac{f_n(q^{m/2})}{m(1-q^m)}x^{nm} \right)\,,$$
where $f_n=f_n(t)$ belongs to $\Z[t^{\pm 1}]$ for all $n\geqslant 1$.
The equivalence of two descriptions follows from the identity
$$\op{log}\left( q^{i/2}x^n;q\right)_\infty=-\sum_{m\geqslant 1} \frac{ \left(q^{m/2}\right)^i}{m(1-q^m)} x^{nm}\,.$$
It is easy to see that any  series which belongs to the multiplicative group
$$1+x\cdot \Z[q^{\pm 1/2}]\,[[x]]$$
is admissible.
Indeed, such a series can be written uniquely as a product
$$\prod_{n\geqslant 1} \prod_{i\in \Z}(1-q^{i/2} x^n)^{b(n,i)}=\prod_{n\geqslant 1}\prod_{i\in \Z}\left( \left(q^{i/2}x^n;q \right)_\infty^{b(n,i)}
\left(q^{i/2+1}x^n;q \right)_\infty^{-b(n,i)}\right),$$
where $b(n,i)\in \Z\,\,\,\forall i,n$, and for any given $n$ we have $b(n,i)=0$ for $|i|\gg 0$.

 Admissibility of a series $F=F(x;q^{1/2})\in \Z((q^{\pm 1/2}))[[x]]$ implies certain divisibility properties. Namely, let us define a new series by  the formula
 $$G(x;q^{1/2}):=\frac{F(x;q^{1/2})}{F(qx;q^{1/2})}\in \Z[q^{\pm 1/2}]\,[[x]]\,.$$
Then the evaluation at $q^{1/2}=1$ of the series $G$ is of the form
$$G(x;1):=\op{lim}_{q^{1/2}\to 1}G(x;q^{1/2})=\prod_{n\geqslant 1} (1-x^n)^{n c(n)}\in 1+x\Z[[x]]\,,\,\,\,c(n)\in \Z\,.$$
Obviously, admissible series form a group under multiplication. The next result shows that the admissibility is preserved by a non-trivial transformation.

\begin{thm}For   a given symmetric integer matrix $\mathcal{B}=(b_{ij})_{i,j\in I}$, a series
$$F=\sum_{\gamma \in \Z_{\geqslant 0}^I} a_\gamma x^\gamma\in B((q^{1/2}))[[(x_i)_{i\in I}]]\,,\,\,\,x^\gamma:=\prod_I x_i^{\gamma^i},\,\,a_\gamma\in B((q^{1/2}))$$
 is admissible if and only if
$$\widetilde{F}:=\sum_{\gamma \in \Z_{\geqslant 0}^I} \left(-q^{1/2}\right)^{\sum_{ij} b_{ij}\gamma^i\gamma^j
} a_\gamma x^\gamma$$
is admissible.
\end{thm}

Proof will be given in Section \ref{sect:endofproof}.

This theorem implies that one can define admissible series in {\it quantum} variables. Namely, let us assume that we are given
 a skew-symmetric integer bilinear form
 $\langle\,,\,\rangle:\Z^I\otimes \Z^I\epi \wedge^2(\Z^I)\to \Z$. We define the quantum torus associated with $\langle\,,\,\rangle$
as an associative  algebra over $\Z((q^{1/2}))$
generated by the elements $\widehat{\mathbf{e}}_\gamma,\,\,\gamma\in \Z^I$ subject to the relations given in Section 5.3 (with the identification $q^{1/2}=\mathbb{L}^{1/2}$):
 $$\widehat{\mathbf{e}}_{\gamma_1}\widehat{\mathbf{e}}_{\gamma_2}=(-q^{1/2})^{-\langle \gamma_1,\gamma_2\rangle}\widehat{\mathbf{e}}_{\gamma_1+\gamma_2}\,,\,\,\,\widehat{\mathbf{e}}_{0}=1\,.$$

Let us choose an ordered basis $(\gamma_1,\dots,\gamma_k),\,k=|I|$ of the lattice $\Z^I$.
 To a ``quantum" series
$$F=\sum_{\gamma \in \Z_{\geqslant 0}^I} a_\gamma\cdot \widehat{\mathbf{e}}_\gamma,\,\,\,a_\gamma\in B((q^{1/2}))$$
we can assign the corresponding ``classical" series in commuting variables $(x_i)_{i\in I}$ in the following way.
 First, we rewrite $F$ as a series in ordered monomials
$$F=\sum_{\gamma=\sum_i n_i \gamma_i \in \Z_{\geqslant 0}^I} \widetilde{a}_\gamma \cdot \widehat{\mathbf{e}}_{\gamma_1}^{n_1}\dots \widehat{\mathbf{e}}_{\gamma_k}^{n_k}
,\,\,\,\widetilde{a}_\gamma\in B((q^{1/2}))$$
and then replace noncommuting variables by commuting ones:
$$F_{cl}:=\sum_{\gamma \in \Z_{\geqslant 0}^I} \widetilde{a}_\gamma \cdot x^\gamma\in B((q^{1/2}))[[(x_i)_{i\in I}]]\,.$$
For different choices of ordered bases the coefficients $\widetilde{a}_\gamma$ differ by the multiplication by a power of $(-q^{1/2})$ with the exponent given by an integer symmetric bilinear form evaluated at $\gamma\otimes \gamma$.

\begin{dfn} We say that the quantum series $F$ is {\bf quantum admissible} if for some (and hence by Theorem 9 any) choice of the ordered basis the corresponding classical series $F_{cl}$ is admissible.

\end{dfn}

\begin{prp}

Let us choose a generic central charge $Z:{\Z}_{\geqslant 0}^I\to \mathbb{H}_+$. Then the set of quantum admissible series coincides with the set of products
$$\prod^{\curvearrowright}_{l=\R_{>0}Z(\gamma_0)}F_l(\widehat{\mathbf{e}}_{\gamma_0})\,,$$
where the product in the clockwise order is taken over all rays generated by primitive vectors $\gamma_0\in \Z_{\geqslant 0}^I$, and $F_l(t)$ is an admissible series in one variable.

\end{prp}

{\it Proof.} Denote by
${\cal R}_+$ the ring of series in $(\widehat{\mathbf{e}}_\gamma)_{\gamma\in\Z^I_{\geqslant 0}}$ with coefficients in $ B((q^{1/2}))$, and
by $\mathfrak{m}\subset {\cal R}_+$ the maximal ideal generated by $\widehat{\mathbf{e}}_\gamma,\, \gamma\ne 0$. Notice that the definition of quantum admissible series makes sense in all quotients ${\cal R}_+/\mathfrak{m}^N, N=1,2,\dots$, and a series $F$ is admissible if and only if all its
 images in ${\cal R}_+/\mathfrak{m}^N, \,N\geqslant 1$ are  admissible.  Because of that we are going to assume that $F\in {\cal R}_+/\mathfrak{m}^N$ is a
 quantum admissible truncated series for a {\it given} $N$. Hence we can write
$$F=\sum_{|\gamma|<N}a_{\gamma}\widehat{\mathbf{e}}_{\gamma}=1+\dots\mbox{ where }|\gamma|:=\sum_i \gamma^i\,.$$
If $F\ne 1$ then there exists a primitive vector $\gamma_0\in \Z_{\geqslant 0}^I$ with $|\gamma_0|<N$ such that for some $k\geqslant 1$ one has $a_{k\gamma_0}\ne 0$, and $Arg(Z(\gamma_0))$ is the largest among all $Arg(Z(\gamma))$ for which $a_{\gamma}\ne 0, \gamma\ne 0$. Quantum admissibility of $F$ implies quantum admissibility of the following truncated series
 in one variable
$$F_{\gamma_0}:=\sum_{|k\gamma_0|<N}a_{k\gamma_0}\widehat{\mathbf{e}}_{k\gamma_0}\,.$$
Let us choose an ordered basis $\gamma_1,\dots,\gamma_k, \,k=|I|$ of $\Z^I$ such that $\gamma_1=\gamma_0$.
Notice that the multiplication from the left by any quantum admissible series $f(\widehat{\mathbf{e}}_{\gamma_0})$ commutes with the correspondence $F\mapsto F_{cl}$. It follows that the series $(F_{\gamma_0})^{-1}F$ is quantum admissible. Then we proceed by induction in the decreasing order of $Arg(Z(\gamma))$ for all non-zero primitive vectors  $\gamma$ which appear in the series $F$ with non-zero coefficient $a_{\gamma}$. Since for a fixed $N$ we have only finitely many rays,  this proves that $F$ is quantum admissible modulo $\mathfrak{m}^N$. Then, as we said above, we take the projective limit $N\to \infty$. This proves that $F$ factorizes into a clockwise ordered product over the set of rays of quantum admissible series in one variable. The converse statement (such a product is quantum admissible) can be proved in a similar way. $\blacksquare$

\begin{prp} The set of quantum admissible series
forms a subgroup under multiplication.

\end{prp}

{\it Proof.} As in the above proof we will consider all series modulo $\mathfrak{m}^N$. Let $F_1,F_2$ be quantum admissible series. Let us choose a generic central charge $Z$. Then, by the previous Proposition we have
$$F_1=\prod_l^{\curvearrowright}F_{1,l}\,,$$
where the product is taken over all rays as above, and $F_{1,l}$ are quantum admissible series in one variable. It suffices to prove that for any ray $l$ and any quantum admissible series $F$ the product $F_{1,l}F$ is quantum admissible (then we proceed by induction). If $l=\R_{>0}Z(\gamma_0)$ then we choose an ordered basis $\gamma_1,\dots,\gamma_k$ of $\Z^I$ such that $\gamma_1=\gamma_0$. As in the proof of the previous Proposition we conclude that the multiplication  from the left by any quantum admissible series in $\widehat{\mathbf{e}}_{\gamma_0}$ preserves the set of quantum admissible series. This implies the desired statement. Similarly one proves that if $F$ is quantum admissible then its inverse $F^{-1}$ is quantum admissible as well. $\blacksquare$.

\begin{cor} For quantum variables $yx=qxy$ the collection of elements
$$\prod^{\curvearrowright}_{(a,b)\in \Z^2_{\geqslant 0}-\{0\}} \prod_{|k|\leqslant const(a,b)}\left((-1)^{ab}q^{k/2}  x^a y^b;q\right)_\infty^{c(a,b;k)}\,\,$$
with $c(a,b;k)\in \Z$,
is closed under the product.
\end{cor}

The corresponding group can be called {\it the quantum tropical vertex group} (cf. \cite{GPS}) since for a quantum admissible series $F$ the automorphism $Ad(F)$ in the limit $q^{1/2}\to 1$ gives rise to a formal symplectomorphism of the symplectic torus considered in \cite{KS2}, \cite{GPS}.

 \subsection{Admissibility of generating series and motivic Donaldson-Thomas invariants}

Let $R$ be an $I$-bigraded smooth algebra endowed with a potential $W$ and satisfying the Assumption of Section 3.3. Let us choose
a central charge  $Z:I\to \mathbb{H}_+$. We also fix a cohomology theory $\mathbf{H}^\bullet$ with values in the category $\mathcal{C}^{\mathbb{Z}-gr}$ associated with a Tannakian category $\mathcal{C}$ (see Section 3.1).
Also we choose a degree $+1$ tensor square root $\mathbb{T}^{\otimes 1/2}$ of $\mathbb{T}\in Ob(\mathcal{C}^{\mathbb{Z}-gr})$ (see  Section 2.7).

Let us denote by $B_{\mathcal{C}}$ the $\lambda$-ring which is the $K_0$-ring  of the subcategory of $\mathcal{C}$ consisting of objects of weight $0$. Hence we have:
 $$K_0(\mathcal{C})=B_{\mathcal{C}}[q^{\pm 1/2}]\,,\,\,\,q^{1/2}:=\mathbb{L}^{1/2}=[\mathbb{T}^{\otimes 1/2}[1]]\,.$$

\begin{thm} For any sector $V\subset \mathbb{H}_+$ the generating series $A_V$ is quantum admissible.
\end{thm}

The proof will be given in Section 6.9.

\begin{cor} Let us assume that for some ray $l=\op{exp}(i \phi)\R_{>0}\subset \mathbb{H}_+$ the restriction of the form $\chi_R$ to the
sublattice $\Gamma_l:=Z^{-1}(\op{exp}(i \phi)\R)\subset \Z^I$ is symmetric. Then for some  elements $\Omega^{mot}(\gamma)\in K_0(\mathcal{C})$ the following formula holds
$$\sum_{\gamma \in Z^{-1}(l)}[\H_{l,\gamma}^{mod}]x^\gamma=\op{Sym}\left(\sum_{\gamma \in Z^{-1}(l)} \Omega^{mot}(\gamma)\cdot
[\mathbf{H}^\bullet(\mathbb{P}^\infty)]\cdot x^\gamma\right).$$

\end{cor}
\begin{dfn} Under the above assumptions the element $\Omega^{mot}(\gamma)$ is called {\bf motivic Donaldson-Thomas invariant}
of the pair $(R,W)$, stability function $Z$ and dimension vector $\gamma$.
\end{dfn}

We can apply homomorphism of $\lambda$-rings $K_0(\mathcal{C})\to \Z[q^{\pm 1/2}]$ and get the so-called {\it refined} (or quantum) Donaldson-Thomas invariants,
 or we can apply Euler characteristic $K_0(\mathcal{C})\to \Z$ (i.e. evaluate at $q^{1/2}=1$) and obtain numerical Donaldson-Thomas invariants (cf. \cite{KS}).

\begin{rmk} The reader should not mix motivic DT-series and motivic DT-invariants $\Omega^{mot}(\gamma)$ with motivic DT-invariants of $3CY$ categories introduced in \cite{KS}, although the latter should be better  called motivic DT-series. In fact motivic DT-series for the {\bf $3CY$ category} which is generated by critical points of $W$ is not necessarily related to $\Omega^{mot}(\gamma)$ and to the motivic DT-series of the {\bf algebra with potential $(R,W)$}. This reflects the difference between the rapid decay cohomology used in the definition of COHA given above and the critical cohomology discussed in Section 7 and used for the definition of ``critical COHA''. The DT-series for the latter are related to those from \cite{KS}.
\end{rmk}

Recall that the assumption that $(\chi_R)_{|\Gamma_l}$ is symmetric (which holds e.g. for any generic central charge $Z$) implies that
 we have  an ordinary associative (super) algebra $\mathcal{H}_{l}^{mod}$ graded by $\Gamma_l\oplus \Z$, where the second $\Z$-grading is the cohomological grading (which differs from the cohomological grading
for non-modified COHA by $\chi_R(\gamma,\gamma),\,\gamma\in \Gamma_l$). The parity is given by  the cohomological grading modulo $2$.

We pointed out  at the end of Section 5.2 an analogy between the product decomposition by slopes from
Theorem 5 and the Poincar\'e-Birkhoff-Witt theorem.
 In particular, the algebra $\mathcal{H}_{l}^{mod}$ can be thought of as an analog of the universal enveloping algebra
 for some  Lie (super) algebra.

\begin{question} Is it true that $\mathcal{H}_{l}^{mod}$ is the universal enveloping algebra $U\mathfrak{g}$, where $\mathfrak{g}=\oplus_{\gamma\in \Gamma_l} \mathfrak{g}_{\gamma}$
 is an infinite-dimensional $\Gamma_l\oplus \Z$-graded Lie super algebra in $\mathcal{C}$,
 which (as $\Z$-graded ind-object of $\mathcal{C}$) has the form $\mathfrak{g}=\mathfrak{g}^{\op{prim}}\otimes (\oplus_{n\geqslant 0}\mathbb{T}^{\otimes n})$, where $\mathfrak{g}^{\op{prim}}$ has finite-dimensional graded components
for each grading $\gamma\in \Gamma_l\cap \Z_{\geqslant 0}^I$?
\end{question}

This question is closely related with Conjecture 1 from Section 2.6.
The main evidence in favor of the  positive answer to the above Question is Theorem 10. In all examples which we have studied so far the algebras
 $\mathcal{H}_{l}^{mod}$ for generic $Z$ were in fact supercommutative. Hence  the hypothetical Lie super algebra $\mathfrak{g}$ was abelian.

\subsection{Equivariant cohomology with respect to maximal tori}

 Let $R$ be a smooth $I$-bigraded algebra over $\kk$, and $\mathbf{H}^\bullet$ be a cohomology theory with values in the tensor category of $\Z$-graded vector spaces over the field $K,\,\,\op{char}(K)=0$
 used in the definition of the Cohomological Hall algebra.
 Then for every $\gamma\in \Z_{\geqslant 0}^I$ the space $\mathcal{H}_{\gamma}$ is a module over ${\bf H}^\bullet(\op{B}\G_\gamma)$.

 Let us introduce a more refined version of the graded component, namely
$$\mathcal{H}'_{\gamma}:={\bf H}^\bullet_{\mathsf{T}_\gamma}(\M_{\gamma})\,,$$
where
$$\mathsf{T}_\gamma= \prod_{i\in I}\left(\mathbb{G}_{m,\kk}\right)^{\gamma^i}\subset \G_\gamma=\prod_{i\in I} \op{GL}_{\kk}(\gamma_i)$$
is the canonical maximal algebraic torus in $\G_\gamma$. The Weyl group of $\G_\gamma$, which is the product of symmetric groups
 $$\op{Sym}_\gamma:=\prod_{i\in I} \op{Sym}_{\gamma^i}$$
acts on $\mathsf{T}_\gamma$ and on $\M_{\gamma}$. Also $\mathcal{H}'_\gamma$ is a module over the ring ${\bf H}^\bullet(\op{B}\mathsf{T}_\gamma)$,
 which is the polynomial ring with $\sum_i \gamma^i$ generators in degree $2$.
We see that $\mathcal{H}'_\gamma$, treated as a super vector space over $K$ (we reduce the cohomological $\Z$-grading on $\mathcal{H}'_\gamma$ to $\Z/2\Z$-grading), is the space of sections of a $\op{Sym}_\gamma$-equivariant super quasi-coherent sheaf $\mathcal{F}_\gamma$
 on
$$\mathbb{A}^\gamma_K:=\prod_{i\in I} \mathbb{A}^{\gamma^i}_K\,\,.$$
The latter affine space can be considered as the space of configurations of $\sum_i \gamma^i$ points in the affine line $\mathbb{A}^1_K$ over $K$.
The original space $\mathcal{H}_\gamma$ coincides with the space of $\op{Sym}_\gamma$-invariants in $\mathcal{H}'_\gamma$ (because ${\bf H}^\bullet(\G_\gamma/(\op
{Sym}_\gamma\ltimes \mathsf{T}_\gamma))={\bf H}^\bullet(pt)$),
and it is a module over
$${\bf H}^\bullet(\op{B}\G_\gamma)=\left( {\bf H}^\bullet(\op{B}\mathsf{T}_\gamma) \right)^{\op{Sym}_\gamma   }\,\,.$$

\begin{lmm} Space $\mathcal{H}'_\gamma$ is a finitely generated ${\bf H}^\bullet(\op{B}\mathsf{T}_\gamma)$-module.
\end{lmm}
Equivalently, the super quasi-coherent sheaf $\mathcal{F}_\gamma$ is coherent. Lemma implies also that $\mathcal{H}_\gamma$
is a finitely generated ${\bf H}^\bullet(\op{B}\G_\gamma)$-module.

{\it Proof.} Notice that the stabilizer  in $\mathsf{T}_\gamma$ of any point in  $\M_{\gamma}$ is a connected subtorus of
 $\mathsf{T}_\gamma$, and for a given $\gamma$ there are finitely many different stabilizers $\mathsf{T}'\subset \mathsf{T}_\gamma$ .
  Consider the finite increasing filtration $\M_{\gamma}^{\leqslant k}$ of $\M_{\gamma}$ by $\mathsf{T}_\gamma$-invariant
Zariski open subspaces consisting of points for which the dimension of the stabilizer is at most $k,\,\,
k=0,\dots,\op{dim} \mathsf{T}_\gamma$. Let us calculate  first term of the corresponding spectral sequence.
 The contribution of the locus $\M_{\gamma}^{=\mathsf{T}'}$ of points in $\M_{\gamma}$ with stabilizer exactly equal to some $\mathsf{T}'$
 is
$${\bf H}^\bullet(\op{B}\mathsf{T}')\otimes {\bf H}^\bullet((\M_{\gamma}^{=\mathsf{T}'}/\mathsf{T}_\gamma)^{coarse})\otimes \mathbb{T}^{\otimes
\op{codim} \M_\gamma^{=\mathsf{T}'}}$$
by Thom isomorphism. Here $(\M_{\gamma}^{=\mathsf{T}'}/\mathsf{T}_\gamma)^{coarse}$ is the scheme-theoretical quotient, not a stack. This is obviously a finitely generated module because $\op{dim}{\bf H}^\bullet((\M_{\gamma}^{=\mathsf{T}'}/\mathsf{T}_\gamma)^{coarse})<\infty$.  In terms of sheaves it corresponds to a trivial
 super vector bundle of finite rank over an affine subspace
$$\op{Spec}{\bf H}^\bullet(\op{B}\mathsf{T}')\subset \op{Spec} {\bf H}^\bullet(\op{B}\mathsf{T}_\gamma)\,\,.$$
Lemma is proven.
$\blacksquare$

Obviously, all the above considerations extend to the case when we have a potential $W\in R/[R,R]$, and also when we are given a stability function $Z$ and a  sector $V\subset \mathbb{H}_+$.
We denote by $\mathcal{F}_{V,\gamma}$ the corresponding super coherent sheaves .
In the next section we axiomatize some algebraic properties of the collection of sheaves $(\mathcal{F}_{V,\gamma})_{\gamma\in \Z_{\geqslant 0}^I-0}
$ (the case $\gamma=0$ is trivial).

\subsection{Factorization systems}

Let $G$ be a pro-algebraic group over a  field $K,\,\,\op{char}(K)=0$, endowed with a homomorphism
$$t^{1/2}:G\to\mathbb{G}_{m,K}\,.$$
For us the basic example will be the motivic Galois group (or its double cover) with $t^{1/2}$ corresponding to
 the square root of Tate motive $K(-1/2)$.
 We define
             $$\widetilde{G}:=G\ltimes \mathbb{G}_a\,,$$
             where $G$ acts on the additive group $\mathbb{G}_a=\mathbb{G}_{a,K}$ via the homomorphism
             $$G\stackrel{t^{-1}}{\to}\mathbb{G}_m=\op{Aut}\mathbb{G}_a\,\,.$$

	     We denote by $\mathbb{T}_0^{1/2}$ the representation of $G$ corresponding to the character
$t^{1/2}$. We set $\mathbb{T}_0=(\mathbb{T}_0^{1/2})^{\otimes 2}$. Also, we denote by $q^{1/2}$ the class of $\mathbb{T}_0^{1/2}$ in $K_0(G-mod)$.

             Group scheme $\widetilde{G}$ maps to $\op{Aff}(1):=\mathbb{G}_m\ltimes \mathbb{G}_a$, and therefore acts on $\mathbb{A}^1_K$ by affine transformations.
      Hence $\widetilde{G}$ acts naturally on any coordinate affine space
              $\mathbb{A}^N_K=(\mathbb{A}^1_K)^N\,,\,\,\,N\geqslant 0$.

In the case of the motivic Galois group
  the algebra of functions $\mathcal{O}(\mathbb{A}^1_K)$ is isomorphic as $G$-module to ${\bf {H}}^\bullet(\mathbb{P}^\infty)=\oplus_{n\geqslant 0}K(-n)$.

\begin{dfn} For a given finite set $I$ and pair $(G,t^{1/2})$ as above, an $I$-colored  {\bf factorization system} is given by the following data:
\begin{itemize}
\item (equivariant sheaf) for a given $\gamma=(\gamma^i)_{i\in I}\in \mathbb{Z}^I_{\geqslant 0}- 0$ a $\widetilde{G}\times \op{Sym}_\gamma$-equivariant super coherent  sheaf $\mathcal{F}_\gamma$ on
 $$\mathbb{A}^\gamma_K:=\prod_{i\in I} \mathbb{A}^{\gamma^i}_K\,,$$
 where the group $\op{Sym}_\gamma=\prod_{i\in I} \op{Sym}_{\gamma^i}$ acts  by the permutation of coordinates
  in $\mathbb{A}^{\gamma^i}_K,\,\,i\in I$.
\item (multiplication) For $\forall \,\gamma_1,\gamma_2\in \mathbb{Z}^I_{\geqslant 0}- 0$ a multiplication
 isomorphism
 $$\mu_{\gamma_1,\gamma_2}:\left(\mathcal{F}_{\gamma_1}\boxtimes \mathcal{F}_{\gamma_2}\right)_{|U_{\gamma_1,\gamma_2}} \stackrel{\sim}{\to}{\mathcal{F}_{\gamma_1+\gamma_2}}_{|U_{\gamma_1,\gamma_2}}\,,$$
 where Zariski open affine subscheme
 $$U_{\gamma_1,\gamma_2}\subset \mathbb{A}^{\gamma_1+\gamma_2}_K
 \simeq \mathbb{A}^{\gamma_1}_K\times \mathbb{A}^{\gamma_2}_K$$
 is defined as
 $$\left\{\left( (x_{i,\alpha})_{i\in I,1\leqslant \alpha \leqslant \gamma_1^i},
 (y_{j,\beta})_{j\in I,1\leqslant \beta \leqslant \gamma_2^j} \right)\,|\,x_{i,\alpha}\ne y_{j,\beta}\,\forall\,\, i,\alpha,j,\beta\right\}\, .$$
 \end{itemize}
 These data are subject to the following conditions:
 \begin{itemize}
 \item (equivariance of multiplication) $\forall \,\gamma_1,\gamma_2\in \mathbb{Z}^I_{\geqslant 0}- 0$
 the isomorphism $\mu_{\gamma_1,\gamma_2}$ is equivariant with respect to
  the action of $\widetilde{G}\times \op{Sym}_{\gamma_1}\times \op{Sym}_{\gamma_2}$;
 \item (commutativity)  $\forall \,\gamma_1,\gamma_2\in \mathbb{Z}^I_{\geqslant 0}- 0$  we have
 $$\mu_{\gamma_2,\gamma_1}=\sigma_{\gamma_1,\gamma_2}\circ \mu_{\gamma_1,\gamma_2}\circ
 \sigma_{12},$$
 where $\sigma_{\gamma_1,\gamma_2}\in \op{Sym}_{\gamma_1+\gamma_2}$ is the collection of permutations for $i\in I$  acting on $\mathbb{A}_K^{\gamma_1^i+\gamma_2^i}$ by
 $$(x_{i,1},\dots,x_{i,\gamma_1^i},y_{i,1},\dots,y_{i,\gamma_2^i})\mapsto
  (y_{i,1},\dots,y_{i,\gamma_2^i},x_{i,1},\dots,x_{i,\gamma_1^i})$$
  and $\sigma_{12}$ is the standard permutation of indices;
  \item(associativity) $\forall \,\gamma_1,\gamma_2,\gamma_3\in \mathbb{Z}^I_{\geqslant 0}- 0$ two isomorphisms
  $$ \left(\mathcal{F}_{\gamma_1}\boxtimes \mathcal{F}_{\gamma_2}\boxtimes \mathcal{F}_{\gamma_3}\right)_{|U_{\gamma_1,\gamma_2,\gamma_3}} \stackrel{\sim}{\to}{\mathcal{F}_{\gamma_1+\gamma_2+\gamma_3}}_{|U_{\gamma_1,\gamma_2,\gamma_3}}  $$
   given by two ways to put brackets coincide with each other, i.e. we have in a simplified notation
   $$\mu_{\gamma_1+\gamma_2,\gamma_3}\circ \left(\mu_{\gamma_1,\gamma_2}\boxtimes  \op{id}_{\gamma_3}\right)=\mu_{\gamma_1,\gamma_2+\gamma_3}\circ
   \left(  \op{id}_{\gamma_1}\boxtimes\, \mu_{\gamma_2,\gamma_2}    \right).$$
 Here $U_{\gamma_1,\gamma_2,\gamma_3}\subset \mathbb{A}^{\gamma_1+\gamma_2+\gamma_3}_K$ is the intersection of pullbacks of  open sets $U_{\gamma_i,\gamma_j},\,\,1\leqslant i<j\leqslant 3$.\end{itemize}
\end{dfn}

We will also need a closely related notion.

\begin{dfn} For given $I$ and $(G,t^{1/2})$, a {\bf non-singular} $I$-colored factorization system is a factorization system
 together with a $\widetilde{G}\times \op{Sym}_{\gamma_1}\times\op{Sym}_{\gamma_2}$-equivariant extension of $\mu_{\gamma_1,\gamma_2}$
to a morphism of super coherent sheaves
$$\mathcal{F}_{\gamma_1}\boxtimes \mathcal{F}_{\gamma_2}\to \mathcal{F}_{\gamma_1+\gamma_2}$$
satisfying (in the notation of the previous definition) the constraints of commutativity and associativity on whole spaces $\mathbb{A}^{\gamma_1+\gamma_2}_K$ (resp.
 $\mathbb{A}^{\gamma_1+\gamma_2+\gamma_3}_K$) .
\end{dfn}

\subsection{Equivariant cohomology as a factorization system}

Let us return to the collection of torus equivariant cohomology for $\gamma\ne 0$:
$$\mathcal{H}'_\gamma={\bf H}^\bullet_{\mathsf{T}_\gamma}(\M_\gamma)\,\,$$
interpreted as global sections of super coherent sheaves $\mathcal{F}_\gamma$.
Here we omit the potential $W$ and  sector $V$  in order to simplify the notation.

Our goal in this subsection is to define a structure of non-singular factorization system on the collection of sheaves $\mathcal{F}_{\gamma}^{\prime}=\mathcal{F}'_\gamma\otimes \mathbb{T}^{\otimes \sum_{i,j}a_{ij}\gamma^i\gamma^j}$, with multiplication induced in a sense by the operation of direct sum of representations.

\begin{lmm} For any $\gamma_1,\gamma_2\in \mathbb{Z}^I_{\geqslant 0}-0$ there is a natural isomorphism
$${\bf H}^\bullet_{\mathsf{T}_{\gamma_1+\gamma_2}}(\M_{\gamma_1+\gamma_2},\M_{\gamma_1+\gamma_2}- \M_{\gamma_1}\times \M_{\gamma_2})\simeq$$
$$\simeq {\bf H}^\bullet_{\mathsf{T}_{\gamma_1}}(\M_{\gamma_1})\otimes {\bf H}^\bullet_{\mathsf{T}_{\gamma_2}}(\M_{\gamma_2})\otimes
\mathbb{T}^{\otimes \sum_{i,j} a_{ij}(\gamma^i_1\gamma^j_2+\gamma^j_1\gamma^i_2)}\,.$$
\end{lmm}
{\it Proof.} This is just the Thom isomorphism. Notice  that $\mathsf{T}_{\gamma_1+\gamma_2}=\mathsf{T}_{\gamma_1}\times \mathsf{T}_{\gamma_2}$ and
$$ \sum_{i,j} a_{ij}(\gamma^i_1\gamma^j_2+\gamma^j_1\gamma^i_2)=\op{dim}\M_{\gamma_1+\gamma_2}-\op{dim} (\M_{\gamma_1}\times
\M_{\gamma_2})     \,\,.\,\,\blacksquare$$
The standard map from the long exact sequence
$$\widetilde{\mu}_{\gamma_1,\gamma_2}:{\bf H}^\bullet_{\mathsf{T}_{\gamma_1+\gamma_2}}(\M_{\gamma_1+\gamma_2},\M_{\gamma_1+\gamma_2}- \M_{\gamma_1}\times \M_{\gamma_2})
\to {\bf H}^\bullet_{\mathsf{T}_{\gamma_1+\gamma_2}}(\M_{\gamma_1+\gamma_2})$$
gives a morphism of sheaves
$$\mathcal{F}_{\gamma_1}\boxtimes \mathcal{F}_{\gamma_2}\otimes
\mathbb{T}^{\otimes \sum_{i,j} a_{ij}(\gamma^i_1\gamma^j_2+\gamma^j_1\gamma^i_2)}\to \mathcal{F}_{\gamma_1+\gamma_2}\,,$$
where $\mathbb{T}$ can be interpreted at the trivial line bundle on $\mathbb{A}^{\gamma_1+\gamma_2}_K$ with the equivariant
$\widetilde{G}$ structure given by the character $t$.

\begin{prp} For any $\gamma_1,\gamma_2\in \mathbb{Z}^I_{\geqslant 0}-0$ the support of ${\bf H}^\bullet(\op{B}\mathsf{T}_{\gamma_1+\gamma_2})$-module
$${\bf H}^\bullet_{\mathsf{T}_{\gamma_1+\gamma_2}}(\M_{\gamma_1+\gamma_2}- \M_{\gamma_1}\times \M_{\gamma_2})$$ belongs to the
closed subset  $\mathbb{A}^{\gamma_1+\gamma_2}_K - U_{\gamma_1,\gamma_2}$.
\end{prp}
{\it Proof.} We will give a proof in the case when $char(\kk)=0$. The general case is similar.
Using the comparison isomorphism we may assume that $\kk=\C$ and ${\bf H}^{\bullet}={\bf H}^{\bullet}_{Betti}$.
We will need the following lemma.

\begin{lmm}
Let $\T^c$ be a compact torus acting on a complex algebraic variety $X$. Then for any locally closed real semi-algebraic subset $Z\subset X$
 the $\Z$-graded dual space to the cohomology
 $$H^\bullet_{c,\T}(Z,\Q)$$
is a finitely-generated $H^\bullet(\B\T^c,\Q)$-module with the support belonging to the union of affine spaces $\op{Spec} H^\bullet(\B\T_\alpha^c,\Q)\subset \op{Spec} H^\bullet(\B\T^c,\Q)$, where $\T_\alpha^c\subset \T^c$ are stabilizers of points in $Z$.
\end{lmm}
{\it Proof of Lemma}. We can stratify $Z$ by locally closed subspaces according to the dimension of the stabilizer in $\T^c$.
 By spectral sequences, the question reduces to the case when the stabilizer of any point $z\in Z$ is a given subgroup $\T_\alpha^c\subset \T^c$.
 Let us choose a ``complementary'' connected subtorus $\widetilde{\T}_\alpha^c\subset \T$ such that the multiplication 
$\T_\alpha^c\times \widetilde{\T}^c_\alpha\to \T^c$
 is surjective with finite kernel. Torus $\widetilde{\T}^c_\alpha$ acts almost freely (with finite stabilizers) on $Z$.
 Hence 
$$H^\bullet_{c,\T^c}(Z,\Q)\simeq H^\bullet_{c,\T_\alpha^c}(Z/\widetilde{\T}^c_\alpha,\Q)\simeq H^\bullet_{c,\T_\alpha^c}(pt,\Q)\otimes H^\bullet_c(Z/\widetilde{\T}^c_\alpha
,\Q)\,.$$
Rational cohomology with compact support of any real semi-algebraic set  is a finite-dimensional space, hence we obtain
 by duality a finitely-generated module supported on $\op{Spec} H^\bullet(\B\T_\alpha^c,\Q)$.  $\blacksquare$

In order to prove the Proposition we will apply the previous Lemma to $X=\M_{\gamma_1+\gamma_2},\,\T=\mathsf{T}_{\gamma_1+\gamma_2}$ and $Z=\M_{\gamma_1+\gamma_2}
-\M_{\gamma_1}\times \M_{\gamma_2}$. In this case $Z$ is an open subset of $X$, hence it is smooth. By Poincar\'e duality
  $H^\bullet_{c,\T}(Z,\Q)$ is dual to $H^\bullet_\T(Z,\Q) $ (up to shift).

For any $\gamma\in \mathbb{Z}^I_{\geqslant 0}$ and
 any point in $\M_{\gamma}$ its stabilizer in $\mathsf{T}_{\gamma}=\mathbb{G}_m^{|\gamma|},\,\,|\gamma|:=\sum_i\gamma^i$,
 corresponds to an equivalence relation $\sim$ on the finite set
$$\{1,\dots, |\gamma|\}=\sqcup_{i\in I}\{1,\dots,\gamma^i\}\,,$$
 and consists
 of sequences of scalars $(\lambda_k\in K^\times)_{1\leqslant k\leqslant |\gamma|}$ such that $\lambda_k=\lambda_{k'}$ if $k\sim k'$.
 The equivalence relation is defined such as follows. The super coherent sheaf corresponding to the contribution of points with such a stabilizer is a super vector bundle
  on the affine subspace consisting of configurations of points $(x_{i,\alpha})_{i\in I,1\leqslant \alpha\leqslant \gamma^i}$ such that
   $$x_{i,\alpha}=x_{i',\alpha'}\mbox{ for } (i,\alpha)\sim (i',\alpha')\,\,.$$
  For $\gamma=\gamma_1+\gamma_2$ and any
point in $\M_{\gamma_1+\gamma_2}- \M_{\gamma_1}\times \M_{\gamma_2}$ its stabilizer corresponds to an equivalence relation
which in terms of points $(x_{i,\alpha})_{i\in I,1\leqslant \alpha \leqslant \gamma_1^i},
 (y_{j,\beta})_{j\in I,1\leqslant \beta \leqslant \gamma_2^j} $ means that
   there exist $i_1,\alpha_1,j_1,\beta_1$ such that
  $$x_{i_1,\alpha_1}=y_{j_1,\beta_1}\,\,.$$
This is exactly the condition which defines the complement to $U_{\gamma_1,\gamma_2}$. $\blacksquare$

From the long exact sequence and the previous lemma it follows that the restriction of $\widetilde{\mu}_{\gamma_1,\gamma_2}$
 to $U_{\gamma_1,\gamma_2}$ is an isomorphism.

\begin{thm} \label{thm:modified-fact-sys} The collection of super coherent sheaves $$\mathcal{F}'_\gamma:=\mathcal{F}_{\gamma}\otimes \mathbb{T}_0^{\otimes \sum_{i,j} a_{ij}\gamma^i\gamma^j}\,,$$
with the multiplication isomorphisms obtained from $\widetilde{\mu}_{\gamma_1,\gamma_2}$
 by tensoring with the identity automorphism of a tensor power of $\mathbb{T}_0$, is
 a non-singular factorization system.
\end{thm}
 {\it Proof.} Straightforward. $\blacksquare$

\subsection{From factorization systems to admissible series}\label{sec:fact-admissible}

Let $\mathcal{F}=(\mathcal{F}_\gamma)$ be a non-singular factorization system. Then
the multiplication morphisms endow the space
$$A^{\mathcal{F}}_+:=\oplus_{\gamma\in \Z^I_{\geqslant 0}-0}A^{\mathcal F}_{\gamma},\,\,\,
\,A^{\mathcal{F}}_{\gamma}:= \Gamma(\mathcal{F}_\gamma)^{\op{Sym}_\gamma}$$
with a super-commutative associative {\it non-unital} product.

\begin{thm} Under the above assumptions, there  exists a collection of finite-dimensional super representations $(E_\gamma)_{\gamma\in
\Z^I_{\geqslant 0}-0}$ of ${G}$ such that two $\Z^I_{\geqslant 0}$-graded super representations of $G$, the spaces $A^{\mathcal{F}}_+$ and
 $$\op{Sym}_+ \left( \oplus_\gamma E_\gamma  \otimes (\oplus_{n\geqslant 0} \mathbb{T}^{\otimes n})\right)$$
have the same classes in $K_0(G-mod)$ in each graded component with respect to the $\Z^I_{\geqslant 0}$-grading.
\end{thm}

\begin{cor} Let $G=G_0\times \mathbb{G}_{m,K}$, $B=K_0(G_0-mod)$. Then the series
$$\sum_{\gamma\in \Z^I_{\geqslant 0}}[A^{\mathcal{F}}_{\gamma}]x^{\gamma}$$
is admissible.

\end{cor}

The proof of the Theorem is quite involved, but in the special case when $\mathcal{F}$ comes from Cohomological Hall algebra associated with some pair $(R,W)$ one can make
 a much shorter and clear argument which we explain in the next section. Nevertheless, we give here the proof in complete generality having in mind
 applications to certain generalizations of the Cohomological Hall algebra (see Section 7.6).

Before giving the proof we collect some preparatory material  concerning {\it Harrison homology}.
Recall (see e.g. \cite{Barr} or in modern language \cite{GetzKapr}) that the homological Harrison complex of any non-unital supercommutative algebra
 $A$ is defined as the free Lie super coalgebra generated by $A[1]$, with the differential arising from the product on $A$.
Explicitly, we  have
$$C_{Harr}(A):=\oplus_{k\geqslant 1}\left(L({k})\otimes (A[1])^{\otimes {k}}\right)^{\op{Sym}_{k}},$$
where $L(k)$ is a representation of $\op{Sym}_{k}$ of dimension $(k-1)!$, the tensor product of the contragredient to the
 representation $Lie(k)$ of the operad of Lie  algebras, and of the sign representation, i.e.
$$L(k)=Lie(k)^\vee\otimes sign_k\,\,.$$
If $A=\op{Sym}_+(V)$ is a free non-unital supercommutative algebra, then Harrison homology of $A$ is isomorphic to $V[1]$.
In general, the  non-unital symmetric algebra generated by the Harrison homology $H_\bullet^{\op{Harr}}(A)$ (shifted by $[-1]$) is
quasi-isomorphic
to $A$, with certain explicit homotopy.
Also we will use the following statement.
\begin{lmm} Let $B,C$ be two non-unital supercommutative algebras, and
 $$A=B\oplus C\oplus (B\otimes C)$$
with the obvious product. Then inclusions $B\mono A$ and $C\mono A$ induce an isomorphism
$$H_\bullet^{\op{Harr}}(B)\oplus H_\bullet^{\op{Harr}}(C)\simeq H_\bullet^{\op{Harr}}(A)\,.$$
\end{lmm}
{\it Idea of the proof.} Using spectral sequence one can see that the statement of  Lemma follows from the special case when multiplication on $B$ and $C$
 vanishes. Hence we have a purely operadic identity. $\blacksquare$

{\it Proof of  Theorem 12.} 
In our case there is a $\Z_{\geqslant 0}^I-\{0\}$-grading of both: of the algebra $A=A^{\mathcal{F}}_+$ and of its Harrison complex.
 This grading ensures an appropriate finiteness, and hence we have the following identity in the completed $K_0$-group:
$$\sum_{\gamma} [A^{\mathcal{F}}_{+,\gamma}]\, x^\gamma=\op{Sym}_+\left(-\sum_{\gamma}[H_{\bullet,\gamma}^{\op{Harr}}(A^{\mathcal{F}}_+)]\,x^\gamma\right).$$

Component of degree $\gamma$ of the Harrison complex of $A^{\mathcal F}$ is the space of global sections of a complex of $\widetilde{G}$-equivariant super coherent
sheaves  on the quotient scheme
$$X_\gamma:=\mathbb{A}^\gamma_K/\op{Sym}_\gamma\,.$$

\begin{prp} For any $\gamma\ne 0$ super vector space $H_{\bullet,\gamma}^{\op{Harr}}(A^{\mathcal{F}}_+)$ is a super coherent sheaf on $X_\gamma$ supported on the main diagonal
$$X_\gamma^{\op{diag}}\simeq \mathbb{A}^1_K\subset X_\gamma\,.$$
\end{prp}
{\it Proof of  Proposition.}
Recall that our factorization systems are defined over the field $K$ of characteristic $0$. Hence, we can safely assume that $K=\C$. In what follows it will be convenient for us to use analytic topology\footnote{A purely algebraic proof exists, but it is  less transparent.} on $\mathbb{A}^1(\C)=\C$.

Let $U\subset \C$ be an open subset. We define a nuclear $\Z^I_{\geqslant 0}$-graded non-unital Fr\'echet super algebra $A^{\mathcal{F}}_U$ by taking spaces
 of analytic sections over powers of $U$. We define analytic Harrison complex of $A^{\mathcal{F}}_U$ by taking completed tensor products.
 Obviously  $H_{\bullet,\gamma,an}^{\op{Harr}}(A_U^{\mathcal{F}})$ coincides with the space of analytic sections of the algebraic super coherent sheaf corresponding to  $H_{\bullet,\gamma}^{\op{Harr}}(A^{\mathcal{F}})$, over the Stein subset $U_{\gamma}\subset X_{\gamma}^{an}(\C)$  consisting of configurations
 of points on $U$.

 Let $x\in X_{\gamma}$ be a point which does not belong to  the main diagonal. Then $x$ corresponds to a finite subset $S\subset \C$ which has at least {\it two} elements. Therefore, we can find two disjoint non-empty sets $U',U''\subset \C$ such that
$$x\notin  U'_\gamma,\,\,x\notin U''_\gamma,\,\,x\in U_\gamma\mbox{ for } U:=U'\cup U'' .     $$
Then we have
$$A^{\mathcal{F}}_{U}\simeq A^{\mathcal{F}}_{U'}\oplus A^{\mathcal{F}}_{U''}\oplus \left(A^{\mathcal{F}}_{U'}\widehat{\otimes}A^{\mathcal{F}}_{U''}\right).$$
Applying Lemma 8 we conclude
that the Harrison cohomology of $A^{\mathcal{F}}_{U}$ in degree $\gamma\in \Z^I$ corresponds to a super coherent sheaf supported on
$U'_\gamma\cup U''_\gamma$. In particular, it vanishes in a neighborhood of $x$. Therefore $H_{\bullet,\gamma}^{\op{Harr}}(A^{\mathcal{F}})$
 vanishes near $x$. We proved that $H_{\bullet,\gamma}^{\op{Harr}}(A^{\mathcal{F}})$ is supported on the main diagonal. $\blacksquare$

Any $\widetilde{G}$-equivariant super coherent sheaf on $X_\gamma$ supported on the main diagonal is invariant under the action of the shift group $\mathbb{G}_a\simeq X_\gamma^{\op{diag}}$. Hence it admits a canonical finite filtration with quotients isomorphic to sums of copies of $\mathcal{O}(X_\gamma^{\op{diag}})$. Moreover, as $G^{(2)}$-modules these quotients are
   of the form
$$E\otimes \mathcal{O}(X_\gamma^{\op{diag}})=E \otimes (\oplus_{n\geqslant 0} \mathbb{T}^{\otimes n})\,,$$
where $E$ is a finite-dimensional representation of $G$. Hence $H_{\bullet,\gamma}^{\op{Harr}}(A^{\mathcal{F}}_+)$ has the same class
 in the completed $K_0$-ring of $G^{(2)}$ as $E_\gamma\otimes (\oplus_{n\geqslant 0} \mathbb{T}_0)^{\otimes n}$ for some finite-dimensional ${G}$-module $E_\gamma$. 
This concludes the proof of Theorem 12. $\blacksquare$

\subsection{Admissibility in the geometric case}

Let us assume that we are given a smooth $I$-bigraded algebra $R/\kk$ with potential $W$ and a cohomology theory $\mathbf{H}^\bullet$.
 All the arguments below will  work also if we fix a central charge $Z$ and a  sector $V$. In order to shorten the formulas we will skip
 $V$ from the notation.

\begin{dfn} A non-trivial representation of $R$ in coordinate spaces  is called
 $T$-{\bf indecomposable} if it can not be decomposed into the direct sum of two non-trivial subrepresentations in coordinate subspaces.
\end{dfn}
The set of $T$-indecomposable representations of $R$ with dimension vector $\gamma$ is a Zariski open subset (hence smooth)
 $$\M^{T-ind}_\gamma\subset \M_\gamma\,,$$
which is invariant under the action of $\op{Sym}_{\gamma}\ltimes\T_{\gamma}$.

Let us consider the ind-Artin stack
$$\sqcup_{\gamma \in\Z^I_{\geqslant 0}} \M_\gamma/\left(\op{Sym}_{\gamma}\ltimes\T_{\gamma}\right).$$
It carries a natural stratification with smooth strata by the dimension of the stabilizer. Strata are labeled by sequences $(\gamma_1,\dots,\gamma_k),\,k\geqslant 1$ of non-zero vectors $\gamma_i\in \Z^I_{\geqslant 0}$ such that $\sum_{i=1}^k \gamma_i=\gamma$, up to the action of the permutation group $\op{Sym}_k$.
  Point $E\in \M_\gamma$ belongs to the stratum corresponding to $(\gamma_1,\dots,\gamma_k)$ if and only if the representation $E$ can be decomposed into a direct sum of $T$-indecomposable subrepresentations $E^j,\,\,j=1,\dots,k$ in coordinate subspaces of dimension vectors $\gamma_j$.

Then we can use spectral sequences and Thom isomorphisms
 as in Section 5.2, and obtain an identity between generating series. At this point  it is convenient to use  cohomology of stacks  with compact support.
We have a constructible equivalence of ind-Artin stacks
$$\sqcup_{\gamma \in\Z^I_{\geqslant 0}} \M_\gamma/\left(\op{Sym}_{\gamma}\ltimes\T_{\gamma}\right)\simeq \op{Sym}\left(\sqcup_{\gamma\in \Z^I_{\geqslant 0}-\{0\}}
\M^{T-ind}_\gamma/\left(\op{Sym}_{\gamma}\ltimes\T_{\gamma}\right)  \right),$$
where for any ind-Artin stack $X$ we define $\op{Sym}(X):=\sqcup_{k\geqslant 0} X^k/\op{Sym}_k$. The constructible equivalence implies the identity
$$\sum_{\gamma \in\Z^I_{\geqslant 0}} [\mathbf{H}^\bullet_{c,\mathsf{T}_\gamma}(\M_\gamma,-W_\gamma)^{\op{Sym}_\gamma} ]x^\gamma= \op{Sym}\left(\sum_{\gamma\in\Z^I_{\geqslant 0}-\{0\}}[\mathbf{H}^\bullet_{c,\mathsf{T}_\gamma}(\M_\gamma^{T-ind},-W_\gamma)^{\op{Sym}_\gamma} ]x^\gamma
\right). $$
Passing to the dual spaces, and making the change of line variables $x^\gamma\mapsto \mathbb{L}^{\op{dim}\mathsf{T}_\gamma}x^\gamma$, we obtain the following formula (see Section 5.3 for the definition of $[\mathcal{H}_\gamma]$):
$$\sum_{\gamma \in\Z^I_{\geqslant 0}}[\mathcal{H}_\gamma] \, \mathbb{L}^{-\op{dim}\M_\gamma} x^\gamma=  \op{Sym}\left(\sum_{\gamma\in\Z^I_{\geqslant 0}-\{0\}}[\mathbf{H}^\bullet_{\mathsf{T}_\gamma}(\M_\gamma^{T-ind},W_\gamma)^{\op{Sym}_\gamma} ]\,\mathbb{L}^{-\op{dim}\M_\gamma}x^\gamma
\right). $$
Notice that for any $\gamma\in\Z^I_{\geqslant 0}-\{0\}$  one has an exact sequence of groups
$$1\to \mathbb{G}_m^{diag}\to \T_{\gamma}\to \mathsf{P}\T_{\gamma}\to 1\,,$$
where $\mathbb{G}_m^{diag}$ is the naturally embedded diagonal subgroup of $\T_{\gamma}$, and $\mathsf{P}\T_{\gamma}$ is the quotient group. Subgroup
  $\mathbb{G}_m^{diag}$ acts trivially on $\M_\gamma^{T-ind} $, and the quotient group $\mathsf{P}\T_{\gamma}$ acts {\it freely} on $\M_{\gamma}^{T-ind}$. Therefore we have
$$\mathbf{H}^{\bullet}_{\T_{\gamma}}(\M_{\gamma}^{T-ind},W_{\gamma})^{\op{Sym}_{\gamma}}\simeq \mathbf{H}^{\bullet}(\M_{\gamma}^{T-ind}/\mathsf{P}\T_{\gamma},W_{\gamma})^{\op{Sym}_{\gamma}}\otimes \mathbf{H}^{\bullet}(B\mathbb{G}_m^{diag})\,,$$
where in RHS we take ordinary cohomology of the quotient variety endowed with a function.
Hence we have proved  the following result which is parallel to Theorem 10 from Section 6.2:

$$\boxed{\mbox {\it The series $\sum_{\gamma \in\Z^I_{\geqslant 0}}[\mathcal{H}_\gamma] \, \mathbb{L}^{-\op{dim}\M_\gamma} x^\gamma$ is admissible.}}$$

\subsection{One-dimensional factorization systems}\label{sec:one-dim-factor}

There is an operation of tensor product (over $\left(\mathcal{O}(\mathbb{A}^\gamma_K)\right)_
{\gamma\in \Z_{\geqslant 0}^I-\{0\}}$) on factorization systems, and moreover the product of two non-singular systems is again non-singular.
In this section we will define a class of invertible factorization systems such that underlying sheaves are super line bundles on configuration spaces.
 These systems form a subgroup of the Picard group of the tensor category of factorization systems.

Let $\mathcal{B}=(b_{ij})_{i,j\in I}$ be an integral symmetric matrix. Denote by $\epsilon:\Z^I\to \Z/2\Z$ the group homomorphism given by
$$\epsilon(\gamma):=\sum_{i\in I} b_{ii} \gamma^i \op{mod} \,2 \,\,.$$
By the arguments from Section 2.6 there exists a homomorphism (bilinear form)
$$\phi:\Z^I\otimes \Z^I\to \Z/2\Z$$ such that
$$\sum_{i,j} b_{ij} \gamma_1^i\gamma_2^j +\epsilon(\gamma_1)\epsilon(\gamma_2)+\phi(\gamma_1,\gamma_2)+\phi(\gamma_2,\gamma_1)=0\,\op{mod}\,2\,\,\,\,\,\,\forall \,\gamma_1,\gamma_2\in \Z^I\,\,.$$
Let us choose such a bilinear form $\phi$.
 We define a collection of super coherent sheaves $\mathcal{L}_\gamma^\mathcal{B}$ as trivialized line bundles on $\mathbb{A}_K^\gamma$
 (with generators denoted by $s_\gamma$) endowed with the parity $\epsilon(\gamma)$. We define the action of $\op{Sym}_\gamma$ to be the standard action on
$\mathcal{L}_\gamma^\mathcal{B}=\mathcal{O}_{\mathbb{A}_K^\gamma}$ via permutation of coordinates. The  group $\widetilde{G}$ acts on the
 canonical generator of $\mathcal{L}_\gamma^\mathcal{B}$ via the representation $\left(K(-1/2)\right)^{\otimes \sum_{i,j} b_{ij}
\gamma^i\gamma^j}$, where $K(-1/2)$ is the square root of the Tate motive.

Next,  the multiplication isomorphism on $U_{\gamma_1,\gamma_2}$ is defined on the standard generators such as follows:
 $$\mu_{\gamma_1,\gamma_2}: s_{\gamma_1}\otimes s_{\gamma_2}\mapsto
\phi(\gamma_1,\gamma_2)\prod_{i,\alpha,j,\beta}(x_{i,\alpha}-y_{j,\beta})^{b_{ij}}  s_{\gamma_1+\gamma_2}\,\,.$$

\begin{lmm} The collection $(\mathcal{L}_\gamma^\mathcal{B})$ of equivariant super coherent sheaves together with isomorphisms $\mu_{\gamma_1,\gamma_2}$ defined above, forms a factorization system. It is a non-singular factorization system if $b_{ij}\geqslant 0$ for all $i,j\in I$.
\end{lmm}
{\it Proof.} The associativity is obvious, the commutativity follows by a straightforward check from the relation between
 $(b_{ij}),\epsilon,\phi$. $\blacksquare$

Similar to Section 2.6, one can show that $(\mathcal{L}_\gamma^\mathcal{B})$ does not depend (up to an isomorphism) on the choice
 of the bilinear form $\phi$, and moreover the identification can be made canonical.

In the non-singular case when $b_{ij}\geqslant 0$ for all $i,j\in I$, this factorization system admits a very explicit geometric interpretation. Namely, let us consider a quiver $Q^\mathcal{B}$ with the set of vertices $I$ and $b_{ij}$ arrows between vertices $i$ and $j$.
 We endow this quiver with the non-degenerate quadratic potential
$$W^{\mathcal{B}}						=\sum_{i\ne j}\sum_{1\leqslant l\leqslant b_{ij}}x_{i,j;l} x_{i,j;l}^*+\sum_{i\in I}\sum_{1\leqslant l\leqslant b_{ii}}(y_{i;l})^2\,.$$
Here $\{x_{i,j;l}\}_{1\leqslant l\leqslant b_{ij}}$ is the set of all arrows in $Q_\mathcal{B}$ connecting vertices $i$ and $j\ne i$, and $*$ denotes an involution
 on arrows reversing orientation. Also $\{y_{i;l}\}_{1\leqslant l\leqslant b_{ii}}$ is the set of loops at the vertex $i$.

Pair $(Q^\mathcal{B},W^\mathcal{B})$  gives  a non-singular factorization system
$$\mathcal{F}^{\mathcal{B}}_\gamma=\mathbf{H}^\bullet_{\mathsf{T}_\gamma}(\M_\gamma^{Q^\mathcal{B}}, W^\mathcal{B}_\gamma)\otimes \mathbb{T}_0^{-\otimes \sum_{i j} b_{ij}\gamma^i\gamma^j}\, .$$

Function $W^\mathcal{B}_\gamma$ is a non-degenerate quadratic form on the affine space $\M_\gamma^{Q^\mathcal{B}}$. By purity, we have
$$\mathbf{H}^\bullet_{\mathsf{T}_\gamma}(\M_\gamma^{Q^\mathcal{B}}, W^\mathcal{B}_\gamma)\simeq \mathbf{H}^\bullet(\mathcal{B}\mathsf{T}_\gamma)\otimes
\mathbf{H}^\bullet(\M_\gamma^{Q^\mathcal{B}}, W^\mathcal{B}_\gamma)\,.$$
 Let us assume that the cohomology theory $\mathbf{H}^\bullet$ does
 not distinguish quadratic forms of the same ranks (like e.g. EMHS). In this case we can write
$$\mathbf{H}^\bullet(\M_\gamma^{Q^\mathcal{B}}, W^\mathcal{B}_\gamma)\simeq \left(\mathbb{T}^{\otimes 1/2}\right)^{\otimes \sum_{i j} b_{ij}\gamma^i\gamma^j},$$
where $\mathbb{T}^{\otimes 1/2}=\mathbf{H}^\bullet(\mathbb{A}^1_\kk,-z^2)=\mathbb{T}_0^{\otimes 1/2}[-1]$.
The action of the permutation group $\op{Sym}_{\gamma}$ on  one-dimensional space $\mathbf{H}^\bullet(\M_\gamma^{Q^\mathcal{B}}, W^\mathcal{B}_\gamma)$ is trivial.

A straightforward check shows that factorizations systems
 $(\mathcal{F}^{\mathcal{B}}_\gamma)$ and $(\mathcal{L}_\gamma^\mathcal{B})$ are isomorphic.

\subsection{Approximation by large integers and end of the proof of Theorems 9,10.
}\label{sect:endofproof}

 The statement of Theorem 9 is a certain universal identity in $\lambda$-rings, where
the  coefficients of polynomials $f_\gamma$ in the notation of Definition 19 are free $\lambda$-variables.
First, we will prove it in the case when $b_{ij}\geqslant 0\,\,\,\,\forall i,j\in I$. This can be done  in a ``model-theory fashion", i.e. by using certain universal identities in $\lambda$-rings.

We start with a non-singular $I$-colored factorization system, where
$G=G_\alpha\times \mathbb{G}_{m,K}$ for some algebraic group $G_{\alpha}$, and $t^{1/2}: G\to \mathbb{G}_{m,K}$ is the projection to the second factor (see Definition 22). Let $B_{\alpha}$ denotes the $\lambda$-ring $K_0(G_{\alpha}-mod)$ and $y_{\gamma}\in B_{\alpha}[q^{\pm 1/2}]$ be a collection of elements parametrized by $\gamma\in \mathbb{Z}_{\geqslant 0}^I-\{0\}$. Then there exists a non-singular factorizations system and the corresponding collection $(E_{\gamma})_{\gamma\in \mathbb{Z}_{\geqslant 0}^I-\{0\}}$ of super representations (see Theorem 12) such that $y_{\gamma}=[E_{\gamma}]\in K_0(G_{\alpha}-mod)$. This can be shown by induction by the norm $|\gamma|$ (for any choice of the norm). Namely, we can consider free non-singular system generated by arbitrary $\widetilde{G}\times \op{Sym}_\gamma$-equivariant super coherent sheaves on $\mathbb{A}^\gamma_K$ supported on the main diagonal for $\gamma\ne 0$.

Next, we can twist the factorization system by the one-dimensional system $(\mathcal{F}^{\mathcal{B}}_\gamma)$ (see Section 6.8) and obtain another non-singular factorization system by Theorem 11. Then Theorem 12 gives a new collection
$(E_{\gamma}^\prime)_{\gamma\in \mathbb{Z}_{\geqslant 0}^I-\{0\}}$ and the corresponding classes
$y_{\gamma}^\prime\in B_{\alpha}[q^{\pm 1/2}]$. By Corollary to Theorem 12 we obtain the admissibility property.

In order to prove Theorem 9 for general $\lambda$-rings it suffices to have it for  infinitely generated free $\lambda$-rings. Those can be ``approximated" by $\lambda$-rings of the type $K_0(G_{\alpha}-mod)$ such as follows.
Let us fix a non-zero element $\gamma_0$ and number $c\geqslant 0$ such that
$$y_{\gamma}=\sum_{|i|\leqslant c}y_{\gamma,i}q^{i/2}
$$
for all $|\gamma|\leqslant |\gamma_0|$.
Then $y_{\gamma_0}$ is a universal $\lambda$-polynomial in variables  $y_{\gamma,i},|\gamma|\leqslant |\gamma_0|, |i|\leqslant c $ and $\frac{1}{1-q},\,q^{\pm 1/2}$. A priori it belongs to $B_{\alpha}((q^{1/2}))$. The fact that $y_{\gamma_0}$ belongs to $B_{\alpha}[q^{\pm 1/2}]$ is  equivalent to a (finite) collection of equations on $y_{\gamma,i}$. We would like to treat $y_{\gamma,i}$ as free $\lambda$-variables. Although this is not the case, we can find a collection of algebraic groups $(G_{\alpha})$ and elements $y_{\gamma,i}^{({\alpha})}\in K_0(G_{\alpha}-mod)$ such that the kernels of the homomorphisms of the free $\lambda$-ring generated by the ``universal" symbols $y_{\gamma,i}^{univ}$ to $K_0(G_{\alpha}-mod), y_{\gamma,i}^{univ}\mapsto y_{\gamma,i}^{({\alpha})}$ are the ideals with the trivial intersection over the set of indices $\alpha$. For example one can take as $G_{\alpha}$ the product of finitely many copies of $GL(N)$ for large $N$. This finishes the proof of Theorem 9 for $b_{ij}\geqslant 0$.

Let $\mathcal{B}=(b_{ij})$ be a matrix with integer but not necessarily non-negative coefficients.
For each $n\in \mathbb{Z}$ we consider the matrix $\mathcal{B}^{(n)}$ with entries
$$b_{ij}^{(n)}:=b_{ij}+n\,,\,\,\, n\in \Z\,.$$
For sufficiently large $n$ all entries are non-negative, and hence we have admissibility by the above considerations.
For a given $\gamma\ne 0$ let us consider the dependence of  $[E_\gamma]$ on $n$. This $\lambda$-polynomial satisfies  the following Proposition.

\begin{prp}
Let $B$ be a $\lambda$-ring, and $P$ be a $\lambda$-polynomial with coefficients in $B$ in three  variables $x,y,z$,
 where $y,z$ are invertible line elements. Suppose that for any sufficiently large integer $n\gg 0$ the value
 $$P\left(\frac{1}{1-q},q^{1/2},(-q^{1/2})^n\right)\in B((q^{1/2}))$$
is a Laurent polynomial in the line variable $q^{1/2}$ . Then the same is true for any integer $n\in \Z$.
\end{prp}

{\it Proof.} First, we have an identity in the $\lambda$-ring $\Z((q))$:
$$\lambda^k\left(\frac{1}{1-q}\right)=\frac{q^{k(k-1)/2}}{(1-q)\dots (1-q^k)}\,,\,\,\,\forall k\geqslant 0\,.$$
Hence the evaluation $P\left(\frac{1}{1-q},q^{1/2},(-q^{1/2})^n\right)$ can be identified with
$$\frac{1}{(1-q^A)^D}\cdot \sum_{i,j:\,|i|+|j|<N} c_{ij} (-q^{1/2})^{i+nj}$$
for some integers $A,D,N\geqslant 1$ and elements $c_{ij}\in B$.
Our assumption means that for any $n\geqslant 0$ the above series is a Laurent polynomial, i.e. its coefficients in sufficiently
large powers of $-q^{1/2}$ vanish.$\Z^I_{\geqslant 0}$-grade

Let us fix two residues $a_1,a_2\in \Z/A\Z$. It is easy to see that the coefficient of $(-q^{1/2})^k$ for $k\gg 1$ and $k=a_1\pmod{A},n=a_2\pmod{A}$
 can be written as
 $$\sum_{i,j} c_{ij}\, F_{i,j,a_1,a_2}(k,n)\,,$$
where $F_{i,j,a_1,a_2}$ is a polynomial in two variables depending on residues $(a_1,a_2)$ modulo $A$, with integer values at $(k,n)\in \Z^2$ .
 Consider the additive subgroup $\Gamma$ of $B$ generated by $(c_{ij})$. It is a finitely generated abelian group.
 Pick any additive functional $\Gamma\to \Z/M\Z$ for some $M\geqslant 0$. Then we use the following obvious statement.
\begin{lmm} Let $f:\Z\times \Z \to \Z$ be a $\Z$-valued polynomial. If $f(x,y)=0\pmod{M}$ for all $x\in \Z$ and all $y\gg 1$ then
 $f(x,y)=0\pmod{M}$ for all $(x,y)\in \Z^2$.
\end{lmm}
Applying the above Lemma we obtain the proof of the Proposition.
$\blacksquare$

Returning to the proof of Theorem 9 we apply the above Proposition to the $\lambda$-polynomial $[E_{\gamma}]$ as a function of $n$ and conclude that it belongs to $B[q^{\pm 1/2}]$ for $n=0$. This finishes
the proof of Theorem 9.

 Theorem 10 now follows from Theorems 9, 11 and 12. $\blacksquare$

\section{Critical COHA}

\subsection{Description of results}

Let $R$ be  a smooth $I$-bigraded algebra over  a field  $\kk$, endowed with a bilinear form $\chi_R$ on $\Z^I$ compatible with the Euler characteristic,
  and a potential $W\in R/[R,R]$, as in Section 4.7.
 Also, suppose that we are given an additional data, consisting of a collection of $\G_\gamma$-invariant  closed subsets $\M_\gamma^\s\subset \M_\gamma$ for all
 $\gamma\in \Z^I_{\geqslant 0}$ (superscript $^\s$ means ``special'')
 satisfying the following conditions:
\begin{itemize}
\item for any $\gamma$ we have $\M_\gamma^\s\subset \op{Crit}(W_\gamma)$, i.e. 1-form $dW_\gamma$ vanishes at $\M_\gamma^\s$,
\item for any short exact sequence $0\to E_1\to E\to E_2\to 0$ of representations of $\overline{\kk}\otimes_\kk R$ with dimension vectors $\gamma_1,\gamma:=\gamma_1+\gamma_2,\gamma_2$ correspondingly, such that all $E_1,E_2,E$ are critical points of the potential,the representation $E$ belongs to $\M_\gamma^\s$ if and only if
 both representations $E_1,E_2$ belong to $\M_{\gamma_1}^\s,\M_{\gamma_2}^\s$ respectively.
\end{itemize}
The last condition implies that the collection of representations from $\M_\gamma^\s(\overline{\kk})$ for all $\gamma\in \Z^I_{\geqslant 0}$ form an abelian category, which is  a Serre subcategory of the abelian category $\op{Crit}(W)(\overline{\kk}):=\sqcup_\gamma \op{Crit}(W_\gamma)(\overline{\kk})$, which is itself a full subcategory of $\overline{\kk}\otimes_\kk R-\op{mod}$.
For example, one can always make the maximal choice
    $$\M_\gamma^\s:= \op{Crit}(W_\gamma)\,\,\,\,\forall \gamma\in \Z^I_{\geqslant 0}\,. $$
We can construct more examples such as follows. Pick an {\it arbitrary} subset  $N\subset R$ and  define
 $\M_\gamma^\s$ as the set of representations belonging to $\op{Crit}(W_\gamma)$ for which all elements $n\in N$ act as {\it nilpotent} operators.

Assume that $\kk=\C$. We will define the {\it critical COHA} as
$$\mathcal{H}=\oplus_{\gamma\in \Z^I_{\geqslant 0}}\H_\gamma\,,$$ 
where
$$\H_\gamma:= \bigoplus_{z\in \C}\left( H^{\bullet}_{\G_\gamma,c} (\M_\gamma^\s\cap W_\gamma^{(-1)}(z), \phi_{W_\gamma-z} \Q_{\M_\gamma})
\right)^\vee \otimes \mathbb{T}^{\otimes \dim M_\gamma/G_\gamma}\,.$$
Here we use equivariant cohomology with compact support with coefficients in the sheaf of vanishing cycles (see below).

We will define   a structure of EMHS of  special type (called {\it monodromic} mixed Hodge structure)  with Betti realization $\H_\gamma$, and a twisted associative product. For the critical COHA in Tannakian category EMHS there are analogs of results from Sections 5 and 6. In particular,  for any stability condition we obtain a factorization of the motivic DT-series, and all motivic DT-series associated with  sectors are admissible.
  The proof of the factorization property for a chosen stability condition is much more complicated than in the case of the exponential cohomology.
 It is based on certain ``integral identity'' involving vanishing cycles and conjectured in our earlier work \cite{KS} (also we proposed there a sketch of the proof). The reason for the complication is the lack of Thom
isomorphism for vanishing cycles. In Section 7.8 we will give a complete proof of the identity.

Under certain assumptions one can show that the critical COHA has the same motivic  DT-series as the one defined for rapid decay cohomology. In particular, we can use this fact in the case of quiver $Q_1$ with polynomial potential. Also, ``critical'' DT-series introduced below matches
 those introduced in \cite{KS} in the framework of ind-constructible 3-dimensional Calabi-Yau categories. We should mention that the idea 
of defining Donaldson-Thomas invariants using compactly supported cohomology of the vanishing cycle complex appeared first in \cite{DimcaSzendroi} and also was implicit in \cite{KS}.

The generality of our construction is still not satisfactory for all applications.  One needs a minor generalization of the theory of vanishing cycles to formal schemes and stacks, which is still
 absent in the existing literature. In the ideal picture, our theory will be applicable to a large class of 3-dimensional Calabi-Yau categories,
 including those associated with quivers with formal potentials, e.g. the cluster ones. 

\subsection{Vanishing cycles in the analytic case}
 Let us first recall basic facts about vanishing cycles for constructible sheaves in analytic geometry (see e.g. \cite{Schuermann}).
 If $X$ is a complex analytic space endowed with a holomorphic function $f:X\to \C$, then we have the functor $\psi_f$ of nearby cycles
from $D^b_c(X)$, the derived category of complexes of sheaves on $X$ constructible with respect to a complex analytic stratification, into a similar one for $X_0:=f^{-1}(0)$, given (in the obvious notation) by
$$\psi_f:=(X_0\to X)^*\circ (X_+\to X)_*\circ (X_+\to X)^*\,,\,\,\,\,X_+:=f^{-1}(\R_{>0})\subset X\,.$$
Also we have the  functor $\phi_f$ of vanishing cycles   given by the cone of the adjunction
$$\phi_f\mathcal{F}\simeq\op{Cone}\left( (X_0\to X)^* \mathcal{F}\to \psi_f \mathcal{F}\right).$$

 The functor of vanishing cycles commutes  with the direct image for proper morphisms $\pi:X_1\to X_2$ for varieties with functions $f_1,f_2$ such that $f_1=\pi^* f_2$,
$$\phi_{f_2}\circ \pi_*\simeq \pi_*\circ \phi_{f_1}$$

Thom-Sebastiani theorem for sheaves (see e.g. \cite{Massey}, \cite{Schuermann}) says that that the shifted functor of vanishing cycles commutes with the external products. In other words, if $X=X^{(1)}\times X^{(2)}$ and $f=f^{(1)}\boxplus f^{(2)}$ for $f^{(i)}\in \mathcal{O}(X^{(i)}),\,i=1,2$ and
$$\mathcal{F}=\mathcal{F}^{(1)}\boxtimes \mathcal{F}^{(2)}:=\op{pr}_{X\to X^{(1)}}^*(\mathcal{F}^{(1)})\otimes \op{pr}_{X\to X^{(2)}}^*(\mathcal{F}^{(2)})$$
 for $\mathcal{F}^{(i)}\in D^b_c(X^{(i)})\,,\,i=1,2$, then
$$\left( X^{(1)}_0\times X^{(2)}_0\to X_0 \right)^* \phi_f [-1] \mathcal{F}\simeq \phi_{f^{(1)}}[-1] \mathcal{F}^{(1)}\boxtimes \phi_{f^{(2)}}[-1] \mathcal{F}^{(2)}\,,$$
 where $X^{(1)}_0,X^{(2)}_0,X_0$ are zero loci of functions $f^{(1)},f^{(2)}, f$ respectively.

\subsection{Relation with the rapid decay cohomology}

For a large class of functions one can relate rapid decay cohomology and vanishing cycles.
Let $X$ be a complex {\it algebraic} variety. For a regular function $f\in \mathcal {O}(X)$ considered as a map $X\to \mathbb{A}^1_\C$ we define its {\bf bifurcation set} $\op{Bif}(f)\subset \C$ as the set of points over which $f$ is not a locally trivial fibration
 (in the analytic topology). The set $\op{Bif}(f)$ is  finite, containing singularities of the direct image $f_*\Z_{X(\C)}$. Also, 
for smooth $X$ the set  $\op{Bif}(f)$ contains the set of critical values of $f$.

Let us choose an isotopy class of a collection of disjoint paths from $-\infty$ to points of
$\op{Bif}(f)$. Then we obtain a canonical isomorphism
$$H^\bullet(X,f)=H^\bullet(X, f^{-1}(-c))\simeq \bigoplus_{z\in \op{Bif}(f)} H^\bullet\left(f^{-1}(B_\delta(z)),f^{-1}(z-\delta)\right)$$
for sufficiently small real $0<\delta\ll 1$ and sufficiently large real $c\gg 0$. Here $B_\delta(z):=\{w\in \C\,|\,\,\,|z-w|\leqslant \delta\}$. Indeed, the pointed topological space
$X(\C)/f^{-1}(-c)$ is homotopy equivalent to the wedge sum (coproduct) of pointed spaces
 $$\bigvee_{z\in \op{Bif}(f)} f^{-1}(B_\delta(z))/f^{-1}(z-\delta)\,,$$ with the homotopy equivalence depending on  the collection of disjoint paths (up to isotopy).

\begin{dfn} A function $f$ on {\bf smooth}  algebraic variety $X/\C$ is called  {\bf topologically isotrivial at infinity} if there exists
a  $C^\infty$-manifold with boundary $U\subset X(\C),\,\,\dim_\R U=2\dim X$ such that $f_{|U}:U\to \C$ is proper, and a homeomorphism
 $$h:\partial U \times [0,+\infty)\simeq X(\C)- \op{int} U\, ,\,\,\,h(x,0)=x\mbox{ for }x\in \partial U\,,$$
such that $f(h(x,t))=f(x)$ for any $x\in \partial U$.
\end{dfn}

In particular, if $f:X\to \C$
is a proper map for smooth $X$, then it is topologically isotrivial at infinity.
We were not able to find a discussion of this property in the literature. The closest analog is a strictly stronger property
  called {\it M-tame} (see \cite{Sabbah2}).
   It is easy to see that for a function $f$ topologically isotrivial at infinity, the bifurcation set coincides with the set of
critical values, and
for any $z\in \op{Bif}(f)$ one has a canonical isomorphism
$$H^\bullet(f^{-1}(z),\phi_{-f+z}\, {{\mathbb{Z}}}_X[-1])\simeq H^\bullet\left(f^{-1}(B_\delta(z)),f^{-1}(z-\delta)\right)$$
for $0<\delta\ll 1$. Hence we see that for such $f$ the direct sum
over critical values
$$\bigoplus_{z\in \C} H^\bullet(f^{-1}(z),\phi_{-f+z}\, {{\mathbb{Z}}}_X[-1])$$
is isomorphic
 to $H^\bullet(X,f)$. The isomorphism is  not canonical, it depends on a choice of an isotopy class of a collection of disjoint paths from $-\infty$ to points of
$\op{Bif}(f)$, as before.

Also we have a (non-canonical) isomorphism between cohomology with compact support for a topologically isotrivial at infinity function $f$
$$H^\bullet_c(X,f)\simeq \bigoplus_{z\in \C} H^\bullet_c(f^{-1}(z),\phi_{-f+z}\, {{\mathbb{Z}}}_X[-1])\,.$$
This follows by Verdier duality from the previous considerations.

Here is an application promised in Section 4.7. Let $W\in \C[x]$ be a polynomial of degree $N\geqslant 3$ such that $W'$ has $N-1$ distinct roots $r_1,\dots,r_{N-1}\in \C$.
We can consider $W$ as a cyclic polynomial in $\C Q_1$. Hence for any $n\geqslant 0$ we get function $W_n$ on $\op{Mat}(n\times n,\C)$ given
 by $W_n(X)=\op{Tr} W(X)$. One can show that for an appropriate $C^\infty$ function $F:\C\to \R_{>0}$ and any $n$, the submanifold $U_n\subset\op{Mat}(n\times n,\C)$ given by the inequality
 $$\op{Tr}(XX^\dagger)\leqslant F(\op{Tr}(W(X)))$$
satisfies an $U(n)$-equivariant version of the conditions from the definition
 of topological isotriviality at infinity. Thus implies that one can identify equivariant rapid decay cohomology
 with the equivariant cohomology of vanishing cycles. The latter is easy to compute.

The set of critical points of the  potential $W_n$
 consists of operators $x\in \op{Mat}(n\times n,\C)$ satisfying the equation
$$W'(x)=0\,\,.$$
This means that $x$ has eigenvalues $r_1,\dots ,r_{N-1}$ with some multiplicities, and no non-trivial Jordan blocks.
 Therefore the space $\C^n$  splits into the direct sum $\oplus_{i=1}^{N-1}V_i$ of eigenspaces.
 We conclude that the set of critical points of $W_n$ is the disjoint union over all ordered partitions
 $$n_1+\dots+n_{N-1},\,\,n_i\geqslant 0$$
of components corresponding to such decompositions with $\op{dim}V_i=n_i$.

Potential $W_n$ has Bott-Morse singularity at any
critical point, i.e. locally (in analytic topology) it is the Thom-Sebastiani  sum $\boxplus$ of a constant function and of a quadratic form.

The rank of the second derivative $W_n''$ at the component corresponding to the  decomposition $n_1+\dots+n_{N-1}$, is
 equal to $\sum_i n_i^2$. Also, this component is a homogeneous space $GL(n,\C)/\prod_i GL(n_i,\C)$. One can check that the sheaf
of vanishing cycles (shifted by $[-1]$) in our case is the constant sheaf $\underline{\Z}$ shifted by $[-\sum_i n_i^2]$.
Therefore, we conclude that
$$\H_n\simeq \bigoplus_{\substack{ n_1,\dots,n_{N-1}\geqslant 0
\\n_1+\dots+n_{N-1}}} \prod_{i=1}^{N-1} H^\bullet(\op{BGL}(n_i,\C))\,\,.$$
Then $\H=\bigoplus_{n\geqslant 0}\H_n$ is isomorphic as a bigraded space to
 $$\left(\bigoplus_{m\geqslant 0}H^\bullet(\op{BGL}(m,\C))\right)^{\otimes (N-1)}\,\,.$$
It looks plausible  that $\H$ considered as an algebra (and not only as a  bigraded vector space),
 is also isomorphic to the $(N-1)$-st tensor power of the  exterior algebra corresponding to the case $N=2$.
 Also it looks plausible that the condition on $W$ to have simple critical points  can be dropped.

\subsection{Vanishing cycles and monodromic Hodge structures}

The definition of vanishing cycles given in the analytic case is not algebro-geometric. To our knowledge there is no satisfactory general geometric approach as well as a proof of Thom-Sebastiani theorem.
 There are individual theories in the \'etale case (see \cite{SGA}), and for mixed Hodge modules (see \cite{Saito}). In what follows we
 will  work over $\kk=\C$ and use mixed Hodge modules, although all arguments hold also for \'etale constructible sheaves as well.

Cohomology groups of sheaves
 of vanishing cycles carry a natural mixed Hodge structure. In fact,  they should be considered as EMHS of special type, which we will call monodromic.
 The origin of the name is explained by the fact that  there is a natural action of monodromy  operator on $\phi_f \mathcal{F}$.
 The monodromy comes from the parallel transport around $0$ for locally constant family of sheaves $\phi_{\lambda f} \mathcal{F}$ depending on a parameter $\lambda\in \C^*$.

\begin{dfn} Tannakian category of {\bf monodromic mixed Hodge structures} $MMHS$ is the full Tannakian subcategory of EMHS consisting
 of objects unramified on $\C^*\subset \C$.
\end{dfn}

As an abelian category $MMHS$ is naturally equivalent to the category of unramified mixed Hodge modules on $\mathbb{G}_{m,\C}$,  i.e.  admissible variations of mixed Hodge structures on $\C^*$. The correspondence (shifted by 1) is given in one way by the restriction  from $\C=\mathbb{A}^1(\C)$ to $\C^*$, and in another way by the direct image with compact supports. 

 Tensor product  on Tannakian category $MMHS$ {\it does not} coincide with the tensor product of corresponding variations of mixed Hodge structures on on $\C^*$, but for the underlying constructible sheaves of $\Q$-vector spaces there is an isomorphism of two tensor products.

\begin{dfn} Let $X/\C$ be an algebraic variety endowed with a function $f\in \mathcal{O}(X)$, and a locally closed subset $X^{\s}\subset X_0:=f^{-1}(0)$. The {\bf critical cohomology with compact support} is defined as the cohomology of an object of $D^b(MMHS)\subset D^b(MHM_{\mathbb{A}^1_\C})$ , given by
\[(\mathbb{G}_m\to \mathbb{A}^1_\C)_! \,(X^\s\times \mathbb{G}_{m}\to \mathbb{G}_m)_!\,\left((X^\s\times \mathbb{G}_m)\to (X_0\times \mathbb{G}_{m})\right)^* \,\phi_{\frac{f}{u}}\Q_{X\times\mathbb{G}_m},\]
where affine line $\mathbb{A}^1_\C$  is endowed with coordinate $u$.
We will denote it by ${H}^{\bullet,crit}_{c,MMHS}(X^{\s},f)$. 
\end{dfn}

The above Definition becomes very transparent when one applies the Betti realization functor at $u=1$. Then the critical cohomology with compact support turns into a space
$$H_c^{\bullet}(X^\s,(X^\s\to X_0)^{\ast}\phi_f(\Q_X))$$
endowed with the action of the monodromy operator.

The reader should notice that (as we already pointed out in Section 4.5) the critical cohomology depends on the ambient space $X$ and the corresponding function $f_X$. This is quite different from the case of rapid decay cohomology.

There is a version of the Thom-Sebastiani theorem  in the setting of mixed Hodge modules
 (M.~Saito, \cite{Saito-TS}). It implies in particular an isomorphism
$${H}^{\bullet,crit}_{c,MMHS}(X^{(1),\s},f^{(1)})\otimes {H}^{\bullet,crit}_{c,MMHS}(X^{(2),\s},f^{(2)})\simeq {H}^{\bullet,crit}_{c,MMHS}(X^{(1),\s}\times X^{(2),\s},f)$$
where
$X=X^{(1)}\times X^{(2)},\,f=f^{(1)}\boxplus f^{(2)},\,f^{(i)}\in \mathcal{O}(X^{(i)}), \, i=1,2$ and the tensor product is the additive convolution $*_{+}$.

Let $\kk$ be a field, $\op{char}\kk=0$ and ${\bf H}^\bullet $ be one of the standard  cohomology theories, namely \'etale cohomology, Betti cohomology,
 or mixed Hodge modules. Then we  can define similarly to what has been done above,  monodromic exponential mixed motives  ${\bf H}^{\bullet,crit}_{c}(X^{\s},f)$
 for  triples $(X,f,X^\s)$ defined over $\kk$. It is expected that the Thom-Sebastiani theorem holds also in the \'etale case (P.~Deligne, private communication)\footnote{After  this pape was finished, a proof of the Thom-Sebastiani theorem in the case of positive characteristic was announced by Lei Fu, see \cite{Fu}.}.
 In what follows we will use the notation ${\bf H}^{\bullet,crit}_{c}$ assuming the Thom-Sebastiani theorem. The cautious reader can always replace it by
 ${H}^{\bullet,crit}_{c,MMHS}$. 

Another unsatisfactory issue which we have already mentioned in Section
7.1  is that vanishing cycles are not yet defined for formal functions. Ideally we would like to have
a definition of the sheaf of vanishing cycles ``$\phi_f\Q_X$'' associated with any smooth formal scheme $\mathfrak{X}/\kk$ endowed with a function $f\in \mathcal{O}(\mathfrak{X})$. This  should be a  constructible   sheaf on the ordinary
reduced scheme $f^{-1}(0)\cup \mathfrak{X}_{red}$, where $\mathfrak{X}_{red}$ is obtained from $\mathfrak{X}$ by killing all topologically nilpotent elements.
 In the case when $\mathfrak{X}$ is a completion of an ordinary scheme $X$ along a closed subset $X^{(0)}\subset X$, and $f$ is a restriction of
$f_X\in \mathcal{O}(X)$,  the sheaf ``$\phi_f\Q_X$'' should be canonically isomorphic to the restriction of $\phi_{f_X} \Q_X$ to
$f_X^{-1}(0)\cap X^{(0)}$. Also we would like to have equivariant cohomology with compact support, and more generally cohomology of stacks.

There is a great variety of partial results which leave no doubts that such a theory should exist, see \cite{Berkovich1},
\cite{Berkovich2},\cite{NicaiseSebag},\cite{Jiang}.
 For a geometric theory formulated in terms of the Grothendieck group of varieties, one has good indications that the corresponding version of the Thom-Sebastiani theorem
 holds, see \cite{DenefLoeser}. Also, in the same spirit, see \cite{Temkin} for appropriate results on the equivariant resolution of singularities in
the formal setting.

\subsection{Rapid decay cohomology and monodromic Hodge structures}

Finally, we would like to return to the comparison with rapid decay cohomology, already discussed in Section 7.3, but this time at the level of Hodge structures. There is an exact faithful tensor functor $F: EMHS\to MHM^{\C-grad}$, where $MHM^{\C-grad}$ consists  of finite direct sums $\oplus_{z_i\in \C}M_{z_i}$ of monodromic mixed Hogde modules. The functor $F$ assigns to $M\in Ob(EMHS)$ the direct sum of the corresponding monodromic mixed Hodge modules over the finite set $S:=\{z_1,...,z_n\}$ of singular points of $M$. In order to define the summands $M_{z_i}$ we consider $M$ as an object of $MHM_{\AC}$.
Each $M_{z_i}$ is given by the formula similar to the one from Definition 27. Namely,
$$M_{z_i}=(\{z_i\}\times \mathbb{G}_{m,\C}\to \AC)_!\,\,\phi_{\frac{z-z_i}{u}}\,(\AC\times \mathbb{G}_{m,\C}\to \AC)^{\ast}M\,,$$
where the map $\{z_i\}\times \mathbb{G}_{m,\C}\to \AC$ is given by $(z_i,u)\mapsto u$.

There is an analog of the functor $F$ in the framework of non-commutative Hodge structures (see \cite{KaKoPa}). Namely, with a nc Hodge structure of exponential type one can associate a $\C$-graded nc Hodge structure with regular singularities by forgetting the gluing data (cf. Theorem 2.35 from loc. cit.).
In terms of $D$-modules the situation is described as forgetting the Stokes data, and hence passing to the formal classification of $D$-modules\footnote{Similar 
story appears in matrix integrals. In that case the critical cohomology corresponds to the formal expansion. For $1$-matrix integrals filling fractions play a role 
of the dimension of a representation of a quiver.}.  Notice that the functor $F$ preserves the weight filtration, and hence it preserves the Serre polynomial. In the case of the function $f$ which is isotrivial at infinity we have:
$$F(H_{c,EMHS}(X,f)^{\bullet})\simeq \bigoplus_{z\in Bif(f)}H_{c,MMHS}^{\bullet,crit}(f^{-1}(z),f-z)\,.$$

\subsection{Definition of the critical COHA}

Let us assume that we are in the situation described in Section 7.1, i.e. given smooth $I$-bigraded algebra $R/\kk$, a bilinear form $\chi_R$,  a potential $W$, and collection of subsets $\M_\gamma^\s$ satisfying the conditions from Section 7.1.
The  cohomology space $\mathcal{H}_\gamma$ is a $\Z$-graded object of the category of monodromic exponential mixed $\bf H$-motives
 defined as
$$ \mathcal{H}_\gamma:=D \left({\bf H}^{\bullet,crit}_{c,\G_\gamma}(\M_\gamma^\s,W_\gamma)
 \right) \otimes \mathbb{T}^{\otimes \dim \M_\gamma/\G_\gamma}\,.$$
 In the above formula we define
the equivariant cohomology with compact support  as the inductive limit, using finite-dimensional approximations to $\op{B}\G_\gamma$
 as in Section 4.5, and $D$ is  duality in the ind-completion of the category of monodromic exponential mixed $\bf H$-motives.
In what follows in order to alleviate the notation  we assume that potentials $W_\gamma$ have only one critical value $0$. The general case is
 similar.

 We have to define the (twisted) multiplication. For any $\gamma_1,\gamma_2\in \Z^I_{\geqslant 0},\,\,\gamma:=\gamma_1+\gamma_2$ we define the dual to the component $m_{\gamma_1,\gamma_2}$ of the product 
$$m_{\gamma_1,\gamma_2}^\vee: {\bf H}^{\bullet,crit}_{c,\G_\gamma}(\M_\gamma^\s,W_\gamma) \to {\bf H}^{\bullet,crit}_{c,\G_{\gamma_1}}(\M_{\gamma_1}^\s,W_{\gamma_1})\otimes
{\bf H}^{\bullet,crit}_{c,\G_{\gamma_2}}(\M_{\gamma_2}^\s,W_{\gamma_2}) \otimes \mathbb{T}^{\otimes\left(-\chi_R(\gamma_2,\gamma_1)\right)} $$
as the composition of the following homomorphisms of groups:
\begin{itemize}
\item  $ {\bf H}^{\bullet,crit}_{c,\G_\gamma}(\M_\gamma^\s,W_\gamma) \to     {\bf H}^{\bullet,crit}_{c,\G_{\gamma_1,\gamma_2}}(\M_\gamma^\s,W_\gamma)   $, which is the  pull-back associated with the embedding of groups $\G_{\gamma_1,\gamma_2}\to \G_{\gamma}$ with {\it proper} quotient.
\item $ {\bf H}^{\bullet,crit}_{c,\G_{\gamma_1,\gamma_2}}(\M_\gamma^\s,W_\gamma) \to
{\bf H}^{\bullet,crit}_{c,\G_{\gamma_1,\gamma_2}}(\M_{\gamma_1,\gamma_2}^\s,W_\gamma)     $, where $\M_{\gamma_1,\gamma_2}^\s:=\M_\gamma^\s\cap \M_{\gamma_1,\gamma_2}$, also given by the pullback for the {\it closed} embedding $\M_{\gamma_1,\gamma_2}^\s\mono\M_\gamma^\s$.
\item ${\bf H}^{\bullet,crit}_{c,\G_{\gamma_1,\gamma_2}}(\M_{\gamma_1,\gamma_2}^\s,W_\gamma)\simeq  {\bf H}^{\bullet,crit}_{c,\G_{\gamma_1,\gamma_2}}(\widetilde{\M}_{\gamma_1,\gamma_2}^\s,W_\gamma) $, where $\widetilde{\M}_{\gamma_1,\gamma_2}^\s\subset \M_{\gamma_1,\gamma_2}$
 is the pullback of $\M_{\gamma_1}^\s\times \M_{\gamma_2}^\s$ under the projection $\M_{\gamma_1,\gamma_2}\to \M_{\gamma_1}\times \M_{\gamma_2}$.
   The isomorphism follows from the fact that by the assumptions from Section 7.1 the space $\M_{\gamma_1,\gamma_2}^\s$ is equal to the intersection $\op{Crit}(W_\gamma)\cap
  \widetilde{\M}_{\gamma_1,\gamma_2}^\s$. Hence the sheaf of vanishing cycles of $W_\gamma$ vanishes on $\widetilde{\M}_{\gamma_1,\gamma_2}^\s
- \M_{\gamma_1,\gamma_2}^\s$.
\item  $ {\bf H}^{\bullet,crit}_{c,\G_{\gamma_1,\gamma_2}}(\widetilde{\M}_{\gamma_1,\gamma_2}^\s,W_\gamma)\to
 {\bf H}^{\bullet,crit}_{c,\G_{\gamma_1,\gamma_2}}(\widetilde{\M}_{\gamma_1,\gamma_2}^\s,W_{\gamma_1,\gamma_2})   $, where $W_{\gamma_1,\gamma_2}$ is the restriction of $W_\gamma$ to $\M_{\gamma_1,\gamma_2}$.
\item $ {\bf H}^{\bullet,crit}_{c,\G_{\gamma_1,\gamma_2}}(\widetilde{\M}_{\gamma_1,\gamma_2}^\s,W_{\gamma_1,\gamma_2}) \simeq
   {\bf H}^{\bullet,crit}_{c,\G_{\gamma_1}\times\G_{\gamma_2}}({\M}_{\gamma_1}^\s\times\M_{\gamma_2}^\s,W_{\gamma_1}\boxplus W_{\gamma_2})\otimes
\mathbb{T}^c$, where the isomorphism comes from homotopy equivalence $\G_{\gamma_1}\times\G_{\gamma_2}\sim \G_{\gamma_1,\gamma_2}$, the fact that
$\widetilde{\M}_{\gamma_1,\gamma_2}^\s$ is a bundle over ${\M}_{\gamma_1}^\s\times\M_{\gamma_2}^\s$ with affine fibers, and the fact that $W_{\gamma_1,\gamma_2}$
 is the pullback of $W_{\gamma_1}\boxplus W_{\gamma_2}$.
 The shift is given by
 $$c=\dim\M_{\gamma_1,\gamma_2}/\G_{\gamma_1,\gamma_2}-\dim\M_{\gamma_1}/\G_{\gamma_1}-\dim\M_{\gamma_2}/\G_{\gamma_2}= -\chi_R(\gamma_2,\gamma_1)
\,.$$
\item   Thom-Sebastiani isomorphism.

\end{itemize}
The proof of associativity  is similar to the one from Section 2.3.

\subsection{Stability conditions and factorization of critical motivic DT-series}

We will use the notation for the Harder-Narasimhan filtration from  Section 5.1. Again, we assume for simplicity that all critical values of $W_\gamma$ are equal to zero.
 Stratification by $HN$-strata gives a spectral sequence converging to ${\bf H}^{\bullet,crit}_{c,\G_\gamma}(\M_\gamma^\s,W_\gamma)$ with the first term
$$ \bigoplus_{n\geqslant 0}\bigoplus_{\substack{\gamma_1,\dots,\gamma_n\in \Z_{\geqslant 0}^I- 0\\
\op{Arg}\gamma_1>\dots>\op{Arg}\gamma_n}  } {\bf H}^{\bullet,crit}_{c,\G_\gamma}(\M_{\gamma,\gamma_\bullet}^{\s,HN},W_\gamma)\,,$$
where $\gamma_\bullet=(\gamma_1,\dots,\gamma_n)$ and $\M_{\gamma,\gamma_\bullet}^{\s,HN}:=\M_{\gamma,\gamma_\bullet}^{\s}\cap
\M_{\gamma,\gamma_\bullet}^{HN}$.
Then we have
 $${\bf H}^{\bullet,crit}_{c,\G_\gamma}(\M_{\gamma,\gamma_\bullet}^{\s,HN},W_\gamma)\simeq
{\bf H}^{\bullet,crit}_{c,\G_{\gamma_1,\dots,\gamma_n}}(\M_{\gamma_1,\dots,\gamma_n}^{\s,ss},W_\gamma)\,,\,\,\,\M_{\gamma_1,\dots,\gamma_n}^{\s,ss}:=
\M_{\gamma}^{\s}\cap\M_{\gamma_1,\dots,\gamma_n}^{ss}\,.$$
Denote by $\widetilde{\M}_{\gamma_1,\dots,\gamma_n}^{\s,ss}$ the pullback of $\M_{\gamma_1}^{\s,ss}\times\dots\times \M_{\gamma_n}^{\s,ss}$
 under the projection $\M_{\gamma_1,\dots,\gamma_n}\to \M_{\gamma_1}\times\dots\times \M_{\gamma_n}$.
 We have an isomorphism
 $${\bf H}^{\bullet,crit}_{c,\G_{\gamma_1,\dots,\gamma_n}}(\M_{\gamma_1,\dots,\gamma_n}^{\s,ss},W_\gamma)\simeq
{\bf H}^{\bullet,crit}_{c,\G_{\gamma_1,\dots,\gamma_n}}(\widetilde{\M}_{\gamma_1,\dots,\gamma_n}^{\s,ss},W_\gamma)\,,$$
because $W_\gamma$ has no critical points on $\widetilde{\M}_{\gamma_1,\dots,\gamma_n}^{\s,ss}-\M_{\gamma_1,\dots,\gamma_n}^{\s,ss}$.
The next step is an isomorphism
 $$ {\bf H}^{\bullet,crit}_{c,\G_{\gamma_1,\dots,\gamma_n}}(\widetilde{\M}_{\gamma_1,\dots,\gamma_n}^{\s,ss},W_\gamma)
\simeq
{\bf H}^{\bullet,crit}_{c,\G_{\gamma_1,\dots,\gamma_n}}(\widetilde{\M}_{\gamma_1,\dots,\gamma_n}^{\s,ss},(W_\gamma)_{|\M_{\gamma_1,\dots,\gamma_n}})\,.$$
This is the key point. This isomorphism follows iteratively from the integral identity proven in the next section, see Corollary 5.
Finally, using the fact that $\widetilde{\M}_{\gamma_1,\dots,\gamma_n}^{\s,ss}$ is an affine bundle over $ \M_{\gamma_1}^{\s,ss}\times\dots\times \M_{\gamma_n}^{\s,ss}$ and that the function
$(W_\gamma)_{|\M_{\gamma_1,\dots,\gamma_n}}$ is the pullback of $W_{\gamma_1}\boxplus\dots \boxplus  W_{\gamma_n}$, we obtain that
$$ {\bf H}^{\bullet,crit}_{c,\G_\gamma}(\M_{\gamma,\gamma_\bullet}^{\s,HN},W_\gamma)\simeq  $$
$$\simeq
 {\bf H}^{\bullet,crit}_{c,\G_{\gamma_1}\times \dots \times \G_{\gamma_n}}({\M}_{\gamma_1}^{\s,ss}\times \dots\times\M_{\gamma_n}^{\s,ss},W_{\gamma_1}\boxplus \dots\boxplus  W_{\gamma_n})  \otimes \mathbb{T}^{\otimes (-\sum_{i<j} \chi_R(\gamma_j,\gamma_i))}\,. $$
 Applying Thom-Sebastiani isomorphism and duality we obtain the Factorization Formula for the critical motivic DT-series, as in Section 5.3.

\subsection{Integral identity}

Let $X$ be an algebraic variety over $\C$, and $\E_1,\E_2$ be two vector bundles over $X$. Denote by $Y$ the total space of the bundle
$\E_1\oplus \E_2$, and by $Y^{(1)}\subset Y$ the total space of $\E_1$. Assume that we are given a function
 $$f\in {\mathcal{O}}(Y)^{{\mathbb{G}}_m}\,,$$
 where the group ${\mathbb{G}}_m$ acts linearly along fibers of $\E_1\oplus \E_2$ with weights
$+1$ on $\E_1$ and $-1$ on $\E_2$.

\begin{thm} Let $Y_0=f^{-1}(0)\subset Y$ be the zero locus of $f$ and $Y_0^{(1)}:=Y_0\cap Y^{(1)}$. Then
$$(Y_0^{(1)}\to X)_!\, (Y_0^{(1)}\to Y_0)^*\, \op{Cone}\left(\phi_f \Z_Y\to \phi_f \left(Y^{(1)}\to Y)_*\Z_{Y^{(1)}}\right)\right)=0\,.$$
\end{thm}

\begin{cor} For any closed subset $X^\s\subset X$ such that its preimage $Y^\s$ under the projection $Y^{(1)}\to X$ lies in $Y_0$, we have
 a canonical isomorphism
$${\bf H}^{\bullet,crit}_c (Y^\s,f)\simeq {\bf H}^{\bullet,crit}_c (Y^\s,f_{|Y^{(1)}})\,.$$
Moreover, the same conclusion holds in the equivariant setting.
\end{cor}
{\it Proof of the corollary}.
 There is an obvious morphism in $D^b(MMHS)$
$$ \alpha:  H^{\bullet,crit}_{c,MMHS}(Y^\s,f)\to  H^{\bullet,crit}_{c,MMHS}(Y^\s,f_{|Y^{(1)}})   \,\,$$
We claim that $\phi$ is an isomorphism, i.e. the cone of $\alpha$ is zero.
By faithfulness, it is sufficient to prove the vanishing using only Betti realization $\alpha_{\op{Betti}}$.
Denote by $\mathcal{E}$ the complex of constructible sheaves whose vanishing is claimed in Theorem 13. Then
we have
$$ (X^\s\to pt)_! (X^\s\to X)^* \mathcal {E}\otimes \Q\simeq \op{Cone} (\alpha_{\op{Betti}})$$
Hence we proved the result in the non-equivariant case. The extension to the equivariant setting is straightforward.

{\it Proof of Theorem 13}.
Let us consider  an affine morphism $Z\to X$ whose fiber $Z_x$ at any point $x\in X$ is $\op{Spec}\mathcal{O}((Y\to X)^{-1}(x))^{{\mathbb{G}}_m}$.
There is an obvious map $Y\to Z$, and the function $f$ is the pullback of a function $f_Z\in \mathcal{O}(Z)$.
We will construct a completion $\overline{Y}\supset Y$ such that the projection $Y\to Z$ extends to a {\it proper} map $\pi:\overline{Y}\to Z$.
  In what follows we will compose $\pi$ with the projection $Z\to X$ and describe $\overline{Y}$ as a fibration over $X$.
 Namely, for any $x\in X$ consider the product of two projective spaces
$$\mathbb{P}(\E_{1,x}\oplus \C)\times \mathbb{P}(\E_{2,x}\oplus\C)\simeq
\left(\E_{1,x}\cup \mathbb{P}\E_{1,x}\right)\times \left(\E_{2,x}\cup \mathbb{P}\E_{2,x}\right)\,,$$
make blow-up at two disjoint submanifolds $\mathbb{P}\E_{1,x}\times \{0_{\E_{2,x}}\}$ and $\{0_{\E_{1,x}}\}\times \mathbb{P}\E_{2,x}$
 (here $0_{\E_{1,x}},0_{\E_{2,x}} $ are zero points in fibers $\E_{1,x}$ and $\E_{2,x}$ respectively), and remove the proper transform
 of divisors $\mathbb{P}\E_{1,x}\times \E_{2,x}$ and $\E_{1,x}\times \mathbb{P}\E_{2,x}$. The  result is by definition the fiber
 of $\overline{Y}$ over $x$. This fiber is smooth, and carries a natural stratification with 8 smooth strata. The incidence diagram of strata is
$$\xymatrix{ S_{7,x}\ar[d]\ar[r] & S_{8,x}\ar[d] & \\
              S_{4,x}\ar[r] & S_{5,x} & S_{6,x}\ar[l]\\
              S_{1,x}\ar[r]\ar[u] & S_{2,x}\ar[u] & S_{3,x}\ar[l]\ar[u]
}$$
 where all arrows are inclusions. Open parts of the strata are
$$\xymatrix{  \mathbb{P}\E_{2,x} & S_{8,x} & \\
\E_{2,x}^{\ne 0} &   \E_{1,x}^{\ne 0}\times   \E_{2,x}^{\ne 0}  & S_{6,x} \\
\{0_{\E_{1,x}\oplus \E_{2,x}}\} & \E_{1,x}^{\ne 0} & \mathbb{P}\E_{1,x}
}$$
where $\E_{i,x}^{\ne 0}:=\E_{i,x}-\{0_{\E_{i,x}}\},\,\,i=1,2$ and $S_{6,x}\simeq S_{8,x}\simeq \left( \E_{1,x}^{\ne 0}\times   \E_{2,x}^{\ne 0}  \right)/\mathbb{G}_{m,\C}$.
Under the projection to the fiber $Z_x$ of $Z\to X$ the strata $S_{1,x},S_{2,x},S_{3,x},S_{4,x},S_{7,x}$ are mapped to the base point in $Z_x$,
 whereas the union $S_{5,x}\cup S_{6,x}\cup S_{8,x}$ is a fibration with fiber $\mathbb{P}^1$ over an open stratum in $Z_x$ isomorphic to $S_{6,x}$.
 Varying point $x\in X$ we obtain strata $(S_i)_{1\leqslant i\leqslant 8}$ of $\overline{Y}$ with open parts $(S_i^\circ)_{1\leqslant i\leqslant 8}$. The open subset $Y\subset
\overline{Y}$ is a union of strata, $Y=S_1^\circ\sqcup S_2^\circ\sqcup S_4^\circ\sqcup S_5^\circ$.

Now return to the proof of identity. By base change it is sufficient to prove that for any $x\in X$ such that $f_X(x)=0$ (here $f_X$ is the pullback
 of $f$ via the embedding $X\mono Y$ as the zero section) we have
$$(\E_{1,x}\to pt)_!\, (\E_{1,x}\to Y_0)^*\, \phi_f\,\mathcal{G}^\bullet=0,\,\,\,\mathcal{G}^\bullet:=\op{Cone}\left( \Z_Y\to  (Y^{(1)}\to Y)_*\Z_{Y^{(1)}} \right).$$
We will construct certain extension $\mathcal{F}^{\bullet}$ of $\mathcal{G}^\bullet$ from $Y$ to
 $\overline{Y}$. Let $\mathcal{F}_1$ be the extension of $\Z_Y$ by $!$ to $S_3^\circ\sqcup S_6^\circ$, and by
  $*$ to $S_7^\circ\sqcup S_8^\circ$. Similarly we define $\mathcal{F}_2$ to be the extension of $\Z_{Y^{(1)}}$  by $!$ to  $\overline{Y}$. Restrictions of $\mathcal{F}_1$ and $\mathcal{F}_2$ to $Y\subset \overline{Y}$ coincide with the sheaves
 $\Z_Y$ and $(Y^{(1)}\to Y)_*\Z_{Y^{(1)}}$ respectively.
 There is a natural morphism $\mathcal{F}_1\to \mathcal{F}_2$, and the restriction of $\mathcal{F}^\bullet:=\op{Cone}(\mathcal{F}_1\to \mathcal{F}_2)$ to $Y$
 is isomorphic to $\mathcal{G}^\bullet$.

\begin{lmm} The complexes
$$(\E_{1,x}\to pt)_!\, (\E_{1,x}\to Y_0)^*\, \phi_f\,\mathcal{G}^\bullet\simeq (\E_{1,x}\to pt)_!\, (\E_{1,x}\to \overline{Y}_0)^*\,
\phi_{f_{\overline{Y}}}\,\mathcal{F}^\bullet$$ and
$$(\overline{Y}_x'\to pt)_!\, (\overline{Y}_x'\to \overline{Y}_0)^* \,\phi_{f_{\overline{Y}}}\, \mathcal{F}^\bullet     $$
are isomorphic. Here $\overline{Y}_x':=S_{1,x}^\circ \sqcup S_{2,x}^\circ\sqcup S_{3,x}^{\circ}\sqcup S_{4,x}^\circ\sqcup S_{7,x}^\circ$ is a closed subset of the fiber of $\overline{Y}$ over $x$, function $f_{\overline{Y}}=(\overline{Y}\to Z)^* f_Z$ is the extension of $f$ by continuity,
 and $\overline{Y}_0$ is zero locus of $f_{\overline{Y}}$.
\end{lmm}
{\it Proof of Lemma}.  Notice that $\overline{Y}_x'=\E_{1,x}\sqcup S_{3,x}^\circ\sqcup S_{4,x}^\circ\sqcup S_{7,x}^\circ$. By spectral sequences it is sufficient to check that
 \begin{enumerate}
\item $(S_{3,x}^\circ\to pt)_! \,(S_{3,x}^\circ\to \overline{Y}_0)^* \,\phi_{f_{\overline{Y}}}\, \mathcal{F}^\bullet =0, $
\item $((S_{4,x}^\circ\sqcup S_{7,x}^\circ)\to pt)_! \,((S_{4,x}^\circ\sqcup S_{7,x}^\circ)\to \overline{Y}_0)^* \,\phi_{f_{\overline{Y}}}\, \mathcal{F}^\bullet =0.$
\end{enumerate}
First vanishing follows from the observation that the closed subspace $V:=S_{3}^\circ\sqcup S_{6}^\circ$ of $\overline{Y}$ has an open neighborhood
  $U:=S_{2}^\circ\sqcup S_{3}^\circ\sqcup S_{5}^\circ\sqcup S_{6}^\circ$ which is the total space of a line bundle over $S_{3}^\circ\sqcup S_{6}^\circ$. Restriction of function $f_{\overline{Y}}$
to $U$ is the pullback of a function on $V$, and the restriction of the sheaf $\mathcal{F}^\bullet$ to $U$ is the extension by zero of
 the constant sheaf on $U-V$. Hence the functor of vanishing cycles  can be interchanged with the restriction  to $V$, and therefore
  $$\left(((S_3^\circ\sqcup S_6^\circ)\cap \overline{Y}_0) \to \overline{Y}_0\right)^* \,\phi_{f_{\overline{Y}}}\, \mathcal{F}^\bullet =0    \,.$$
This implies 1).

 Similarly, $V':=S_{7}^\circ\sqcup S_8^\circ$ is a base of a line bundle with the total space $U'=S_4^\circ\sqcup S_5^\circ\sqcup S_7^\circ \sqcup S_8^\circ$. Restriction of function $f_{\overline{Y}}$
to $U'$ is the pullback of a function on $V'$, and the restriction of the sheaf $\mathcal{F}^\bullet$ to $U$ is the direct image of
 the constant sheaf on $U'-V'$. Hence  we have
 $$  (U'_0 \to V'_0)_!\, (U'_0\to \overline{Y}_0)^*\,  \mathcal{F}^\bullet=0\,,\,\,\,\,\,U_0':= U\cap\overline{Y}_0,\,V_0':=V\cap \overline{Y}_0  \,.$$
The reason is that $(\C\to pt)_! \,(\C^*\to \C)_*\, \Z_{\C^*}=0$. This proves the second vanishing.
$\blacksquare$

Using the fact that the functor of vanishing cycles commutes with proper morphisms (e.g. with $\pi:\overline{Y}\to Z$), we obtain an isomorphism
$$(\overline{Y}_x'\to pt)_!\, (\overline{Y}_x'\to \overline{Y}_0)^* \,\phi_{f_{\overline{Y}}}\, \mathcal{F}^\bullet
\simeq (pt \to f_Z^{-1}(0))^*\,\phi_Z  \left( \overline{Y}\to Z\right)_* \, \mathcal{F}^\bullet, $$
where the inclusion of the point to $ f_Z^{-1}(0)$ is given by the base point at the fiber $Z_x\subset f_Z^{-1}(0)$.
Now we apply the following lemma.
\begin{lmm} The direct image $\left( \overline{Y}\to Z\right)_* \, \mathcal{F}^\bullet $ vanishes.
\end{lmm}
{\it Proof of Lemma}. The morphism $\pi:\overline{Y}\to Z$ is proper, hence we can use the base change, and it is sufficient to prove that
 for any point $z \in Z$ the direct image
$$(\pi^{-1}(z)\to pt)_*\,\, (\pi^{-1}(z)\to \overline{Y})^*  \, \mathcal{F}^\bullet=0\,. $$

Denote by $i: X\mono Z$  the canonical embedding given by base points of cones $(Z_x)_{x\in X}$.

 For $z\in Z-i(X)$  the fiber $\pi^{-1}(z)$
 is the projective line $\mathbb{P}^1_\C$. 
The  restriction of the sheaf $\mathcal{F}^\bullet$ to $\pi^{-1}(z)$ is a constant sheaf on $\mathbb{A}^1_\C$ extended by zero to the
point $\infty \in \mathbb{P}^1_\C$. Hence $R\Gamma$ of this sheaf on $\pi^{-1}(z)$ vanishes.

For a point $z\in i(X)$ , the fiber $\pi^{-1}(z)$ is the joint of two projective spaces identified at zero:
$$    \pi^{-1}(z)=  \mathbb{P}(\E_{1,x}\oplus \C)\cup_{0} \mathbb{P}(\E_{2,x}\oplus\C)=S_{1,x}^\circ\sqcup S_{2,x}^\circ \sqcup S_{3,x}^\circ \sqcup
 S_{4,x}^\circ \sqcup S_{7,x}^\circ\,,$$
where the point $x\in X$ corresponds to $z$.
The restriction of $\mathcal{F}^\bullet$ to $\pi^{-1}(z)$ is the constant sheaf on the  open part $S_{1,x}^\circ \sqcup S_{2,x}^\circ\sqcup S_{4,x}^\circ$
 extended by zero to $S_{3,x}^\circ$, and by the direct image to $S_{7,x}^\circ$. Then a direct calculation shows that $R\Gamma$ of $\mathcal{F}^\bullet_{|\pi^{-1}(z)}$ vanishes.
This proves the lemma. $\blacksquare$

Two above lemmas  imply that for any point $x\in {f_X}^{-1}(0)$ we have
$$(\E_{1,x}\to pt)_!\, (\E_{1,x}\to Y_0)^*\, \phi_f\,\mathcal{G}^\bullet=0\, . $$
By base change this implies the Theorem. $\blacksquare$

\subsection{Factorization systems in the critical case}

The proof of the integrality property is based on the theory of factorization systems from Sections 6.5 and 6.6\footnote{We do not know how to make a shortcut similar to the one in Section 6.7.}.

First, for $\gamma\in \Z^I_{\geqslant 0}$  we define  the following graded space:
$$ \mathcal{H}_\gamma':=D \left({\bf H}^{\bullet,crit}_{c,\T_\gamma}(\M_\gamma^\s,W_\gamma)
 \right) \otimes \mathbb{T}^{\otimes \dim \M_\gamma/\T_\gamma}\,.$$
Our goal is to define a structure of a non-singular factorization system on the collection $( \mathcal{H}_\gamma')_{\gamma\in \Z^I_{\geqslant 0}}$.
 Then applying Theorem 12 and Corollary 4 from Section 6.6 we will deduce admissibility of the motivic Donaldson-Thomas series in the critical case.

First, we claim that $ \mathcal{H}_\gamma'$ 
 is a {\it finitely generated} module over ${\bf H}^{\bullet}(\B\T_\gamma)$. This is a particular case of a more general result proven below, where we identify cohomology with coefficients in the  sheaf of vanishing cycles with cohomology of certain real semi-algebraic sets.

We define the multiplication map (for $\gamma=\gamma_1+\gamma_2$) as the dual to the composition
$$ {\bf H}^{\bullet,crit}_{c,\T_\gamma}(\M_\gamma^\s,W_\gamma)\to {\bf H}^{\bullet,crit}_{c,\T_\gamma}(\M_{\gamma_1}^\s\times \M_{\gamma_2}^\s,W_\gamma)      \to {\bf H}^{\bullet,crit}_{c,\T_\gamma}(\M_{\gamma_1}^s\times \M_{\gamma_2}^\s,W_{\gamma_1}\boxplus W_{\gamma_2})\,,$$ 
where the first arrow is the restriction to the closed subspace $\M_{\gamma_1}^\s\times \M_{\gamma_2}^\s\subset \M_\gamma^\s$, and the second arrow
 comes from the canonical morphism of sheaves. It is clear that the constraints of commutativity, associativity, and of the equivariance of the product are satisfied.

\begin{thm}
In the notation from Section 6.4 the multiplication map induces an isomorphism of coherent sheaves on the open sets $U_{\gamma_1,\gamma_2}$ .
\end{thm}

{\it Proof} \footnote{ We thank to Joerg Sch\"urmann for a discussion which led to the  proof below.}.  We will use Betti realization of the critical cohomology cohomology (with coefficients in $\Q$).
  We have to prove that the dual $\Z$-graded vector space over $\Q$ to
$$H^\bullet_{c,\T_\gamma}\left(\M^\s_\gamma,\left(\M^\s_\gamma\mono W_\gamma^{-1}(0)\right)^*\phi_{W_\gamma}\left(\op{Cone}(\Q_{\M_\gamma} \to
 i_*\Q_{\M_{\gamma_1}\times \M_{\gamma_2}} )\right)   \right),$$
where $i:\M_{\gamma_1}\times \M_{\gamma_2}\to \M_\gamma$ is the obvious embedding, is a finitely generated ${H}^{\bullet}(\B\T_\gamma,\Q)$-module
with support on $\op{Spec}{H}^{\bullet}(\B\T_\gamma,\Q)-U_{\gamma_1,\gamma_2}$.

First, notice that in the above formula we can replace complex torus \begin{large}$\T_\gamma$                                                                         \end{large} by its maximal compact real subtorus $\T_\gamma^c$ (hence the shift is replaced by the real dimension
of the contractible group $\T_\gamma/\T^c_\gamma$).

Also  notice that we can replace the sheaf of vanishing cycles by the sheaf of nearby cycles. 
Indeed, these two choices  (i.e. $\phi_{W_\gamma}$ and $\psi_{W_\gamma}$) are related by the exact triangle from Section 7.2 with the cohomology of the third vertex of the triangle given by
$$H^\bullet_{c,\T_\gamma}\left(\M^\s_\gamma, (\M^\s_\gamma\to \M_\gamma)^*\left(\op{Cone}(\Q_{\M_\gamma} \to
 i_*\Q_{\M_{\gamma_1}\times \M_{\gamma_2}} )\right)\right)\simeq H^\bullet_{c,\T_\gamma}\left(\M^\s_\gamma-(\M_{\gamma_1}^\s\times \M_{\gamma_2}^\s)\right).$$
Then we can use   Proposition 11 from Section 6.5.

Now we have reduced the question to the one about equivariant cohomology with coefficients in the sheaves of nearby cycles.
Let us consider the following general situation. Suppose we are given a complex smooth algebraic variety $X$ endowed with an algebraic action of a complex torus $\T$ as well as 
a $\T$-invariant function $f\in \mathcal{O}(X)$, and a closed $\T$-invariant algebraic set $S\subset f^{-1}(0)$.
 We want to calculate
 $$H_{c,\T}(S,(S\to f^{-1}(0))^* \psi_f \Q_X)\,.$$
As above, we can replace $\T$ by the maximal compact real subtorus $\T^c\subset \T$.
Let us choose $\T^c$-invariant continuous  real semi-algebraic maps $r,s:X\to [0,+\infty)$  such that: 

1) the map $r$ is proper, 

2)  $d^{-1}(0)=S$ (the map $d$ can be thought of as  a ``distance to $S$''). 

Then we claim that $H_{c,\T^c}(S,(S\to f^{-1}(0))^* \psi_f \Q_X)$ is isomorphic to the iterated limit\footnote{More precisely, we should work in the derived category and take {\it derived} projective limits. In our case all cohomology spaces will have bounded finite dimension, hence the derived projective limit coincides with the
 usual one.}
$$\varinjlim_{r_1\to +\infty}\,\varprojlim_{r_2\to +\infty}\,\varinjlim_{d_1\to +0}\,\varinjlim_{\epsilon_1\to +0}\,\varinjlim_{\epsilon_2\to+0}\,\varprojlim_{\xi_1\to +0}\,\varprojlim_{\eta_1\to +0}\,\varprojlim_{\eta_2\to +0}\,
\varprojlim_{f_1\to +0} H^\bullet_{\T^c}(S_1,S_2;\Q)\,,
$$
where $(S_1,S_2)$ is the following pair of $\T^c$-invariant compact real semi-algebraic sets in $X$ depending on positive real parameters $r_1,r_2,d_1,\epsilon_1,\epsilon_2,\xi_1,\eta_1,\eta_2,f_1$:
$$S_1=\{x\in X\,|\,r(x)\leqslant r_2+\epsilon_2-\eta_2,\,f(x)\geqslant f_1,\,d(x)\leqslant d_1-\xi_1\}\,,$$
$$S_2=\{x\in X\,|\,r_1-\epsilon_1+\eta_1\leqslant r(x)\leqslant r_2+\epsilon_2-\eta_2,\,f(x)\geqslant f_1,\,d(x)\leqslant d_1-\xi_1\}\,.$$

This description of cohomology is similar to ``semi-global'' Milnor fibrations from \cite{Schuermann}.

Let us explain the above sequence of direct and inverse limits. First pair of limits $\varinjlim_{r_1\to +\infty}\,\varprojlim_{r_2\to +\infty}$ corresponds to the fact that
 the cohomology of $S$ with compact support with coefficient in {\it any} sheaf can be calculated as the double limit of cohomology of pair $(K_1,K_2)$ of compact sets depending on parameters $r_1,r_2$ given by
$$K_1:=\{x\in S\,|\, r(x)\leqslant r_2 \}\,,\,\,\,\,K_2:=\{x\in S\,|\, r_1\leqslant r(x)\leqslant r_2\}\,.     $$
By definition of $\psi_f \Q_X$, we should calculate cohomology of  pair $(K_1,K_2)$ with coefficients in the sheaf $(X_+\to X)_* \Q_{X_+}$ on $X$, where $X_+:=f^{-1}(\R_{>0})$.
Then we consider a  system of open neighborhoods $U_1,U_2$ of $K_1,K_2$ in $X$, depending on small positive parameters $d_1,\epsilon_1,\epsilon_2$ and given by 
$$U_1:=\{x\in X\,|\,r(x)< r_2+\epsilon_2),\,d(x)< d_1\}\,,$$
$$U_2:=\{x\in X\,|\,r_2-\epsilon_1< r(x)< r_2+\epsilon_2,\,d(x)< d_1\}\,.$$
This is a {\it fundamental} system of neighborhoods (as follows from the properness of the map $(r,d):X\to \R^2$). Hence the cohomology of the pair $(K_1,K_2)$ is the inductive limit of the cohomology of  the pair $(U_1,U_2)$ as $d_1,\epsilon_1,\epsilon_2\to +0$. 
Now we should calculate the cohomology  $H^\bullet(U_1\cap X_+,U_2\cap X_+;\Q)$.
  We approximate  locally compact spaces $U_i\cap X_+,\,i=1,2$ from inside by a growing family of compact subsets $S_i,\,i=1,2$ depending on small
 positive parameters $\xi_1,\eta_1,\eta_2,f_1$. This finishes the explanation of  the multiple limit formula for $H_{c,\T^c}(S,(S\to f^{-1}(0))^* \psi_f \Q_X)$.
  One can simplify slightly this formula, by reducing to the case 
$$d_1=\epsilon_1=\epsilon_2,\,\,\,\xi_1=\eta_1=\eta_2=f_1  \,. $$

The above cohomology  depends on all parameters in a constructible way in the  sense of the standard real semi-algebraic structure on $\R$.
 This is clear for non-equivariant cohomology. In the case of equivariant cohomology we should use the stratification given by the dimension of stabilizer in $\T^c$ of a point $x\in X$. 
Therefore, 
 at every stage the limit stabilizes for sufficiently small values of the parameters $1/r_1,1/r_2,d_1,\epsilon_1,\epsilon_2,\xi_1,\eta_1,\eta_2,f_1$. Hence, the space $H_{c,\T^c}(S,(S\to f^{-1}(0))^* \psi_f \Q_X)$ coincides with the cohomology with compact support of the semi-algebraic set $S_1-S_2$ for certain values of parameters,
 and therefore the $\Z$-graded dual space is a finitely generated $H^\bullet(\B\T^c,\Q)$-module.

Now we can apply the iterated limit formula in {\it two} cases: one with
$$\T=\T_\gamma,\,X=\M_\gamma,\,f=W_\gamma,\,S=\M^\s_\gamma\,\,\,\mbox{ and another one with}$$
$$\T'=\T=\T_\gamma,\,X'=\M_{\gamma_1}\times \M_{\gamma_2}\subset X,\,f=W_{\gamma_1}\boxplus W_{\gamma_2}=f_{|X'},\,
S'=\M^\s_{\gamma_1}\times \M^\s_{\gamma_2}=S\cap X'.$$
We  also choose the functions $r,d$ on $X$ as above, and define $r',d'$ to be the restrictions of $r,d$ to $X'$.
Then we have equalities 
 $$S_1'=S_1\cap X',\,\,S_2'=S_2\cap X'$$
and a long exact sequence
$$\dots \to H^\bullet_\T(S_1,S_2\cup_{S_2'} S_1';\Q)\to H^\bullet_\T(S_1,S_2;\Q)\to H^\bullet_\T(S_1',S_2';\Q)\to \dots\,\,\,. $$
For sufficiently small values of parameters as above we can replace $H^\bullet_\T(S_1,S_2;\Q)$ and $ H^\bullet_\T(S_1',S_2')$ by 
$$H^\bullet_{c,\T_\gamma}\left(\M^\s_\gamma,\left(\M^\s_\gamma\mono W_\gamma^{-1}(0)\right)^*\phi_{W_\gamma}\Q_{\M_\gamma}\right)\,\,\,\mbox{and}$$
$$H^\bullet_{c,\T_\gamma}\left(\M^\s_{\gamma_1}\times \M^\s_{\gamma_2},\left(\M^\s_{\gamma_1}\times \M^\s_{\gamma_2}\mono (W_{\gamma_1}\boxplus W_{\gamma_2})^{-1}(0)\right)^*\phi_{W_{\gamma_1}\boxplus W_{\gamma_2}}\Q_{\M_{\gamma_1}\times \M_{\gamma_2}}\right)$$
respectively.

The equivariant cohomology $ H^\bullet_\T(S_1,S_2\cup_{S_2'} S_1';\Q)$ coincides with the equivariant cohomology with compact support
 $$H^\bullet_{c,\T}\left(S_1-(S_2\cup_{S_2'} S_1'),\Q\right)$$ of
 a locally closed real semi-algebraic set which is {\it contained} in $X-X'=\M_\gamma-(\M_{\gamma_1}\times \M_{\gamma_2})$. Applying  Proposition 11 we finish the proof of  theorem.
 $\blacksquare$

\subsection{Comparison with motivic DT-invariants from \cite{KS}}

Results of Sections 7.1-7.9 ensure the existence of the critical COHA $\mathcal{H}^{crit}=\oplus_{\gamma}\mathcal{H}^{crit}_{\gamma}$ as well as the corresponding graded space $\mathcal{H}^{crit}_V$ associated with a sector $V$ in the upper-half plane. Hence we have the corresponding critical DT-series $A_V^{crit}$ as well as critical motivic DT-invariants $\Omega^{crit,mot}(\gamma)$. In this subsection we would like to compare $A_V^{crit}$ with the series $A_V^{mot}$ introduced in \cite{KS}.

We will be  sketchy. 
We continue to use the notation introduced in Section 7.1. 
 We will also assume that the reader is familiar with the set-up of the paper \cite{KS}.

 Let us make a simplifying assumption\footnote{All the conclusion of this section hold in the general case of $I$-bigraded smooth algebra.} that $R/\kk$ is the path algebra of a finite quiver $Q$ with the set of vertices $I$. Also we assume that the ground field is algebraically closed, $\kk=\overline{\kk}$. Then the potential $W$ gives rise to a homologically smooth dg $3CY$-algebra $B/\kk$ with  generators in degrees $-2,-1,0$ (Ginzburg algebra), and we can consider triangulated dg-category ${\cal C}={\cal C}_{(Q,W)}$ of finite-dimensional (over $\kk$) $A_\infty$-modules over $B$.
 It is an ind-constructible triangulated 3-dimensional Calabi-Yau $\kk$-linear category in the sense of \cite{KS}. It has a canonical bounded $t$-structure with the heart equivalent to 
 the abelian category $(\op{Crit}W)({\kk})$ in the notation from Section 7.1. Define ${\cal C}^\s$ to be the full triangulated subcategory of ${\cal C}$ generated by objects from $\sqcup_{\gamma}\M^\s_\gamma(\kk)$. It is also an ind-constructible triangulated category. Any choice of a central charge $Z:\Z^I\to \C$ as in Section 5.1 gives an ind-constructible stability condition on ${\cal C}^s$.

Calabi-Yau category ${\cal C}^\s$ has a canonical orientation data in the sense of \cite{KS}, given by the super line 
$$sdet\left(\op{Ext}^\bullet_{R-\op{mod}}(E_1,E_2)\right)$$
associated with any two objects $E_1,E_2\in {\cal C}^\s$. Here we us the fact that $R$ is a subalgebra of $B$.

Let us fix a central charge $Z$ and a sector $V$ in the upper-half plane.
 Then, as we recalled above, the theory developed in \cite{KS} gives rise to a series $A^{mot}_V$.
 This is a series in (a completion of) the quantum torus with coefficients lying in a localization of the Grothendieck ring $\mathcal{M}^\mu(\op{Spec}(\kk))$ of varieties over $\kk$ endowed with a good action of the pro-algebraic group ${\mu}=\varprojlim \mu_n$ where $\mu_n$ is the finite group scheme of $n$-th roots of one.
   For $\kk=\C$ there is a natural homomorphism of rings
$$\Phi:\mathcal{M}^\mu(\op{Spec}(\C))\to K_0(MMHS_{\AC})\,,$$
which associates with $X$ (endowed with a $\mu_n$-action for some $n\ge 1$) the $K$-theory class of the following object in $D^b(MHM_{\AC})$:
$$({\mathbb{G}_{m,\C}}\to {\AC})_! ((X\times \mathbb{G}_{m,\C})/\mu_n\stackrel{g}{\to} \mathbb{G}_{m,\C})_! \Q_{(X\times\mathbb{G}_{m,\C})/\mu_n}(0)\,,$$ 
where the action of
 $\mu_n$ on $X\times \mathbb{G}_{m,\C}$ is the product of the given action on $X$ and the multiplication on $\mathbb{G}_{m,\C}$. Map $g$ is given
 by $g(x,t)=t^n$.

One can check that $\Phi$ maps the motivic Milnor fiber used in \cite{KS} to the class of cohomology with compact support of the sheaf of vanishing cycles. Therefore  we have proved the following result.

\begin{prp}
One has an equality of generating series with coefficients in $ K_0(MMHS_{\AC})$:
$$\Phi(A^{mot}_V) =D(A_V^{\op{crit}})\,,$$
where  $D$ is the duality.

\end{prp}
The above Proposition ensures that the theories developed in \cite{KS} and in the present paper  give essentially the same result at the level of generating series.

\section{Some speculations}

There are plenty of  conjectures and speculations about Cohomological Hall algebra and related motivic DT-series. From that variety we have chosen just two.

\subsection{Categorification of the critical COHA}

Let us keep the set up of Section 7.4. The critical cohomology $H_c^{\bullet,crit}(X^\s,f)$ are related to a certain $2$-periodic triangulated category. 
Namely, for any Zariski closed subset $X^\s\subset f^{-1}(0)$ let us define the category of matrix factorizations supported on $X^\s$ in the following way:
$$MF_{X^\s}(f):=D^b_{X^\s}(Coh(f^{-1}(0)))/Perf_{X^\s}(f^{-1}(0))\,,$$
where the subscript $_{X^\s}$ denotes the category of bounded complexes of coherent sheaves on the closed subscheme $f^{-1}(0)$ (resp. of perfect complexes
 on $f^{-1}(0)$), with cohomology sheaves supported on the closed subset $X^\s$.
Then there is a Chern character homomorphism
$$ch: K_0(MF_{X^\s}(f))\to (H_c^{ev,crit}(X^\s,f))^{\vee}$$
where the r.h.s. in de Rham realization is (hypothetically) the  periodic cyclic homology of $MF_{X^\s}(f)$.
One has also  an equivariant version $MF_{X^\s,G}(f)$ of the above category and of the Chern character (here $G$ is an algebraic group acting on $X$ and preserving $X^\s$ and $f$). We expect that the multiplication on the critical COHA comes from a monoidal structure on the direct sum of categories
$$\bigoplus_{\gamma \in \Z_{\ge 0}^I}\bigoplus_{z\in \C}MF_{\M_{\gamma}^\s \cap \left(W_\gamma \right)^{-1}(z),\G_{\gamma}}(W_{\gamma}-z)\,.$$
The correspondences $\M_{\gamma_1,\gamma_2}$ should be upgraded to functors between different summands, and the multiplication in the critical COHA will be the induced morphism of periodic cycllic homology spaces.
The monoidal structure could be thought of as a categorification of the critical COHA (which is itself a categorification of DT-invariants).

\subsection{Invariants of $3$-dimensional manifolds}

As we already mentioned in this paper, the comparison with \cite{KS} suggests that there exists a critical COHA associated with any ind-constructible $3CY$ category endowed with orientation data. The potential in this case is a (partially) formal function. Any compact oriented $3$-dimensional $C^{\infty}$ manifold $X$ gives an example. Namely, let us consider the triangulated category $D^b_{constr}(X)$ of complexes of sheaves with locally constant cohomology. This category has a $t$-structure with the heart equivalent to the category of finite-dimensional complex representations of the fundamental group $\pi_1(X,x_0), x_0\in X$. For a given $n\ge 0$ the stack ${ Rep}_n(X)$ of  representations of dimension $n$ is an Artin stack of finite type over $\C$. Locally (in analytic topology) we can represent ${Rep}_n(X)$ as the set of critical points of the  Chern-Simons functional:
$$CS (A)=\int_X Tr\left({\frac{dA\cdot A}{2}}+{\frac{A^3}{3}}\right)\in \C/(2\pi i)^2\Z,$$
where $A\in \Omega^1(X)\otimes Mat(n,\C)$, modulo action of the gauge group.
It looks plausible that the corresponding $3CY$ category admits orientation data in the sense of \cite{KS}.
Therefore we obtain a topological invariant of $X$ given by the motivic DT-series in one variable. For $X=S^3$ the invariant coincides with the motivic DT-series for the quiver $A_1=Q_0$ endowed with the trivial potential (essentially it is the quantum dilogarithm, see Section 2.5). For $X=(S^1)^3$ it is given by  Proposition 7 from Section 5.6 (quantum MacMahon function).
One can also compute the invariant e.g. for $X=S^1\times S^2$, but in general the answer is not known to us.

Another interesting story is a relation of COHA with holomorphic Chern-Simons functional for $\overline{\partial}$-connections on $C^\infty$ complex vector bundles on a compact complex $3CY$ manifold $ X_\C$ endowed with a non-zero holomorphic $3$ form $\Omega^{3,0}_{X_\C}$. In this case $CS$ is defined modulo the abelian subgroup of $\C$ consisting (up to a torsion)  of integrals  of $\Omega^{3,0}_{X_\C}$ over integral cycles. 

 Also for both $C^\infty$ and holomorphic Chern-Simons functionals one can try to define a {\it rapid decay}  version of COHA. The 
latter is achieved by  taking into account Stokes data  (the same as gluing data in \cite{KaKoPa}) and counting  gradient lines of the real part of $\exp(i\phi) CS$ for various $\phi\in\R/2\pi\Z$, which connect different critical points of $CS$.
 In order to do this we have to use an appropriate infinite covering of the space of connections, in order to have a globally defined holomorphic  functional $CS$. This goes beyond the formalism of  $3CY$ categories, as the gradient lines are trajectories in the space of  non-flat connection in the case of a real oriented $3$-dimensional manifold, or
non-holomorphic $\overline{\partial}$-connections in the complex case. Geometrically these lines correspond to self-dual non-unitary connections on the 4-dimensional Riemannian manifold $X\times \R$ in the real case, or on the $Spin(7)$-manifold $X_\C\times \R$ in the complex case, with appropriate boundary conditions at infinity (compare with \cite{Witten}). The resulting structure  is in a sense an exponential mixed Hodge structure of infinite rank. We hope to discuss it in the future.

\vspace{3mm}

Addresses:

M.K.: IHES, 35 route de Chartres, F-91440, France

{maxim@ihes.fr}\\

Y.S.: Department of Mathematics, KSU, Manhattan, KS 66506, USA

{soibel@math.ksu.edu}

     \end{document}